\documentclass[10pt,a4paper, draft]{article}
\usepackage[utf8]{inputenc}
\usepackage{amsmath}
\usepackage{amsfonts}
\usepackage{amssymb,bbm,fullpage,subfigure}
\usepackage{enumerate}
\usepackage[left=2cm,right=2cm,top=2cm,bottom=2cm]{geometry}
\usepackage[colorlinks=true]{hyperref}
\usepackage{pgfplots}
\usepackage{tikz}
\usetikzlibrary{calc,patterns,angles,quotes}
\usepackage{amsthm}
\usepackage{thmtools}
\everymath{\displaystyle}
\usepackage[capitalise]{cleveref}
\newtheorem{theorem}{Theorem}[section]
\newtheorem{proposition}[theorem]{Proposition}
\newtheorem{corollary}[theorem]{Corollary}
\newtheorem{fact}[theorem]{Fact}
\newtheorem{lemma}[theorem]{Lemma}
\theoremstyle{definition}
\newtheorem{definition}{Definition}[section]

\newcommand{\qede}{\hspace*{\fill}$\Diamond$\medskip}
\declaretheorem[style=remark,qed=\qede,Refname={Remark,Remarks},sibling=theorem]{remark}
\declaretheorem[style=remark,qed=\qede,Refname={Example,Examples},sibling=theorem]{example}

\usepackage{framed}

\DeclareMathOperator{\cone}{cone}

\DeclareMathOperator{\dom}{dom}

\DeclareMathOperator{\Argmax}{Argmax}
\DeclareMathOperator{\argmin}{argmin}

\usepackage{todonotes}
\usepackage{cleveref}

\newcommand{\al}{\alpha}
\newcommand{\be}{\beta}

\newcommand{\la}{\lambda}
\newcommand{\de}{\delta}
\newcommand{\eps}{\varepsilon}
\newcommand{\bx}{\bar x}
\newcommand{\by}{\bar y}

\newcommand {\R} {\mathbb R}
\newcommand {\N} {\mathbb N}

\newcommand {\B} {\mathbb B}

\newcommand {\epi} {{\rm epi}\,}

\newcommand {\cl} {{\rm cl}\,}

\newcommand {\conv} {{\rm conv}\,}
\newcommand {\bd} {{\rm bd}\,}

\newcommand {\Int} {{\rm int}\,}

\newcommand{\menge}[2]{\big\{{#1} \mid {#2}\big\}}%for defining sets
\newcommand{\norm}[1]{\left\Vert#1\right\Vert}
\newcommand{\abs}[1]{\left\vert#1\right\vert}

\newcommand{\ang}[1]{\left\langle #1 \right\rangle}

\newcommand{\AND}{\quad\mbox{and}\quad}
\newcounter{mycount}

\makeatletter
\begin{document}
\title{Single-Projection %\todo{MKT: ``Procedure''? because we say SPP}
Procedure for Infinite Dimensional Convex Optimization Problems}
\author{Hoa T.\ Bui\footnote{ARC Training Centre for Transforming Maintenance through Data Science, and 
Curtin Centre for Optimisation and Decision Science, Curtin University, Australia. Email: hoa.bui@curtin.edu.au}
\and
Regina S.\ Burachik\footnote{STEM Mathematics, University of South Australia, Australia. Email: Regina.Burachik@unisa.edu.au}
\and
Evgeni A.\ Nurminski\footnote{Center for Research and Education in Mathematics (CREM), Far-Eastern Federal University, Vladivostok, Russia. Email: nurminskiy.ea@dvfu.ru}
\and
Matthew K.\ Tam\footnote{School of Mathematics and Statistics, University of Melbourne, Australia. Email: matthew.tam@unimelb.edu.au}}

\maketitle

%\todo[inline]{MKT: Please check AU vs US spelling for consistency. At the moment, there is both in the manuscript. Since we seems to be following US spelling, I have changed anything I have found to US (although it pains me to do! :p)}

\begin{abstract}
In this work, we consider a class of convex optimization problems in a real Hilbert space that can be solved by performing a single projection, i.e., by projecting an infeasible point onto the feasible set. Our results improve those established for the linear programming setting in Nurminski (2015) by considering problems that: (i) may have multiple solutions, (ii) do not satisfy strict complementary conditions, and (iii) possess non-linear convex constraints. As a by-product of our analysis, we provide a quantitative estimate on the required distance between the infeasible point and the feasible set in order for its projection to be a solution of the problem. Our analysis relies on a ``sharpness'' property of the constraint set; a new property we introduce here.
\end{abstract}



\section{Introduction}\label{intro}
%\textcolor{green}{EN: 02.08. Looking over the article I noticed that the term "projection" is used many times but never formally defined. Also it will be better to introduce the projection operator which makes statements shorter. Who volunteers :) ?}https://www.overleaf.com/project/5e3bbcae7896760001761cf3
Let $\mathcal H$ be a real Hilbert space
with inner product $\ang{\cdot,\cdot}$ and induced norm denoted
$\|x\| = \sqrt{\ang{x,x}}$. Consider a problem of the form
\begin{equation}
    \label{LP}\tag{P}
    \min_{x\in A}\ang{x^*,x},
\end{equation}
where $A\subseteq\mathcal{H}$ is a nonempty, closed and convex set and $\ang{x^*,\cdot}$ is a linear function ($x^*\neq 0$). Without loss of generality, we assume that  $\|x^*\|=1$. In \cite{nurminski2015single}, the author shows that if $\mathcal H$ is finite dimensional, $A$ is a polyhedron and strict complementarity conditions hold, then the linear programming problem \eqref{LP} has the property:
\begin{quote}
    \emph{For every $x^0\in \mathcal H$, there exists $\theta_0>0$ such that the projection onto $A$ of the point $x^0-\theta x^*$ (which is a steepest descent step from the initial guess $x^0$) solves \eqref{LP} for any $\theta \geq \theta_0$.}
\end{quote} We refer to the aforementioned procedure as the {\em single projection procedure} (SPP) for Problem~\eqref{LP}. 
{Since the publication of \cite{nurminski2015single}, the SPP for solving linear programming in finite dimensions has been considered in \cite{BehBellSan,bui2021}, as a consequence of the finite convergence property of the alternating projection method performed on a polyhedral set and a closed half space in the case when these sets are not intersecting.} 

%\textcolor{purple}{RSB 29 Sept: The previous sentence is not clear, especially the use of "by-product". Do you mean "consequence" of the finite convergence properties?}
%\textcolor{blue}{HB: I am not sure if it reads better now?}
% \todo[inline]{EN: The idea to single out this from the regular text looks fine. The only problem is that it is called ''property'' first, and ''procedure'' later. Better stick to ''procedure'' and reword a bit accordingly. \todo[inline]{MKT: This ``problem'' pre-dates my edits :) I'd propose for minor things like this to edit yourself.}}

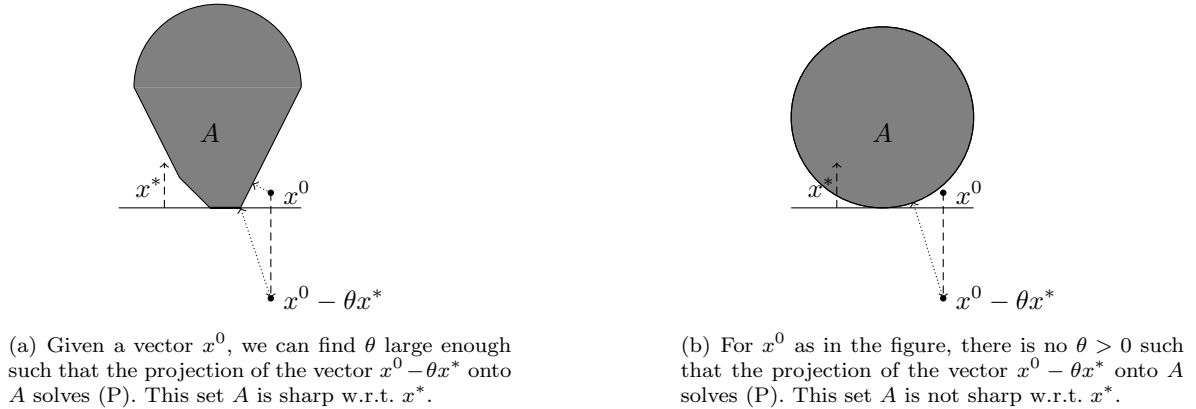
\begin{figure*}[h]
	\centering
		\subfigure[Given a vector $x^0$, we can find $\theta$ large enough such that the projection of the vector $x^0-\theta x^*$ onto $A$ solves (P). This set $A$ is sharp w.r.t.\ $x^*$. ]{\makebox[6.5cm][c]{
		\begin{tikzpicture}[scale= 0.4]
		\fill[gray,opacity=0.1] (0,0)--(1,0)--(2,2)--(3,4)--(-2.5,4)--(-1,1)--(0,0);
        \fill[gray,opacity=0.1] (-3,4)--(3,4) arc(0:180:2.75) --cycle;
        \draw[black,opacity=0.5] (3,4) arc(0:180:2.75);
        \draw[black,opacity=0.5] (0,0)--(1,0)--(2,2)--(3,4);
        \draw[black,thick,opacity=0.5] (0,0)--(1,0);
        \draw[opacity=0.5] (-2.5,4)--(-1,1)--(0,0);
        \draw[black] (-3,0) -- (3,0);
        \draw[->,black,densely dashed] (-1.5,0) -- (-1.5,1.5); 
        \filldraw (0,2.5) node {$A$};
        \filldraw (-2,0.7) node {$x^*$};
        \node[circle,draw=black, fill = black, inner sep=0pt,minimum size=2pt, label=right:{$x^0-\theta x^*$}] (b) at (2,-3) {};
        \draw[densely dotted,black,->](2,-3) -- (1,0);
         \node[circle,draw=black, fill = black, inner sep=0pt,minimum size=2pt, label=right:{$x^0$}] (b) at (2,0.5) {};
         \draw[densely dotted,black,->](2,0.5) -- (1.4,0.8);
        \draw[densely dashed,->](2,0.5) -- (2,-3);
		\end{tikzpicture}
	}
}
\quad \quad \quad \quad\quad\quad
\subfigure[For $x^0$ as in the figure, there is no $\theta >0$ such that the projection of the vector $x^0-\theta x^*$ onto $A$ solves (P). This set $A$ is not sharp w.r.t.\ $x^*$.]{\makebox[6.5cm][c]{
	\begin{tikzpicture}[scale= 0.4]
		\draw[fill = gray,opacity=0.1] (0,3) circle (3cm);
		\draw[black,opacity=0.5] (0,3) circle (3cm);
		\draw[black] (-3,0) -- (3,0);
		\node[circle,draw=black, fill = black, inner sep=0pt,minimum size=2pt, label=right:{$x^0-\theta x^*$}] (b) at (2,-3) {};
        \draw[densely dotted,black,->](2,-3) -- (1,0.2);
        \node[circle,draw=black, fill = black, inner sep=0pt,minimum size=2pt, label=right:{$x^0$}] (b) at (2,0.5) {};
        \draw[densely dashed,->](2,0.5) -- (2,-3);
         \draw[->,black,densely dashed] (-1.5,0) -- (-1.5,1.5); 
         \filldraw (0,2.5) node {$A$};
         \filldraw (-2,0.7) node {$x^*$};
	\end{tikzpicture}
}
}
	\caption{The single-projection procedure extends beyond linear programming. As shown in Figure \ref{Fig:F1}(a), the set $A$ does not need to be polyhedral for the procedure to work, but the sharpness property at $x^*$ is required  (see Definition~\ref{D0}).}\label{Fig:F1}
	% \todo[inline]{MKT: We have two different colour schemes in figures (grey vs blue). I think we should make them the same for consistency. I can do this, but which do we want? I think gray is better for grayscale printing.}
	% \todo[inline]{EN:Vote for gray}
\end{figure*}

%\textcolor{blue}{EN: Does in make sense to introduce the abbreviation SPP ? Read SiProPro :) ? Save sibirian %trees !}
%
%to the set $A$, where $x^0$ and $\theta\in \R$ are suitably chosen. 
The aims of the present paper can be summarized as follows: (i) extend the SPP to the infinite dimensional setting, (ii) show that the SPP remains valid in more general settings including non-linear convex constraints, non-unique solutions and/or in the absence of the strict complementary property, and (iii) provide quantitative estimates on the value of $\theta_0$ needed for the SPP to work. To this end, we show that the SPP is valid for problems satisfying a new property called {\em sharpness} (see Figure~\ref{Fig:F1}) which we introduce here. In the particular case of problem~\eqref{LP}, the sharpness property holds when there is a positive lower bound for the distance between the (unit) vector $-x^*$ and the normal cone $N_A(\cdot)$ at every non-optimal point of Problem~\eqref{LP} (see Definition~\ref{D0}).
In this context, Problem~\eqref{LP} can be solved by the SPP (see Theorem~\ref{T3} and Lemma~\ref{rem:trans}) whenever $\theta$ is {\em sufficiently large} and $x^0 $ is {\em sufficiently far from being optimal} in the sense that
$$\ang{x^*,x^0-\theta x^*} =\ang{x^*,x^0} - \theta <\inf_{x\in A}\ang{x^*,x}.$$
Furthermore, as polyhedral sets are sharp with respect to every unit vector (see Proposition~\ref{active}), our analysis shows that every solvable linear programming problem can be solved by the SPP. As a consequence, this work extends the main results of \cite[Theorem 1]{nurminski2015single} to Hilbert spaces. The setting of potentially convex objective functions is dealt with in Theorems~\ref{T6} and \ref{T4}.

As mentioned above, our analysis relies on a new notion of ``sharpness'' for sets, which this work also studies in its own right in Section~\ref{RP}. Namely, a set $A$ is {\em sharp with respect to $x^*$} if and only if 
$$
\inf_{\substack{x\in A\setminus F_A(x^*)}}d(x^*, N_A(x)) >0,
$$
where $F_A(x^*)$ denotes the face of $A$ defined by a point $x^*$ (see Definition~\ref{D1}). In order to contextualize sharpness in the broader literature, we explore the property from three
perspectives:
\begin{itemize}
    \item \emph{Sharpness of the epigraph:} When the set $A$ is the epigraph of a convex function, we provide characterizations of sharpness in terms of its subdifferential. This includes establishing that the epigraph of a function is sharp with respect to the vector $(0_{\mathcal{H}},-1)$ if and only if it satisfies a global {\em Kurdyka-\L{}ojasiewicz (KL)} property with exponent $0$ (Proposition~\ref{EXP1}). Moreover, we show that a set $A$ is sharp with respect to $x^*$ if and only if the function $\mathbbm{1}_A(\cdot)-\ang{x^*,\cdot}$ has the global KL property with exponent $0$, where $\mathbbm{1}_A(\cdot)$ is the indicator function of the set $A$. 
    \item \emph{Dual characterizations:} The analysis in~\cite{nurminski2015single} relies on the presence of strict complementary. This condition, which is often not satisfied, is known to imply uniqueness of solutions of the linear optimization problem, and that the interior of the normal cone of the feasible set at this optimal solution contains the vector $-x^*$. On the other hand, we show in Proposition~\ref{PP3} that sharpness holds under much weaker conditions, in particular we allow non-unique optimal solutions. Namely, we require only that the union of the normal vectors of $A$ at all optimal solutions must contain $-x^*$. Consequently, although the strict complementary condition can easily fail for general linear programming problems, our sharpness condition holds for every polyhedral set and for every nonzero vector $x^*$ (see Proposition~\ref{active}).
    \item \emph{Metric characterizations:} We show that sharpness with respect to $x^*$ is equivalent to a \emph{subtransversality} property between the set $A$ and its supporting hyperplane at $x^*$ (see Corollary~\ref{C2}). %As a consequence, an amenable cone (see \cite{Lourenco}) is sharp with respect to \todo{MKT: This sentence needs revision. What is the qualification here? all unit vectors?} an unit vector in the range of $N_A(\cdot)$ if the linear span of the face $F_A(x^*)$ is the supporting hyperplane $F$ (see Corollary~\ref{P6}).
As such, the sharpness property can be connected with well-known geometric properties in the literature.
\end{itemize}

The remainder of the present paper is organized as follows. In Section~\ref{Pre}, we provide essential results that will be used in subsequent sections. Section~\ref{RP} formally introduces the new notion of the sharpness property and its connections with some existing geometry properties of sets. Section~\ref{S2} contains our main results on SPP for solving Problem~\eqref{LP}, and its extension to general convex problems in which the objective function is not necessary linear. 
Finally, Section~\ref{sec:conclusion} lists some open questions and discussion.
%\todo{Reminder: Don't forget to add something about ``conclusions and open questions section'' if something is written there.}

\section{Preliminaries}\label{Pre}
We start this section by setting the theoretical framework and recalling the standard definitions for future use.
As stated in the introduction, $\mathcal H$ is a real Hilbert space with inner-product $\ang{\cdot,\cdot}$ and induced norm $\|\cdot\|$. Given a set $C\subset \mathcal{H}$, the distance from $C$ to $x$ is denoted $d_C(x):=\inf_{z\in C} \|x-z\|$.

We use the notation $\R_{\infty}:=\R\cup \{\infty\}$. Given $A\subset \mathcal H$, we denote by ${\rm int}A$, ${\rm cl}A$, ${\rm bd}A$, its topological interior, closure, and boundary, respectively.
Unless specifically mentioned, we consider the strong topology in $\mathcal H$.
We will denote by $\mathcal B:=\{u\in \mathcal H\::\: \|u\|\le 1\}$ the closed unit ball in $\mathcal H$, and by
$\mathcal S := \{u\in \mathcal B\::\: \|u\|=1\}$ the boundary of $\mathcal B$. Consequently, the open unit ball is $\mathcal B\!\setminus{\mathcal S}=\{ u\in \mathcal B\::\: \|u\|<1\}$. Therefore, the open ball of radius $r>0$ and center $x_0\in \mathcal H$ is $x_0+r\, (\mathcal B\!\setminus{\mathcal S})$, and the closed ball of radius $r>0$ and center $x_0\in \mathcal H$ is $x_0+r\, \mathcal B$. 
%Dual norm here

We will consider the product space $\mathcal{H}^2=\mathcal{H}\times \mathcal{H}$ with the max norm $\|\cdot\|_{\infty}$. Namely, given $(x,y)\in \mathcal{H}^2$, we consider $\|(x,y)\|_{\infty}:=\max\{\|x\|,\|y\|\}$. Given this norm in $\mathcal{H}^2$, it is well known that its {\em dual norm} (i.e., the norm in the dual space of $(\mathcal{H}^2,\|\cdot\|_{\infty})$ is the norm 
\begin{equation}\label{def:dual norm}
    \|(x,y)\|_{*}:=\|x\|+\|y\|.
\end{equation}
For the closed unit ball in $\mathcal{H}^2$ induced by the sum norm, we will use the notation $B^*_{\mathcal{H}^2}$. Namely,
\begin{equation}
    \label{eq:L1}
   B^*_{\mathcal{H}^2}:=\{(u,v)\in \mathcal{H}^2\::\: \|u\|+\|v\|\le 1\}. 
\end{equation}

% We will also use the notation $B_\al(x^*):=\{v\in \mathcal H\::\: \|x^* - v\|< \al\}$, and $\overline B_\al(x^*):=\{v\in \mathcal H\::\: \|x^* - v\|\le \al\}$.
\begin{definition}\label{def:projection onto a set}
Let $C\subset \mathcal H$ be a nonempty, closed and convex set. Let $x\in \mathcal{H}$. By 
\cite[Theorem 3.16]{bauschke2011convex}, there exists a unique element in $C$ that minimizes the distance from $C$ to $x$. 
% \textcolor{green}{EN: I am slightly worried about the norm $\max\{\vert x \vert, \vert y \vert\}$ which is NOT strictly convex,
% and therefore min distance problem may have not-unique solution}
% \textcolor{blue}{HB: This norm max is only for the proof of Theorem 2.1, not to solve the problem $max_{x\in C}f(x)$
% }
We denote this element by $P_C(x)$ (it is also called the best approximation to $x$ from $C$). Namely,
\[
d_C(x)=\inf_{z\in C} \|z-x\| = \| x-P_C(x)\|.
\]
When convenient, we may also use the notation $d(x,C):=d_C(x)$.
\end{definition}

Let $f\colon \mathcal H\to \R_{\infty}$. The set
$\dom f:=\{x\in \mathcal H\::\: f(x)<\infty\}$ is the \emph{domain} (or \emph{effective domain}) of $f$, and the function $f^*\colon \mathcal H\to\R_{\infty}$ defined by
\[
f^*(x^*):=\sup_{x\in \mathcal H}\{\ang{x,x^*}-f(x)\},
\]
is the \emph{Fenchel conjugate} of $f$ at $x^*\in \mathcal H$.  The \emph{epigraph} of $f$ is $\epi f := \menge{(x,r)\in
\mathcal H\times\R}{f(x)\leq r}$.
We say that $f\colon \mathcal H\to \R_{\infty}$ is proper if $\dom f\neq\varnothing$.
A function $f$ is said to be (strongly) {\em lower semicontinuous (lsc)} when its epigraph is (strongly) closed. In the latter situation, we say that $f$ is {\em closed}. If $f$ is convex with a closed epigraph, then $f$ is also {\em weakly} lsc (i.e., $\epi f$ is closed in the weak topology).  Given a convex function $f$, recall that the \emph{subdifferential} of
$f$ is the point-to-set mapping $\partial f\colon \mathcal H\rightrightarrows \mathcal H$ defined by
 \begin{equation}\label{pre:subdiff}
     \partial f(x):= \begin{cases}\{x^*\in
\mathcal H \::\: (\forall y\in \mathcal H)\; \ang{y-x,x^*} + f(x)\leq
f(y)\}\,&\text{if}\, x\in \dom f;\\
\varnothing\,&\text{if } x \notin \dom f.
\end{cases} 
\end{equation} 
Note that for points at the boundary of $\dom f$, the subdifferential may or may not exist (i.e., the set in the first line in \eqref{pre:subdiff} may be empty).
% \textcolor{green}{We dropped the interior condition in the 1st line. And strictly speaking for $x$ at the boundary of $X$ the subdifferential may or may not exists, see $-\sqrt{x}$ for $x \geq 0$}

\begin{definition}\label{def:pts}
Given a point-to-set map $T:\mathcal{H}\rightrightarrows \mathcal{H}$, we consider the following sets.
\begin{itemize}
    \item[(a)] The {\em domain of $T$} is the set
    $$D(T):=\{x\in \mathcal{H}\::\: T(x)\not=\emptyset\}.$$
    \item[(b)] The {\em range of $T$} is the set
    $$\operatorname{R}(T):=\{v\in \mathcal{H}\::\: v\in T(x) \hbox{ for some }x\in\mathcal{H}\}.$$ 
    \item[(c)] The {\em graph of $T$} is the set
   $$ G(T):=\{(x,v)\in \mathcal{H}\times \mathcal{H}\::\: v\in T(x)\}.$$
\end{itemize}
\end{definition}

For a fixed nonzero vector $u\in \mathcal H$,  we will denote by
\(
\R_+(u):=\{tu\::\: t\ge 0\}={\rm cone}[u],
\)
the {\em cone} (and also called the {\em ray}) {\em generated by $u$}.   Given $J$ a nonempty set, and a collection of elements $(u_i)_{i\in J}\subset \mathcal H$ indexed by $J$, we denote by
${\rm cone}[u_i, i\in J]$ the convex cone generated by the collection. If $J=\{1,\ldots,l\}$ is finite, then 
\begin{equation}\label{eq:cones1}
    {\rm cone}[u_1,\ldots,u_l]=\sum_{i=1}^l {\rm cone}[u_i] = \sum_{i=1}^l \R_+(u_i),
\end{equation}
where we are using the fact that
${\rm cone}[A \cup B] = {\rm cone}[A] + {\rm cone}[B]$. Note that in these definitions we are using the notation ${\rm cone}[C]$ for the convex cone generated by a set $C$.
% \textcolor{blue}{3 Jun RSB: In this equality we are using the \underline{{\bf convex}} cone instead of cone, because  ${\rm cone}[A \cup B]={\rm cone}[A] \cup {\rm cone}[B]  $, so maybe we should clarify our notation for cone.}
% Given a set $D\subset H$, the {\em dual cone of $D$} is given by $D^+ := \{ v \in H \::\: 
% \ang{v,x} \ge 0,\, \forall\, x \in D\}$.

\begin{definition}\label{def:Normal}
Given a subset $A\subset \mathcal H$ and the point $x\in \mathcal H$, the point-to-set map $N_A:\mathcal H \rightrightarrows \mathcal H$ defined by
\[
N_A(x):= \begin{cases}
           \{x^*\in \mathcal H\::\:(\forall a\in A)\; \ang{a-x,x^*}\leq 0\} &\text{if~}x\in A; \\
           \emptyset & \text{otherwise,}
          \end{cases}
\]
is called the {\em normal cone of a set $A$ at the point $x$}.
%Adding the Normal cone of epigraph

\end{definition}

We will make use of the following result, which is \cite[Proposition 16.35]{bauschke2011convex}.

\begin{fact}\label{Normal cone of epigraph}
Let $f\colon \mathcal H\to \R_{\infty}$ be a convex, proper and lsc function. Denote by $D:=\dom f$ and by $E:=\epi f$. Fix $x\in D$. Then
\[
N_E((x,f(x))=\R_{++}(\partial f(x)\times \{-1\}) \bigcup (N_D(x)\times \{0\}),
\]
where $\R_{++}V:=\{t\,v\::\: t>0, v\in V\}$, denotes the positive cone generated by $V$. Equivalently, 
\begin{equation}
    N_E((x,f(x))=\{ (u,\eta)\in \mathcal H\times \R_- \::\: \eta<0 \hbox{ and } u/(-\eta) \in \partial f(x), \hbox{ or }\,
\eta=0 \hbox{ and } u \in N_D(x)\}.
\end{equation}
By \cite[Proposition 16.17(i)]{bauschke2011convex}, for every $x\in {\rm bd}D$, we can either have $\partial f(x)=\emptyset$ or $\partial f(x)$ unbounded. In the latter case, i.e., when $x\in {\rm bd}D$ and $\partial f(x)\neq\emptyset$, we have
\[
\partial f(x) +N_D(x) \subset \partial f(x).
\]
\end{fact}

Recall that the indicator function of a subset $A\subset \mathcal H$ is the function $\mathbbm{1}_A:\mathcal H\to \R_\infty$ such that $\mathbbm{1}_A(x)=0$ when $x\in A$ and $\mathbbm{1}_A(x)=+\infty$ otherwise. Hence, $\dom \mathbbm{1}_A=A$. If $A$ is a nonempty, closed and convex set, $\mathbbm{1}_A$ is a proper, lsc and convex function. In this situation, \eqref{pre:subdiff} yields
\begin{equation}\label{eq:F1}
   \partial (\mathbbm{1}_A (x))=N_A(x) \hbox{ for all } x \in \mathcal{H},
\end{equation}
and note that $\partial (\mathbbm{1}_A (x))=\emptyset$ otherwise. In this case, we can use the maximal monotonicity of the map $\partial (\mathbbm{1}_A)(\cdot)$ to deduce that the graph of $N_A$, given by the set
\[ 
G(N_A):=\{ (x,x^*)\in  A\times \mathcal H\::\: x^* \in N_A(x) \},
\]
is closed w.r.t.\ the strong-weak topology (i.e., w.r.t.\ the strong topology in the first coordinate and w.r.t. the weak topology in the second coordinate). We call this type of closedness {\em demi-closedness} (note that $G(N_A)$ is also closed when considering the weak topology in the first coordinate and the strong topology in the second one).
%\textcolor{red}{EN: SO, REGINA, I STOPPED HERE :) 8/14/22.}

\begin{remark}\label{rem:norm}
Fix $z_0\in \mathcal H$. Consider the function $\varphi_{z_0}:\mathcal H\to \mathbb{R}$ defined as 
$\varphi_{z_0}(x):=\|x-z_0\|$ where $\|\cdot\|$ is any given norm in $\mathcal H$.
In this situation, the Fenchel-conjugate of $\varphi_{z_0}$ is given by $\varphi_{z_0}^* = \mathbbm{1}_{\mathcal B_*} +\ang{z_0,\cdot}$ where $\mathcal B_*$ is the dual unit ball, i.e., the unit ball with respect to the norm which is dual to the given norm $\|\cdot\|$.  Moreover, by \cite[Example 16.62]{bauschke2011convex} (see also \cite[Corollary 2.4.16]{Za}):
\begin{equation}
    \label{E1}\partial \varphi_{z_0}(x)=\left\{
    \begin{array}{lr}
        {(x-z_0)}/\norm{x-z_0} & \hbox{ if } x\not=z_0, \\
        \mathcal B_* & \hbox{ if } x=z_0.
    \end{array}
    \right.
\end{equation}
Consequently, $\varphi_{z_0}$ is smooth (i.e., the set $\partial \varphi_{z_0}(x)$ is a singleton) at every $x\not= z_0$ see \cite[Proposition~17.32]{bauschke2011convex}. 
\end{remark}
%\textcolor{green}{$\nabla \varphi(x)= {x}/\norm{x}$ will look better (EN)}

\begin{remark}
\label{rem:max}
For future use, we recall here a fact involving the subdifferential of a maximum of two norms in the product space. Fix $\hat z=(\hat z_1,\hat z_2)\in \mathcal H\times \mathcal H$.
Consider the function $\Theta_{\hat z}: \mathcal H\times \mathcal H\to \mathbb{R}$ defined as
\[
\Theta_{\hat z}(x,y):=\max\{\|x-\hat z_1\|,\|y-\hat z_2\|\}=\max\{\varphi_{\hat z_1}(x),\varphi_{\hat z_2}(y)\},
\]
where we are using the notation introduced in Remark \ref{rem:norm} in the second equality. Note that
\begin{equation}\label{eq:DM}
   \Theta_{\hat z}(\hat z_1,\hat z_2)=0=\varphi_{\hat z_1}(\hat z_1)=\varphi_{\hat z_2}(\hat z_2).
\end{equation}
Denote by $\partial_1$ and $\partial_2$ the partial subdifferentials w.r.t. the first and second variables, respectively. Define also $\Theta^1(x,y):= \varphi_{\hat z_1}(x)$ for every $x,y\in \mathcal{H}$ and $\Theta^2(x,y):= \varphi_{\hat z_2}(y)$ for every $x,y\in \mathcal{H}$. So that $\Theta_{\hat z}(x,y)=\max\{\Theta^1(x,y), \Theta^2(x,y)\}$. Using Remark \ref{rem:norm} and the fact that $\Theta^2$ does not depend on the first variable,
we have that
\[
\partial \Theta^1(\hat z_1,\hat z_2)= \left( \partial_1 {\Theta^1} (\hat z_1,\hat z_2), \partial_2 {\Theta^1} (\hat z_1,\hat z_2)\right)=\left( \partial_1 {\Theta^1} (\hat z_1,\hat z_2), 0\right)=\left( \partial\varphi_{\hat z_1}(\hat z_1), 0\right) = B_*\times \{0\},
\]
where we used \eqref{E1} in the last equality. Similarly, we obtain
\[
\partial \Theta^2(\hat z_1,\hat z_2)= \left( \partial_1 {\Theta^2} (\hat z_1,\hat z_2), \partial_2 {\Theta^2} (\hat z_1,\hat z_2)\right)= \{0\}\times B_*,
\]
The classical Theorem of Dubovitskii and Milyutin (see \cite[Theorem 18.5]{bauschke2011convex}), which computes the subdifferential of a maximum of functions, together with \eqref{eq:DM}, yield
\[
\begin{array}{rcl}
   \partial \Theta_{\hat z}(\hat z_1,\hat z_2)  &  = &
   \overline{\rm conv}\{B_*\times \{0\}, \{0\}\times B_*\}=\{(u,v)\in  \mathcal H\times \mathcal H\::\: \|u\|+\|v\|\le 1\}=B^*_{\mathcal{H}^2}.
\end{array}
\]
(Here, $\conv(A)$ denotes the \emph{convex hull} of the set $A$, which is the smallest convex set that contains $A$; whereas $\overline{\rm conv}(A)$ is the closure of its convex hull.)
The second equality can be easily checked, while in the last equality we use the notation introduced in \eqref{eq:L1}.
\end{remark}

\begin{definition}\label{def:Adjoint}
% {\color{red}[MKT: IMO, the style of definition doesn't fit with the rest. I can understand why we have a special definition for polyhedral sets, but having the adjoint in it as well makes it clumsey for my taste.]}
Fix $m\in \mathbb{N}^*$ and $b\in \R^m$.  Let $\mathbb{A}:\mathcal H\to \R^m$ be a {bounded} linear operator. Recall that the {\em adjoint operator of $\mathbb{A}$} is the (bounded) linear map $\mathbb{A}^*:\R^m\to \mathcal H$ defined by the equality
\begin{equation}\label{def:A*}
 \ang{\mathbb{A}x,u}=\ang{x,\mathbb{A}^*u}\quad\forall (x,u)\in \mathcal H\times\R^m.
\end{equation}
\end{definition}

Note that, in \eqref{def:A*}, we are using the same notation for the inner-products in $\mathcal H$ and in $\R^m$. From \cite[Remark 16]{Brezis}, we have that, when $\mathbb{A}$ is bounded, $\mathbb{A}^*$ is also bounded and both maps have the same norm.

\begin{definition}\label{def:Polyhedral set}
Take $\mathbb{A}$ and $b$ as in Definition \ref{def:Adjoint}. Denote by ${\cal B}_m:=\{e_1,\ldots, e_m\}$ the canonical basis in $\R^m$. A {\em polyhedron or polyhedral set} induced by a linear map $\mathbb{A}$ and a vector $b$ is the set 	
\begin{equation}\label{def:polyhedral set eq}
    C(\mathbb{A},b):=\{x\in \mathcal H: \mathbb{A} x\le b\}:=\{x\in \mathcal H: \ang{\mathbb{A}x,e_j}\le b_j \hbox{ for all }j=1,\ldots,m\}.
\end{equation}
 Hence, a polyhedral set is a finite intersection of
%\textcolor{green}{$b$-level ? (EN)}
level sets of linear maps.	
Given $x\in \mathcal H$, the set $I(x):=\left\{j\in \{1,\ldots,m\}\::\: \ang{\mathbb{A}x,e_j}=b_j\right\}$ identifies the inequality constraints that are active at $x$. We say that a function $f:\mathcal{H}\to \R_{\infty}$ is {\em polyhedral} when its epigraph is a polyhedral subset of $\mathcal{H}\times \R$. 
% 	Note that $x\in {\rm bd\,}C(\mathbb{A},b)$ if and only if $I(x)\not=\emptyset$. 
% 	{\color{red} [MKT: I think the last claim regarding the boundary is not true. For example, if $\mathbb{A}=0$ and $b=0$, then $I(x)=\{1,\dots,m\}$ but $C(\mathbb{A},b)=H$.]   {\color{magenta} In this case ${\rm bd\,}C(\mathbb{A},b)=\emptyset$ so we need to add the assumption that the polyhedral has nonempty boundary.}}
\end{definition}

The normal cone to a polyhedral set will have an important role in our analysis, so we recall \cite[Corollary 4.1]{ng-song-2003}, valid in a locally convex space. Polyhedral sets in these general spaces are defined in a similar manner as in Definition \ref{def:Polyhedral set}, where $\mathcal H$ is replaced by a locally convex space denoted by $X$.

\begin{proposition}(\cite[Corollary 4.1]{ng-song-2003}) \label{fact:Normal Cone}
Let $X$ be a locally convex space, and let $C_1,\ldots, C_m$ be polyhedral
subsets of $X$ with $C:=\cap_{i=1}^m C_i\neq\emptyset$. Then, for all $x\in C$, we have
\[
	N_{C}(x)=\sum_{i=1}^m N_{C_i}(x).
\]
% {\color{red}[MKT: Should we remind the reader of the definition of a polyhedral set in $X$? Since the standard definition is only generally done in finite dimensions.]}
 \end{proposition}

We recall Ekeland's Variational Principle, which holds in the setting of metric spaces.

\begin{lemma}[{Ekeland's Variational Principle~\cite[Theorem~1.1]{Eke74}}]
\label{EVP}
Let $X$ be a complete metric space, $\psi: X\to\R_{\infty}$ be lsc, $\bar w\in X$ and $\varepsilon>0.$
If
$$
\psi(\bar{w})<\inf_{w \in X}\psi(x)+\varepsilon,
$$
then, for any $\la>0$, there exists $\hat{w}_{\lambda}\in X$ such that
%\textcolor{red}{EN: shell we better denote it as $\hat{w}_\lambda$ ?}
\begin{enumerate}[(i)]
\item $d(\hat{w}_\lambda,\bar w)<\lambda$;
\item $\psi(\hat{w}_\lambda)\le \psi(\bar w)$;
\item
$\psi(w)+(\varepsilon/\lambda)d(w,\hat{w}_\lambda)>\psi(\hat{w}_\lambda)$ for all $w\in X\setminus\{\hat{w}_\lambda\}.$
\end{enumerate}
\end{lemma}
% \textcolor{red}{there is smth strange with the last statement: by setting $x = \hat x$ we obtain contradicton (EN)}
%{\color{red}[MKT: I'm a bit confused about something with Lemma~\ref{L2}. If $u=v$ with $\|u\|=1$, then $\langle u,v\rangle =1 >0$ but $d(u,\R_+(v))^2=0$, but the lemma says it should be $1$?]}
%{\color{cyan}[HB: I haven't used this lemma anywhere. But I will revise it.]}

The next simple fact will be used in later sections.
\begin{fact}
\label{L2}
Let  $u,v\in \mathcal H$ be such that $\norm{u}=\norm{v}=1$. Then \(
d(u,{\rm cone}[v])^2= 1 - \max(0, \ang{u,v})^2.
\)
\end{fact}
\begin{proof}
Denote $ t_* = \ang{u,v}$. Then,
\begin{align*}
    d(u,{\rm cone}[v])^2 = \min_{t \geq 0} \|u - t v\|^2
    = \min_{t \geq 0} (1 + t^2 - 2 t_* t) = 1 - \max_{t \geq 0} ( -t^2 + 2 t_* t) =
    1 - \max (0, t_*)^2.
\end{align*}
\end{proof}
% \textcolor{blue}{I made it much shorter now.}
% \textcolor{green}{EN:18/6 Sure. But what about this:}
% \textcolor{green}{I found it easier to grasp right away :)}
%Denote by $\bd C$ the boundary of a set w.r.t. the strong topology.

\section{Sharp sets and their characterizations}\label{RP}
Figure~\ref{Fig:F1} illustrates that, for the SPP to work, the set $A$ must possess certain geometric properties. In this section, we formally define a property that allows problem \eqref{LP} to be solved by the SPP. Our focus is on a geometric property associated with the presence of ``sharp corners" of a convex subset of a Hilbert space. 
%\textcolor{magenta}{RSB: my changes in magenta below}
%{\color{green} EN: need smth to say here}
\subsection{Definition and examples}
%\vspace{5mm}
%\noindent In this section, we define a 

\begin{definition}
\label{D1}
We define
 the {\em face of a convex set $A\subset\mathcal{H}$ with respect to a nonzero vector $x^*$} as
$$F_A(x^*):=\{x\in A: \ang{x^*,x}=\sup_{A}\ang{x^*,\cdot}\}=\Argmax_A\ang{x^*,\cdot} .$$
A closed convex subset $F^0\subset A$ is said to be an {\em exposed face of $A$} if there is $x^*\in \mathcal H$ such that 
$
F^0=F_A(x^*).
$
\end{definition}
% We denote $\widehat N_A(x)$ the set of normal vectors in $N_A(x)$ that have norm $1$, i.e., $$\widehat N_A(x):=N_A(x)\cap S.$$
%{\color{red}[MKT: Do we need a special notation for this or can we just write the intersection when we use it? Because the ``hat'' notation creates confusion with the Frechet normal cone (or some other normal cone) which commonly use this notation.]}
The next definition will be crucial in our analysis.
\begin{definition}
\label{D0}
%{\color{red}[MKT: Do we also want to require that the set is convex? Since we use the (convex) normal cone?]}
% Given a closed convex set $A\subset H$, a vector $x^*\in S$ and a positive number $\al >0$, we say that the set $A$ is {\em $\al$-sharp with respect to $x^*$} if, for all $x\in A$ such that $x^*\notin N_A(x)$, we have
Let $A\subset \mathcal H$ be closed and convex, $x^*\in \mathcal S$ and  $\al >0$. We say that $A$ is {\em $\al$-sharp with respect to $x^*$} if, for all $x\in A$ such that $x^*\notin N_A(x)$, we have
\begin{equation}
    \label{DC}
   d\left(x^*, N_A(x)\right)\ge \al.
\end{equation}
%
%{\color{green} EN: is then the half-space $\{ x_2 >= 0 \}$ in $E^2$ sharp wrt $x^* = (o, -1)$ ? }
%
%{\color{blue} HB: yes.}
%
The {\em modulus of sharpness of $A$ w.r.t.\ $x^*$}, denoted $\text{sr}[A,x^*]$, is defined as 
\begin{equation}
    \label{D0.1}
    \text{sr}[A,x^*]:= \inf_{\substack{x\in A,\\ x^*\notin N_A(x)}}d(x^*,N_A(x)).
\end{equation}
We say that $A$ is {\em sharp w.r.t.\ vector $x^*$} if $\text{sr}[A,x^*] >0$.
%{\color{red}[MKT: I think the face definition should not be in the definition of our main property.]}
%{\color{red}[MKT: Should we introduce the notation $\text{sr}[A,x^*]$ as something like the ``modulus of sharpness of $A$ at $x^*$''?]}
\end{definition}

%{\color{magenta}RSB 1 Mar: What happens when the set $T_A(x^*):=\{ x\in A \::\;x^*\notin N_A(x) \}=\emptyset$? This may happen is $A$ has no interior. According to our convention, the infimum in this case is $+\infty$. So $A$ is sharp in that case?}

\begin{remark}\label{rem:sharp 1}
From the definition, $A$ is $\al$-sharp w.r.t\ to the vector $x^*$ if and only if $\text{sr}[A,x^*] \ge \al$. 
The definition of sharpness involves taking an infimum over the set $\{ x\in A \::\;x^*\notin N_A(x) \}$. When the latter set is empty, we have that $x^*\in N_A(x)$ for every $x\in A$. This implies that $\ang{x^*,x}=\sup_{A}\ang{x^*,\cdot}$ for all $x \in A$. Hence, $A= F_A(x^*)$. In this case, $A$ is trivially sharp w.r.t. $x^*$ by vacuity, and since we are taking the infimum over the empty set in \eqref{D0.1}, we deduce that $\text{sr}[A,x^*]=+\infty$. We also note that, since $\|x^*\|=1$ and $0\in N_A(x^*)$, we always have that $\alpha\le 1$. 
\end{remark}

% \textcolor{red}{ Do we really need that (above) ? This is a little bit confusing as if $A$ has nonempty interior then $x^* = 0$ which contradicts the earliest assumption that $\|x^*\| = 1$}
% \textcolor{blue}{ This comment is here to explain what happens when the infimum in the def of sharpness is taken over the empty set. If $x\in {\rm int}A$ then $N_A(x)=\{0\}$ we will never have $x^*\in S$. So this choice of $x^*$ never happens. I am explaining this a bit better in the remark. Let me know if it is clear now.}

%Then, $x^*\in N_A(x)$ if and only if $x\in A\cap F$. {\color{magenta} the equivalence is different for me. See the following remark.}

%%% MKT: The following sentence doesn't really fit with the fact, so I have changed it.
%Faces are connected with optimal solutions of LP problems, as well as with the sharpness property. 
The following result relates faces of sets with the sharpness property in Definition~\ref{D0}.

\begin{fact}\label{fact 1}
%{\color{red}[MKT: $A$ closed and convex?]}
Let $A$ be a closed convex set. Fix $x\in A$ and $x^*\in \mathcal S$. The following statements are equivalent.
\begin{itemize}
    \item[(i)] $x^*\in N_A(x)$.
    \item[(ii)] $x\in F_A(x^*)$.
    \item[(iii)] $x\in \Argmax_A{\ang{x^*,\cdot}}$.
\end{itemize}
% \textcolor{red}{ People are writing $\Argmax$ to denote the whole set of max-solutions and $\argmax$ to denote some solution. Shell we stick to this ?} \textcolor{blue}{Done!}
% Hence, we also have the following equivalent statements.
% \begin{itemize}
%     \item[(i)] $\text{sr}[A,x^*]\ge \al$.
%     \item[(ii)] $d(x^*,N_A(x)) \ge \al,\quad \forall x\in A\setminus F_A(x^*)$.
% \end{itemize}
Consequently, for any $\alpha\ge 0$, we have
$$
\text{sr}[A,x^*]\ge \al \Longleftrightarrow d(x^*,N_A(x)) \ge \al, \forall x\in A\setminus F_A(x^*).
$$
\end{fact}

\begin{proof} 
Using the definitions and the notation of Definition \ref{D0}, we can write
\[
x^*\in N_A(x) \Longleftrightarrow \ang{x^*,y-x}\le 0,\, \forall\, y\in A  \Longleftrightarrow \ang{x^*,x}= \sup_{y\in A}\ang{x^*,y} \Longleftrightarrow x\in F_A(x^*),
\]
which proves the equivalence between (i) and (ii) in the first statement. The equivalence between (ii) and (iii) in the first statement follows from the fact that $F_A(x^*)=\Argmax{\ang{x^*,\cdot}}$. To complete the proof, note that the equivalence between (i) and (ii) implies
\[
x\in A\setminus F_A(x^*) \Longleftrightarrow x\in A \hbox{ and } x^*\notin N_A(x).
\]
Therefore,
\[
\begin{array}{rcl}
  \text{sr}[A,x^*]\ge \al   & \Longleftrightarrow  & d(x^*,N_A(x)) \ge \al,\quad \forall x\in A, x^*\notin N_A(x),  \\
      &&\\
      &\Longleftrightarrow& d(x^*,N_A(x)) \ge \al,\quad \forall x\in A\setminus F(x^*),\\
\end{array}
\]
establishing the last claim.
\end{proof}

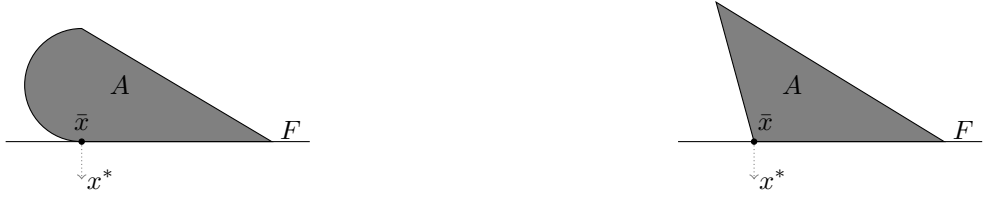
\begin{figure}[h]
	\centering
	\subfigure[The set $A$ is not sharp w.r.t.\ the vector $x^*$: there exist a sequence  $(x_n,x^*_n)\subset \bd A\times \mathcal H$ with $x_n^*\in N_A(x_n)$, $x^*_n\neq x^*$ ($n\in \N$) such that $(x_n)$ converges to $\bx$, and $(x_n^*)$ converges to $x^*$.]{\makebox[6.5cm][c]{
		\begin{tikzpicture}[scale= 0.5]
        \fill[gray,opacity=0.1] (3,0)--(-2,3) arc(90:270:1.5) --cycle;
        \draw[black,opacity=0.5] (3,0)--(-2,3) arc(90:270:1.5) --cycle;
        \draw[black] (4,0)--(-4,0);
         \draw[->,densely dotted,gray] (-2,0)--(-2,-1);
        \filldraw (-2,0) circle (2pt);
        \draw (-1,1.5) node {$A$};
         \draw (-2,0.5) node {$\bar x$};
         \draw (-1.5,-1) node {$x^*$};
         \draw (3.5,0.3) node {$F$};
		\end{tikzpicture}
}
}
\quad \quad \quad \quad\quad\quad
	\subfigure[The set $A$ is sharp w.r.t.\ the vector $x^*$: there exists $r>0$ such that, for any $y^*\in B_r(x^*)$ and $x\in A$ with $y^*\in N_A(x)$, it must hold that $x^*\in N_A(x)$.]{\makebox[6.5cm][c]{
		\begin{tikzpicture}[scale= 0.5]
        \fill[gray,opacity=0.1] (3,0)--(-3,3.7)--(-2,0)--cycle;
        \draw[black,opacity=0.5] (3,0)--(-3,3.7)--(-2,0)--cycle;
        \draw[black] (4,0)--(-4,0);
         \draw[->,densely dotted,gray] (-2,0)--(-2,-1);
        \filldraw (-2,0) circle (2pt);
        \draw (-1,1.5) node {$A$};
         \draw (-1.7,0.5) node {$\bar x$};
         \draw (-1.5,-1) node {$x^*$};
        \draw (3.5,0.3) node {$F$};
		\end{tikzpicture}
	}
}
	\caption{An illustration of the sharpness property. In both figures, $F$ is a supporting hyperplane of $A$ at $\bx$. The set $A$ is not sharp w.r.t.\ the vector $x^*$ in Figure~\ref{F1}(a). In Figure~\ref{F1}(b), there exist no sequence $(x_n)$ outside the face $A\cap F$, such that $x^*\in \cl \bigcup_{n}N_A(x_n)$. See Remark \ref{rem:F1}.
	%{\color{red}MKT: I find the captions a bit confusing here. Why do we need to talk about faces etc when definition 3.1 only mentions normal cones?]}
	}
	\label{F1}
\end{figure}

\begin{remark}\label{rem:F1}
Figure~\ref{F1} illustrates two-dimensional examples of the sharpness property with respect to a given vector, and its connection with faces of the set $A$. 
%In the figure, we denote by $F$ the face $F_A(x^*)$. 
In Figure~\ref{F1}(a), the ``rounded" section of the boundary of $A$ approaches $\bx$ {\em smoothly}. We note, however, that sharpness can also fail for sets whose boundary is not rounded as in Figure~\ref{F1}. We illustrate this situation in the next example.
\end{remark}
\begin{example}\label{F2a}
Consider the set $A\subset \R^2$ as the epigraph of the convex function $f:\R\to \R_{\infty}$ defined as follows.
\[
f(x):=\left\{
\begin{array}{cl}
     +\infty&\hbox{ for }x<-1,  \\
     &\\
  -\frac{(3n^2+3n+1)x+(2n+1) }{n^2(n+1)^2}   & x\in [-1/n, -1/(n+1)], n\in \mathbb{N}, \\
  &\\
    0 &x\ge 0.
\end{array}
\right.
\]
The boundary of the set $A$ to the left of 0 is determined by a piece-wise linear function which at the points $a_n := -1/n$ attains the value $1/n^3$ and is linear between $a_n$ and $a_{n+1}$ ($n = 1,2, \dots$). This function is non-smooth because it has kinks at each $a_n$. It is convex because its slopes monotonically increase to zero. The latter fact also implies that $A$ is not sharp w.r.t $x^*:=(0,-1)$. Therefore, a set with \emph{nonsmooth} boundary need not satisfy the sharpness condition.
\end{example}
% \begin{figure}
%     \centering
%     \begin{tikzpicture}
	
% 	\draw[scale=5,thick,domain=-0.6:-1/2, smooth, variable=\x, blue] plot ({\x}, {-(7*(\x)+3)/4});
% 	\draw[scale=5,thick,domain=-1/2:-1/3, smooth, variable=\x, blue] plot ({\x}, {-(19*(\x)+5)/36});
% 	\draw[thick,scale=5, domain=-1/3:-1/4, smooth, variable=\x, blue] plot ({\x}, {-(37*(\x)+7)/144});
% 	\draw[thick,scale=5, domain=-1/4:-1/5, smooth, variable=\x, blue] plot ({\x}, {-(61*(\x)+9)/400});
% 	\draw[thick,scale=5,blue] (-1/5,0.008)--(0,0.008);
% 	 \draw (0,2) node {epi$(\abs{x})$};
% 	 \draw[->,red] (1,0)--(1,-0.5);
% 	 \draw (1.2,-0.8) node {$(0,-1)$};
%     %\addplot[domain=-3:3,thick,black,no marks,smooth,fill=blue!30] {x^2};
%     \end{tikzpicture}
%     \caption{Caption}
%     \label{fig:my_label}
% \end{figure}

\begin{remark}
When a closed and convex set $A$ is bounded, the lack of sharpness implies a continuity property of the restriction of the point-to-set map $N_A(\cdot)$ to the boundary of $A$. Indeed, if $\text{sr}[A,x^*]=0$, then the definition directly implies the existence of a sequence $(x_n,x^*_n)\subset G(N_A)$ with $x^*\not\in N_A(x_n)$ and \ $x_n^*\to x^*$ strongly. The boundedness of $A$ implies that there exists $\bx$ and a subsequence $(x_{n_k})$ of $(x_n)$ weakly converging to $\bx$. For simplicity, we still call this subsequence $(x_n)$. Since $G(N_A)$ is demi-closed,  $(\bx,x^*)\in G(N_A)$.  The latter fact, together with $x^*\not\in N_A(x_n)$, implies that $x_n\not=\bx$ for all $n$. Altogether, we have a sequence $(x_n,x^*_n)\subset G(N_A)$ such that
\begin{itemize}
    \item[(i)] $x_n^*\to x^*$ strongly with $x^*\not\in N_A(x_n)$,
    \item[(ii)] there exists $\bx\in A$ such that $x_n\to \bx$ weakly, with $x_n\not=\bx$ for all $n$.
\end{itemize}
In particular, in Figure~\ref{F1}(a), for any positive number $r>0$, there is $y^*\in B_r(x^*)$ and $y\in B_r(\bx)\setminus F$ such that $y^*\in N_A(y)$. Note here that $\bx$ is an extreme point of the set $A$ in both Figures~\ref{F1}(a), and (b). {However, the point $\bx$ in Figure~\ref{F1}(b) is an {\em exposed point} (i.e., there exists a hyperplane $H$ such that $A\cap H=\{\bx\}$), whereas the point $\bx$ in Figure~\ref{F1}(a) is not an exposed point (since the only supporting hyperplane of $A$ at $\bx$ is $F$).} In Figure~\ref{F1}(b) it is easy to note that, there is $\al \in (0,1)$ such that for every $y\in A$ with $x^*\notin N_A(y)$ we will have that $d(x^*, N_A(y))\ge \alpha.$
\end{remark}

\begin{remark}\label{rem:Sharp and bounded}
From Fact \ref{fact 1}(i)-(ii), we know that $x^*\notin {\rm R}(N_A)$ if and only if $F_A(x^*)=\emptyset$. We note that this situation cannot hold when $A$ is bounded. Indeed, in this case we can use \cite[Corollary 21.25]{bauschke2011convex} to obtain $\text{R}(N_A)= \mathcal H$. Namely, there is no $x^*$ verifying $x^*\notin {\rm R}(N_A)$. We also note that, in general, the set $\text{R}(N_A)$ is neither closed nor convex. However, it is well known (see, e.g.,  \cite[Corollary 21.25]{bauschke2011convex} or \cite[Theorem 4.4.9]{BurachikIusemBook}) that ${\rm cl }\text{ R}(N_A)$ is a closed and convex set.
\end{remark}
%\textcolor{blue}{HB: I found the book, but couldn't find Corollary 21.25}
% 
% \begin{enumerate}
%     \item[1.] The set $A$ is bounded, i.e., there is $M>0$ such that $\norm{x} <M$ with $x\in A$.
%     Recall that for a bounded set, by Bishop-Phelps theorem, the range of the normal cone operator is dense (w.r.t. strong topology) in $H$. So, there is a sequence $(x_n^*)\subset \text{Range}(N_A)$ such that $x_n^*\to x^*$. Because $x^*\notin \text{Range}(N_A)$, we have $x^*\notin N_A(x_n)$, and by density property,
% $$
% \lim_{n\to \infty} d(x^*,N_A(x_n))=0.
% $$
% In other words, when
% $A$ is bounded, $A$ is not sharp with respect to any vector $x^*\notin \text{Range}(N_A)$. 
% \item[2.] 

% The set $A$ is unbounded. For this case, the set $A$ might or might not be sharp with respect to $x^*\notin \text{Range}(N_A)$. 
% \end{enumerate}

For an unbounded set $A$, we next characterize for which $x^*\notin {\rm R}(N_A)$  we have that $A$ is sharp w.r.t.\ $x^*$.

\begin{proposition}\label{pro:range of NA}
Fix $x^*\in \mathcal S$ and assume that $A\subset \mathcal H$ is a closed and convex set such that $x^*\notin \text{\rm R}(N_A)$. Then $A$ is sharp w.r.t. $x^*$ if and only if {$x^*\notin \bd \text{\rm R}(N_A)$}. In this situation, $A$ must be unbounded.
%Consequently, if $A$ is bounded, then there is no $x^*$ s.t. $A$ is sharp w.r.t. $x^*$. 
\end{proposition}
%\textcolor{blue}{HB 3/3/2022: maybe it is better to put it as boundary instead of closure}
\begin{proof}
The fact that $A$ is unbounded follows from Remark \ref{rem:Sharp and bounded}.
 We show first that, if {$x^*\in {\rm bd }\text{ R}(N_A)$} then $A$ is not sharp w.r.t. $x^*$. Indeed, in this case, we can take a sequence $(x_n)\subset A$ such that $x^*\not=x^*_n\in N_A(x_n)$ with $x^*=\lim_n x^*_n$. The fact that $x^*\not=x^*_n$ holds because $x^*\notin {\rm R}(N_A)$. Hence,
\[
\inf_{x\in A,\, x^*\notin N_A(x)}d(x^*,N_A(x) )\le \lim_{n\to \infty} d(x^*,x^*_n)=0,
\]
so $A$ is not sharp w.r.t. $x^*$. Conversely, assume that $x^*\notin \bd \text{\rm R}(N_A)$. Since we also have that $x^*\notin \text{\rm R}(N_A)$, we deduce that $x^*\notin {\rm cl }\text{ R}(N_A)$. Then, there exists $\al>0$ such that
\[
\al=d(x^*,{\rm cl }\text{ R}(N_A))= \inf_{v\in \text{R}(N_A)} d(x^*,v) \le \inf_{x\in A,\, x^*\notin N_A(x)}d(x^*,N_A(x) ),
\]
so $A$ is sharp w.r.t. $x^*$. 
\end{proof}

We illustrate the last result in the next two examples. In the first example, we have  $x^*\notin  \text{R}(N_A)$ but $x^*\in {\rm bd }\text{ R}(N_A)$. While in the second example $x^*\notin {\rm cl }\text{ R}(N_A)$.

\begin{example}\label{EX:2.8}
% By \cite[Corollary 21.25]{bauschke2011convex}, if $A$ is unbounded, then ${\rm range }(N_A)\subsetneq H$, so there exists $x^*\notin {\rm range }(N_A)\cap S$.
If $A:=\epi{g}$ where $g:\mathbb{R}\to \mathbb{R}_{\infty}$ defined as $g(t)=1/t$ if $t>0$ and $g(t)=+\infty$ otherwise. Then it is easy to check that $x^*:=(0,-1)\notin \text{R}(N_A)$ but $x^*\in \bd \text{\rm R}(N_A)$. Hence, we are not in the conditions of Proposition \ref{pro:range of NA}. Let us check that $A$ is not sharp with respect to $x^*$. Indeed, for $x_n:=(n,1/n)$ we have $(-1/n^2,-1)\in N_A(x_n)$ and
\[
x^*=(0,-1)= \lim_{n\to \infty} (-1/n^2,-1),
\]
so $\lim_{n\to \infty}d(x^*,N_A(n,1/n))=0$ and hence $A$ is not sharp with respect to $x^*$. 
\end{example}

\begin{example}\label{exa:sharp1}
Take $A=\R^2_+$ the nonnegative orthant in $\R^2$ and $x^*=(1/\sqrt{2})(1,1)$. Then $x^*\not\in {\rm cl }\text{ R}(N_A)$ and we are in the conditions of Proposition \ref{pro:range of NA}. It can be verified that $A$ is $1$-sharp w.r.t. $x^*$. Indeed, if  $x\in {\rm int} A$ then $N_A(x)=\{0\}$ so  $d\left(x^*, N_A(x)\right)= d(x^*,0)=1$. If $x=t(1,0)$ for $t>0$,  then $N_A(x)=\{s(0,-1)\::\: s\ge 0\}$. Using Lemma \ref{L2} with $u:=x^*$ and $v:=(0,-1)$, we obtain $d\left(x^*, N_A(x)\right)=1$.  Identical argument shows that $d\left(x^*, N_A(x)\right)=1$ for $x=t(0,1)$ for $t>0$. In the latter case we use  Lemma \ref{L2} with $u:=x^*$ and $v:=(-1,0)$. Finally, if $x=(0,0)$ then $N_A(x)=\R^2_{-}$ so again we have $d\left(x^*, N_A(x)\right)= d\left(x^*, \R^2_{-}\right)=1$.
% If $A$ is the unit ball in $\R^2$ then it can be seen that there is no positive $\al$ that verifies condition \eqref{DC} for all $x\in A$ such that $x^*\not\in N_A(x)$.
%Further illustration of this property is provided in Figure \ref{fig0}.}
\end{example}

We show in the next result that when the set $A$ is bounded, the constant $\text{sr}[A,x^*]$ is exactly the infimum of the distance between $x^*$ and every vector $y^*$ such that the face $F_A(y^*)$ does not intersect with the face $F_A(x^*)$. This observation holds true even in the infinite dimensional setting.

%RSB29Sept

\begin{proposition}
\label{PP01} Suppose $A$ is a nonempty bounded closed convex set. 
Fix $x^*\in \mathcal S$ and $\al \in (0,1)$.
Then, $A$ is $\al$-sharp w.r.t. $x^*$ if and only if, for every $y^*\in \mathcal H\setminus \{0\}$ such that $F_A(y^*)\cap F_A(x^*)=\emptyset$, we have
\begin{equation}
    \label{PP01-1}
    \norm{x^*-y^*} \ge \al.
\end{equation}
Consequently,
\begin{equation}\label{2nd eq}
    \text{\rm sr}[A,x^*] = \inf_{\substack{F_A(y^*)\cap F_A(x^*)=\emptyset\\ y^*\neq 0}}\norm{x^*-y^*}.
\end{equation}
\end{proposition}

\begin{proof}
Note that because the set $A$ is nonempty and bounded, from \cite[Corollary 21.25]{bauschke2011convex}, the set $\Argmax_A\ang{y^*,\cdot}$ is nonempty for any $y^*\in \mathcal H\setminus\{0\}$, therefore from Fact~\ref{fact 1}, the face $F_A(y^*)$ is also nonempty.

First suppose that the set $A$ is $\al$-sharp w.r.t. $x^*\in \mathcal S$. If $F_A(x^*)\cap F_A(y^*)$ is nonempty for every $y^*\in \mathcal H
\setminus\{0\}$, then equation \eqref{PP01-1} in Proposition~\ref{PP01} holds true immediately. Otherwise, assume that there exists $y^*\in \mathcal H\setminus \{0\}$ such that $F_A(y^*)\cap F_A(x^*)=\emptyset$. By Fact~\ref{fact 1}, for any $y\in F_A(y^*)$,  we have that $y^*\in N_A(y)$. Because $F_A(y^*)\cap F_A(x^*)=\emptyset$, we must have $y\not\in F_A(x^*)$, and again Fact~\ref{fact 1} yields $x^*\notin N_A(y)$. Then, the $\al$-sharpness property implies
\begin{equation}
    \label{1rst eq}
    \norm{x^*-y^*}\ge d(x^*,N_A(y))\ge \inf_{y'\in A,\, x^*\notin N_A(y')} d(x^*,N_A(y'))= \text{\rm sr}[A,x^*]\ge \al,
\end{equation}
which shows \eqref{PP01-1}. 

%$$\norm{x^*-y^*}\ge d(x^*,N_A(y))\ge \inf_{y'\in A,\, x^*\notin N_A(y')} %d(x^*,N_A(y'))= \text{\rm sr}[A,x^*]\ge \al,$$
%The second statement follows directly by taking infimum in the expression above.

Conversely, suppose that for every $y^*\in \mathcal H\setminus\{0\}$ satisfying $F_A(y^*)\cap F_A(x^*)=\emptyset$, it holds that $\norm{x^*-y^*} \ge \al$. We now show that the set $A$ must be $\al$-sharp w.r.t. $x^*$. 
In fact, suppose to the contrary that the set $A$ is not $\al$-sharp w.r.t.\ $x^*$. 
Then, there is $y\in A$ such that $x^*\notin N_A(y)$ and 
$$
d(x^*,N_A(y)) <\al.
$$
Take $y^*\in N_A(y)$ such that $\norm{x^*-y^*}<\al$, and take $\eps \in (0,\al - \norm{x^*-y^*})$. Note that because $x^*\notin N_A(y)$ and $y^*\in N_A(y)$, we have 
$$\ang{x^*,y} < \sup_{z\in A}\ang{x^*,z},\AND\ang{y^*,y-x}\ge 0,\,\forall x\in A.$$ 
Then, for every $x\in F_A(x^*)$, we have $\ang{x^*,x} = \sup_{z\in A}\ang{x^*,z}$, and 
\begin{align*}
    \ang{y^*-\eps x^*,y-x} &= \ang{y^*,y-x}-\eps\ang{x^*,y-x}\\
    &= \ang{y^*,y-x}-\eps(\ang{x^*,y}-\sup_{z\in A}\ang{x^*,z}) >0.
\end{align*}
The above expression yields $\ang{y^*-\eps x^*,x} < \sup_{z\in A}\ang{y^*-\eps x^*,z}$, which implies $x\notin F_A(y^*-\eps x^*)$. Because $x\in F_A(x^*)$ is chosen arbitrarily, we have $F_A(x^*)\cap F_A(y^*-\eps x^*)=\emptyset$. Now we can use the hypothesis \eqref{PP01-1} for $y^*-\eps x^*$ instead of $y^*$, to obtain
$$
\norm{(y^*-\eps x^*)-x^*} \ge \al.
$$
On the other hand, the definition of $\eps$ yields
$$
\norm{(y^*-x^*)-\eps x^*}\le \norm{y^*-x^*}+\norm{\eps x^*} = \norm{y^*-x^*}+\eps <\al,
$$
in contradiction with assumption \eqref{PP01-1}. Therefore, we must have $A$ is $\al$-sharp w.r.t. $x^*$. 

To complete the proof, we need to establish \eqref{2nd eq}. From \eqref{1rst eq} we deduce by taking infimum that
\[
  \text{\rm sr}[A,x^*] \le \inf_{\substack{F_A(y^*)\cap F_A(x^*)=\emptyset\\ y^*\neq 0}}\norm{x^*-y^*}.
\]
To prove the opposite inequality, assume that the inequality above is strict, so there is $\gamma$ such that 
\[
  \text{\rm sr}[A,x^*] <\gamma< \inf_{\substack{F_A(\hat{y})\cap F_A(x^*)=\emptyset\\ \hat{y}\neq 0}}\norm{x^*-\hat{y}}.
\]
By definition of $\text{\rm sr}[A,x^*]$ this means that there exists $x\in A$ such that $x^*\not\in N_A(x)$ with $d(x^*, N_A(x))<\gamma$. Thus, we can find $v\in N_A(x)$ s.t. $\|x^*-v\|<\gamma$. We have two possibilities: either $F_A(v)\cap F_A(x^*)=\emptyset$ or $F_A(v)\cap F_A(x^*)\not=\emptyset$. If $F_A(v)\cap F_A(x^*)=\emptyset$ then using the rightmost expression above we deduce that $\|x^*-v\|>\gamma$, contradicting the definition of $v$. Hence, the only possibility is that $F_A(v)\cap F_A(x^*)\not=\emptyset$. In this case, we will again obtain a contradiction. Indeed, as in the earlier part of the proof, take $\varepsilon\in (0, \gamma - \|x^*-v\|)$. Following the same arguments as above we will arrive at 
$F_A(x^*)\cap F_A(v-\eps x^*)=\emptyset$. Now we can use the rightmost expression above and the definition of $\varepsilon$ to deduce that
\[
\gamma<\|x^*-(v-\varepsilon x^*)\|\le \|v-x^*\| + \varepsilon <\gamma,
\]
a contradiction. Hence, we must have that \eqref{2nd eq} holds. This completes the proof.
\end{proof}

\begin{remark}
If the set $A$ is unbounded, the statements in Proposition~\ref{PP01} need not hold. Indeed, if $x^*\notin \text{\rm R}(N_A)$ then Fact~\ref{fact 1} yields $F_A(x^*)=\emptyset$ so $F_A(x^*)\cap F_A(y^*)=\emptyset$ for all $y^*\in \mathcal H$. In this situation,  inequality \eqref{PP01-1} will not hold for every nonzero $y^*\in \mathcal H$. This situation is illustrated by Example~\ref{exa:sharp1}. In this example, $x^*=(1/\sqrt{2})(1,1)\not\in \text{\rm R}(N_A)$, $F_A(x^*)=\emptyset$, and $A$ is $\alpha$-sharp at $x^*$ with $\alpha=1$. In particular, $F_A(x^*)\cap F_A(y^*)=\emptyset$ if $y^*=x^*$, in which case \eqref{PP01-1} is clearly false.
% \textcolor{magenta}{MKT 27/06: I'm have trouble seeing how the explanation in this remark related to unbounded sets. For instance, why do we assume ``if $x^*\notin \text{\rm Range}(N_A)$'' in relation to unbounded sets? Perhaps it would clearer to give a concrete counter-example here, rather than just argue plausability?}\textcolor{blue}{I added Example~\ref{exa:sharp1}}
\end{remark}

% A

% \textcolor{red}{ Did we define ''angle'' between faces ? What is ''positive'' angle ?} \textcolor{blue}{RSB 2 June: I deleted the confusing sentence, which doesn't seem to be necessary.}

We show next that a polyhedron is sharp with respect to every vector $x^*\in \mathcal S$. 
 
%Recall that the constant $\text{sr}[A,x^*]$ is the lower bound of all angles between the supporting hyperplane with respect to $x^*$, and all other supporting hyperplanes of $A$ that do not share common points with the face defined by $x^*$. 
%{\color{magenta} RSB 1 Mar: Hoa can you please express the above property mathematically? Maybe we need to prove it so we can use it comfortably. If you write down the statement I can try to prove it.}

%Therefore, if the set $A$ is a polyhedron, then each face of $A$ makes a positive angle with another face of $A$. Making use of the fact that $A$ has just finite number of faces, the lower bound of these angles, which is $\text{sr}[A,x^*]$, must be positive. This is the idea of our next result. Namely, if $A$ is a polyhedron, then $A$ is sharp w.r.t every vector $x^*\in S$. 

% {\color{magenta}RSB 12/01:It would be good to introduce a better notation for $\R_{+}(v)$. Maybe ${\rm cone}[v]$?}
% \begin{theorem}
% \cite[Theorem 4.1]{van2021generalized}
% Let $P$ and $\Omega$ be nonempty convex subsets of a locally convex topological vector space $X$, where $P$ is a convex polyhedron. Assume that the following quasi-relative interior qualification condition
% $$
% P\cap \qri(\Omega) \neq \emptyset,
% $$
% is satisfied. 
% Then we have the normal cone intersection rule
% $$
% N_{P\cap \Omega}(x)= N_P(x)+N_\Omega(x),\quad \forall x\in P\cap \Omega.
% $$
% \end{theorem}

\begin{proposition}
\label{active}
%{\color{red}[MKT: Revise ``dual vector'' terminology?]}
	Let $A$ be a polyhedron defined by
	$A:=\{x: \mathbb{A} x\le b\}$,
	where $\mathbb{A}$ is a bounded linear operator $\mathbb{A}: \mathcal H \to \R^m$ ($m >0$), and $b\in \R^m$. Denote by ${\cal B}_m:=\{e_1,\ldots, e_m\}$ the canonical basis in $\R^m$. 
	For any vector $x^*\in \mathcal S$, let $J(x^*)$ be a collection of subsets of indexes defined as $J(x^*) :=\left\{J\subset \{1,\ldots,m\} \::\: x^*\notin \cone[\mathbb{A}^*\,e_i]_{i\in J}\right\}$, and
	\begin{equation}\label{al1}
	\al_0:=\min\left\{1,\inf_{J\in J(x^*)} d \left(x^*, \cone[\mathbb{A}^*\,e_i]_{i\in J}\right) \right\},
	\end{equation}
	  with the convention that $\inf\emptyset = +\infty$. Then, $\al_0 >0$ and $A$ is $\al_0$-sharp w.r.t. $x^*$.
\end{proposition}

\begin{proof}
%\todo[inline,caption={}]{MKT: I'm bit unsure about the first line. If we take $c=-e_1$, $M=-1$, $A=e_1$ and $b=0$, then
% $$ A=\{x:x_1=1\},\quad B=\{x:x_1\leq 0\}, $$
%which are disjoint, convex sets. Then for all $x\in A$, we have $y=P_B(x)=(0,x_2,\dots)$ which means that $d(y,A)=1=d(A,b)$ (ie was don't have $d(y,A) > d(A,B)$). 
%}
%We are going to show that $A,B$ verifies Condition~\ref{con3} with $\al$ defined as in \eqref{al} and $\beta = 0$.
%Take $x\in A$ and $y \in P_B(x)$ such that $d(y,A) > d(A,B)$ and let $x'\in P_A(y)$. Note that if there is no such $x$ and $y$, Condition 2 holds trivially.\todo[inline]{MKT: I struggling to see why it holds. Maybe its better to briefly explain the trivial verification?}
%\todo[inline]{MKT: There are a lot of conditions that $x^*$ and $y^*$ satisfy. I struggling to see why such a point exists. What the reason? Maybe it is worth adding a short explanation in the proof?}
%\todo[inline]{HB 28/05/21: I will update more explanations. It is not so clear I admit.}

%Suppose $d(x_0,B)<\frac{d(A,B)}{1-\al^2}$, but $d(x_0,B) > d(A,B)$. Let $y\in P_B(x_0)$.
%Suppose $x_0\notin P_A(y)$, it holds that $y-x_0 \notin \R_+(c)$. Furthermore, $y\in P_B(x_0)$ implies $x_0-y \in N_B(y)$. Together these inclusions imply $-c\notin N_B(y)$, or equivalently, \todo{MKT: Why do we need to exclude $v=0$? Won't the distance be positive in that case too since $\|c\|=1$?}
We first prove that the constant $\al_0$ defined in \eqref{al1} is positive. If the set $J(x^*)$ is empty, then \eqref{al1} becomes $\al_0 = \min\{1,\inf \emptyset\}=1 >0$. So, it suffices to consider the case where $J(x^*)\neq \emptyset$. To this end, let $J\in J(x^*)$ be arbitrary. By definition, $x^*\notin\cone[\mathbb A^*\,e_i]_{i\in J}$ and hence
\begin{equation}\label{ax}
d \left(x^*,\cone[\mathbb{A}^*\,e_i]_{i\in J}\right) >0.
\end{equation}
% \textcolor{green}{EN: 18/06  \( \inf \emptyset = \infty \) is inconsistent with the generic definition of \(\inf\), looks awkward.}
% \textcolor{red}{RSB 21/06: If C is the empty set, then every
% real number is both an upper and a lower bound of C, hence the largest lower bound is $\infty$, the same for the sup of C, so the smallest upper bound is $-\infty$. I hope this responds your concern?}
% \
From inequality \eqref{ax} and the fact that the set $J(x^*)$ is finite, the constant $\al_0$ in \eqref{al1} must be positive. 
% \textcolor{green}{EN:20/6  What I had in mind is that $\inf_{x \in \emptyset \subset \R } x  = \infty $ suits  better the definition of $\inf$ as we have no definition of inequalities between sets. There is no definition of a set {\bf lower} than the other, so we can not define the {\bf lowest} set which would be the $\inf$. Writing $\inf \emptyset = \infty$ kind of implies that $\emptyset = \infty$ :) Of course, any professional will understand right away what you mean but better stick to the rules. Somewhere in the future this paper will be read by robots and they might be confused !  }
% \textcolor{magenta}{MKT 22/06: The convention ``$\inf\emptyset=\infty$'' is common place (eg first page of Rockafellar and Wets). In to RSB's explanation, I'm not concerned about any misunderstanding and think we should leave as is.}
%{\color{red}[MKT: What if $J(x^*)=\emptyset$? Perhaps it's better to write this in a different way?]}
% Since $\R^m$ has finite basic[{\color{red}MKT: finite basis?}], $\al = \frac{1}{2} \min\left\{1,\min_{\substack{
% 		e^*\in E^m;\\ x^*\notin \R_{+}(-\mathbb{A}^*e^*)}} d \left(x^*,\R_{+}(-\mathbb{A}^*e^*)\right) \right\}$ must be positive [{\color{red}MKT: I think we don't need ``strictly'' here. Also should we give details on the case when there is no $e^*$ such that $x^*,\R_{+}(-\mathbb{A}^*e^*$? In this case, are we using the convention that $\min\emptyset=+\infty$?}].

To complete the proof, we show that the constant $\al_0$ defined in \eqref{al1} is a lower bound of all constants $\al$ such that $A$ is $\alpha$-sharp w.r.t.\ $x^*$. In other words, we will show that $\text{sr}[A,x^*]\ge \al_0 >0$. Observe that the polyhedral set $A$ is the intersection of $m$ closed half spaces. Namely, using Definition \ref{def:Polyhedral set} we can write
\begin{equation}\label{P8.3P1}
A = \bigcap_{i=1}^m \left\{x:\; \ang{e_i,\mathbb{A}x} \le \ang{e_i,b}\right\}=  \bigcap_{i=1}^m \left\{x:\; \ang{\mathbb{A}^*\,e_i,x} \le \ang{e_i,b}\right\},
\end{equation}
where we used \eqref{def:A*} in the second equality.
%aca
%Proposition~\ref{fact:Normal Cone} implies that
%{\color{blue}Indeed, take $x\in A$, then $\mathbb{A}\, x \le b$, which by the definition implies $(\mathbb{A}\, x)_i \le b_i$ for $i=1,\ldots,m$. From this, we obtain $\ang{e_i^*,\mathbb{A}\,x} \le \ang{e_i^*,b}$ for $i=1,\ldots,m$. Therefore, $x\in  \bigcap_{i=1}^m \left\{x:\; \ang{e_i^*,\mathbb{A}\,x} \le \ang{e_i^*,b}\right\}$; and so, $A \subset \bigcap_{i=1}^m \left\{x:\; \ang{e_i^*,\mathbb{A}\,x} \le \ang{e_i^*,b}\right\}$. Conversely, if vector $x$ belongs to the intersection $\bigcap_{i=1}^m \left\{x:\; \ang{e_i^*,\mathbb{A}\,x} \le \ang{e_i^*,b}\right\}$, then we must have $(\mathbb A\,x)_i \le b_i$ for $i=1,\ldots,m$. So, the inverse inclusion holds.}
%{\color{red}[MKT: Maybe I'm missing something but I don't see why we need the blue stuff. Isn't immediate from the definition?]}
Because $\mathbb A$ is a bounded linear map, so is $\mathbb{A^*}$, and hence for each $i=1,\ldots,m$, the set $H_i := \left\{x:\; \ang{\mathbb{A}^*\,e_i,x} \le \ang{e_i,b}\right\}$ is a closed half space. Thus, the normal cone operator of $H_i$ at a point $x\in H_i$ is given by
$$
N_{H_i}(x) = \begin{cases}
\cone [\mathbb{A}^*\,e_i],&\quad \text{if }\ang{\mathbb{A}^*\,e_i,x} = \ang{e_i,b};\\
\{0\},&\quad \text{otherwise}.
\end{cases}
$$
%{\color{red}[See Theorem 6.46 from Rockafellar and Wets textbook]}
We now apply the intersection rule in Proposition~\ref{fact:Normal Cone} and \eqref{eq:cones1} to derive the normal cone of $A$ as
\begin{equation}
    \label{P3.8P2}
    N_A(x) =\cone[\mathbb{A}^*e_i]_{i\in I(x)}=\sum_{i\in I(x)}\cone[\mathbb{A}^*e_i],\quad \forall x\in A,
\end{equation}
where $I(x) := \{i:\,\ang{x,\mathbb{A}^*\,e_i} = b_i\}$. 
%{\color{magenta}RSB 12/01: Why is the formula above true?}
% {\color{magenta} In other words, $i\in I(x)$ then we must have $j\notin I(x)$ and
% consequently $e^*_j\in E^*$.} 
Consider any $x\in A$ such that $x^*\notin N_A(x)$. Then, by \eqref{P3.8P2} and the definition of $J(x)$, we have $I(x)\in J(x)$, and 
$$
d\left(x^*,N_A(x)\right)= d\left(x^*,\cone[\mathbb{A}^*\,e_i]_{i\in I(x)}\right) \ge \al_0,
$$
where we used the definition of $\al_0$ in the last inequality.
% {\color{magenta}RSB 12/01: If the formula for $N_A$ is correct, then the equality above holds for all $x$.}
% {\color{magenta}Taking infimum over all $x\in A$ such that $x^*\notin N_A(x)$ we obtain
% \[
% \inf_{\substack{x\in A\\ x^*\notin N_A(x)}} d(x^*,N_A(x))\le \min_{e^*_j\in E^*} d\left(x^*,\R_{+}(-\mathbb{A}^*e^*_j)\right).
% \]
% }
%By \eqref{al1}, this implies $\al_0 \le d\left(x^*,N_A(x)\right)$ for any $x\in A$ such that $x^*\notin N_A(x)$. 
Hence,
$$
\al_0 \le \inf_{\substack{x\in A\\ x^*\notin N_A(x)}} d(x^*,N_A(x)) = \text{sr}[A,x^*],
$$
which implies the set $A$ is $\al_0$-sharp w.r.t. vector $x^*$.
\end{proof}
% Sharpness of a set, as introduced above, has a close connection to other regularity properties. In what follows, we show that the {\em Kurdyka-\L{}ojasiewicz} (KL) property for a function can be characterized in terms of the sharpness of its epigraph. Before showing this fact, we recall the definition of the KL property from \cite{li2018calculus}. Even though the latter reference is set in a finite dimensional framework, the definition below can be stated without changes in our more general setting.
% \begin{definition}[{\cite[Definition 2.3]{li2018calculus}}]
% \label{def:KL}
% %{\color{red}[MKT: Why do we have the function $\psi$ in the definition? Also the general definition (for nonconvex $f$ uses the limiting subdifferential but we use convex here? I think it would also be good to put a reference to the definition.]}
% Consider a proper closed convex function $f$, $\bx \in \dom \partial f$ and $p \in [0,1)$.
% %and a function %$cs^{1-\al}$ for some $c>0$
% If there exist $\al,\epsilon>0$ and ${v} \in (0,\infty]$ such that
% $$
% d(0,\partial f(x)) \ge \al (f(x) - f(\bar x))^p,
% $$
% whenever $\norm{x-\bx} \le \epsilon$ and $f(\bar x) < f(x) < f(\bx) +{v}$, then we say that $f$ has the {\em Kurdyka-\L{}ojasiewicz (KL) property at $\bx$ with exponent $p$}. 
%\end{definition}
\subsection{Sharpness of the epigraph}
Consider the following optimization problem 
\begin{equation}
    \label{prob:CP}
    \min_{x\in A} f(x),
\end{equation}
where $A$ is nonempty, closed and convex and $f:\mathcal H\to \R_{\infty}$ is proper, convex and lsc.
Assume that \eqref{prob:CP} has solutions and denote its solution set by $\mathbb{S}$. In connection with this problem, we consider the following property. Assume that there exists $\beta>0$ such that
\begin{equation}
    \label{eq:extended-KL}
    \inf_{x\notin \mathbb{S}} d(0, \partial f(x))\ge \beta>0.
\end{equation}
\begin{remark}
\label{rem:pwl}
We note that \eqref{eq:extended-KL} holds, for instance, if $f$ is piece-wise linear. Indeed, let $f(x):=\max_{i=1,\ldots,m} \ang{x_i^*,x}+b_i$. Call $J:=\{i,\ldots,m\}$. For every $x\in \mathcal{H}$, define $I(x):=\{i\in J\::\: f(x)=\ang{x_i^*,x}+b_i\}$. From the formula of the subdifferential of a supremum (see, e.g., \cite[Theorem 18.5]{bauschke2011convex} or \cite[Theorem 2.4.18]{Za}), we know that  $\partial f(x)={\rm co}[x_i^*]_{i\in I(x)}$. Hence, $0\notin \partial f(x)$ if and only if $0\notin {\rm co}[x_i^*]_{i\in I(x)}$. 
Consider
\[
{\cal F}:=\{I\subset J \::\: 0\notin {\rm co}[x_i^*]_{i\in I}\}\subset {\cal P}(J),
\]
where ${\cal P}(J)$ denotes the collection of all possible subsets of $J$. Recall that  ${\cal P}(J)$ has cardinality $2^m$. So there is a finite number of possible subdifferential sets such that $0\notin \partial f(x)$. Moreover, $x\notin \mathbb{S}$ if and only if $I(x)\in {\cal F}$.
Altogether,
\[
\inf_{x\notin \mathbb{S}} d(0, \partial f(x))= \inf_{0\notin \partial f(x)}d(0, \partial f(x))\ge \min_{I\in {\cal F} } d(0, {\rm co}[x_i^*]_{i\in I}),
\]
where we have a minimum in the rightmost expression because the infimum is taken over the finite set ${\cal F}$. Since ${\rm co}[x_i^*]_{i\in I}$ is a closed convex cone which doesn't contain zero, $d(0, {\rm co}[x_i^*]_{i\in I})>0$ for every $I\in {\cal F}$. Therefore, the minimum in the right-hand side of the above expression is attained at some positive value $\beta$. Note that the fact that we have a finite supremum of affine functions is essential. The function in Example \ref{F2a} is an infinite supremum of affine functions, for which \eqref{eq:extended-KL} does not hold.
\end{remark}

The next result characterizes the case in which the epigraph of a function is sharp w.r.t.\ the vector
$(0_{\mathcal H},-1)$. The latter property turns out to be important in Section \ref{S2}. 

% \textcolor{magenta}{MKT 27/06: Just an FYI in case you don't realise. The reference you use here (ie \cite[Definition 2.3]{li2018calculus}) defines this property only in the context of $R^n$.}

\begin{proposition}
\label{EXP1}
Let $f:\mathcal H\to \R_{\infty}$ be a proper, convex and lsc function. The following statements are equivalent.
\begin{itemize}
    \item[(i)] $\epi f$ is $\al$-sharp w.r.t.\ the vector $z^*=:(0_{\mathcal H},-1)$. 
    \item[(ii)] $\al<1$ and $f$ verifies \eqref{eq:extended-KL} with parameter $\beta:=\frac{\al}{\sqrt{1-\alpha^2}}$.
\end{itemize}
\end{proposition}
\begin{proof} 
Let $E:=\epi f$, and $D:=\dom f$. From Fact~\ref{Normal cone of epigraph}, we have
    \begin{equation}\label{normal cone epi}
    \begin{array}{rcl}
    N_E((x,f(x))=\R_{++}(\partial f(x)\times \{-1\} ) \bigcup (N_D(x)\times \{0\}).
    \end{array}
  \end{equation}
For calculating the sharpness of $E$ at $z^*$, it is enough to consider only points of the form $(x,f(x))$, since \textcolor{black}{otherwise for $(x,y) \in E$, with $y> f(x)$, then $N_E((x,y)) = N_{\dom f}(x)\times \{0\}$, and hence $d(z^*,N_E((x,y)) \ge 1$}. Thus, for computing the sharpness of $E$ we need to take the infimum over the set:
\[
K:=\{(x,f(x)): (0,-1)\not\in N_E((x,f(x))\}=\{(x,f(x)): \partial f(x)\not=\emptyset \hbox{ and }0\not\in \partial f(x), \hbox{ or } \partial f(x)=\emptyset\},
\]  
where we are using \eqref{normal cone epi} in the characterization of $K$. If $(x,f(x))\in K$ and $\partial f(x)=\emptyset$, then \eqref{normal cone epi} gives 
$ N_E((x,f(x))=N_D(x)\times \{0\}$,
so in this case we have
\begin{equation}
    \label{eq:case empty}
    d\left(z^*,N_E((x,f(x)))\right)   = d((0_{\mathcal H},-1), N_D(x)\times \{0\})=\inf_{w\in N_D(x)} \sqrt{ \|w\|^2 +1 }=1.
\end{equation}
Fix now any $(x,f(x))\in K$ such that $\partial f(x)\not=\emptyset$ and denote by $\de(x):=d(0,\partial f(x))$. Since $\partial f(x)$ is a nonempty, closed and convex set (\textcolor{blue}{see \cite[Proposition~20.31]{bauschke2011convex}}), and $0\notin \partial f(x)$, we have that $\de(x)>0$. \textcolor{black}{We first prove that for every $(x,f(x))\in K$ such that $\partial f(x)\not=\emptyset$, we have
\begin{equation}\label{eq:claim2}
     d\left(z^*,N_E((x,f(x)))\right)   =  \frac{\delta(x)}{\sqrt{\delta(x)^2+1}}.
  \end{equation}
  Indeed, fix $(x,f(x))\in K$. Using \eqref{normal cone epi} we have 
\begin{align*}
    d\left(z^*,N_E((x,f(x)))\right)   
    &= \inf_{(u,\eta)\in N_E((x,f(x))} d(z^*, (u,\eta))
     = \min\left\{ \inf_{w\in N_D(x)} d(z^*, (w,0)), \inf_{u\in \partial f(x),\, t>0} d(z^*, (tu,-t))\right\}\\
     &= \min\left\{ \inf_{w\in N_D(x)} \sqrt{ \|w\|^2 +1 }, \inf_{{u\in\partial f(x)},\,t> 0}   \sqrt{t^2(\|u\|^2+1)-2t+1}\,\right\}\\
      & = \min \left\{  1,  \inf_{u\in\partial f(x)} \left[ \inf_{t>0}\sqrt{t^2(\|u\|^2+1)-2t+1}\right] \right\}
       = \min \left\{  1,  \inf_{u\in\partial f(x)} \frac{\|u\|}{\sqrt{\|u\|^2+1}}  \right\}\\
      &= \frac{\delta(x)}{\sqrt{\delta(x)^2+1}}, 
   \end{align*}
where we are using \eqref{normal cone epi} in the second equality and the definition of $z^*$ in the third one. In the fourth equality we are using the fact that the infimum over $N_D(x)$ is attained at $w=0$. In the fifth equality we are using the fact that the infimum in the expression between square brackets is attained at $t^*:= 1/(\|u\|^2 +1)$, and in the sixth one the fact that the latter infimum value is smaller than $1$ as well as the fact that the function $g(s):=s/\sqrt{s^2+1}$ is increasing over $\R_+$ and attains its minimum at $s^*:=\de(x)$. Hence, \eqref{eq:claim2} holds for every $(x,f(x))\in K$.
  }

Now, we assume that (i) holds. Then by \eqref{eq:case empty} and \eqref{eq:claim2} this means that, for every $(x,f(x))\in K$, we have
\begin{align*}
  \al   &\le   \inf_{(x,f(x))\in K} d\left(z^*,N_E((x,f(x)))\right)  \\
     & = \min\left\{ \inf_{(x,f(x))\in K, \partial f(x)=\emptyset} d\left(z^*,N_E((x,f(x)))\right), \inf_{(x,f(x))\in K, \partial f(x)\not=\emptyset} d\left(z^*,N_E((x,f(x)))\right) \right\}\\
     &= \min\left\{1, \inf_{(x,f(x))\in K, \partial f(x)\not=\emptyset} \frac{\delta(x)}{\sqrt{\delta(x)^2+1}} \right\} \le \frac{\delta(y)}{\sqrt{\delta(y)^2+1}} <1,
\end{align*}
for every $(y,f(y))\in K$ such that $\partial f(y)\not=\emptyset$. Note that we are using 
\eqref{eq:case empty} and \eqref{eq:claim2} in the second equality. So $\al<1$, and the above expression re-writes as
\begin{equation}
    \label{eq:equiv}
    \de(y)\ge \frac{\al}{\sqrt{1-\alpha^2}},\,\forall (y,f(y))\in K \hbox{ s.t. } \partial f(y)\not=\emptyset,
\end{equation}
which means that (ii) holds with parameter $\frac{\al}{\sqrt{1-\alpha^2}}$. Note that we are using the convention that the infimum of the empty set is $+\infty$, so (ii) automatically holds if $\partial f(x)=\emptyset$.

Conversely, if condition (ii) holds with parameter $\beta=\frac{\al}{\sqrt{1-\alpha^2}}$, then \eqref{eq:equiv} holds. The latter rewrites as
\[
\frac{\delta(y)}{\sqrt{\delta(y)^2+1}}\ge \al,\,\forall (y,f(y))\in K \hbox{ s.t. } \partial f(y)\not=\emptyset,
\]
which, together with \eqref{eq:case empty} and \eqref{eq:claim2}, give (i). Therefore, it is enough to establish \eqref{eq:claim2}.
\end{proof}

\begin{figure*}[h]
	\centering
	\subfigure[The function $f(x) = x^2$ does not have the global KL property, and the set $\epi f$ is not sharp w.r.t. the vector $(0,-1)$.]{\makebox[6.5cm][c]{
	\begin{tikzpicture}
	\fill[scale=0.7, domain=-2:2, smooth, variable=\x, gray,opacity=0.1] plot ({\x}, {(\x)^2});
	\draw[thick,scale=0.7, domain=-2:2, smooth, variable=\x, black,opacity=0.5] plot ({\x}, {(\x)^2});
	\draw[thick,black] (2,0)--(-2,0);
	 \draw (0,2) node {$\epi f$};
	 \draw[->,densely dotted, black] (1,0)--(1,-0.5);
	 \draw (1.2,-0.8) node {$(0,-1)$};
%  \addplot[domain=-3:3,thick,black,no marks,smooth,fill=blue!30] {x^2};
\end{tikzpicture}
}
}
\quad \quad \quad \quad\quad\quad
	\subfigure[The function $g(x) = \abs{x}$ has the global KL property, and the set $\epi g$ is sharp w.r.t. the vector $(0,-1)$.]{\makebox[6.5cm][c]{
	\begin{tikzpicture}
	\fill[scale=0.7, domain=-2:2, smooth, variable=\x,gray,opacity=0.1] plot ({\x}, {2*abs(\x)});
	\draw[thick,scale=0.7, domain=-2:2, smooth, variable=\x, black,opacity=0.5] plot ({\x}, {2*abs(\x)});
	\draw[thick,black] (2,0)--(-2,0);
	 \draw (0,2) node {$\epi g$};
	 \draw[->,densely dotted, black] (1,0)--(1,-0.5);
	 \draw (1.2,-0.8) node {$(0,-1)$};
    %\addplot[domain=-3:3,thick,black,no marks,smooth,fill=blue!30] {x^2};
    \end{tikzpicture}
	}
}
	\caption{Illustration of Proposition~\ref{EXP1}: $\epi f$ is sharp w.r.t. the vector $(0,-1)$ if and only if $f$ has the global KL property. The rationale for these figures is identical to the one given earlier for Figure \ref{F1}.}
% \todo[inline]{Todo: revise captions}
		\label{F3}
\end{figure*}
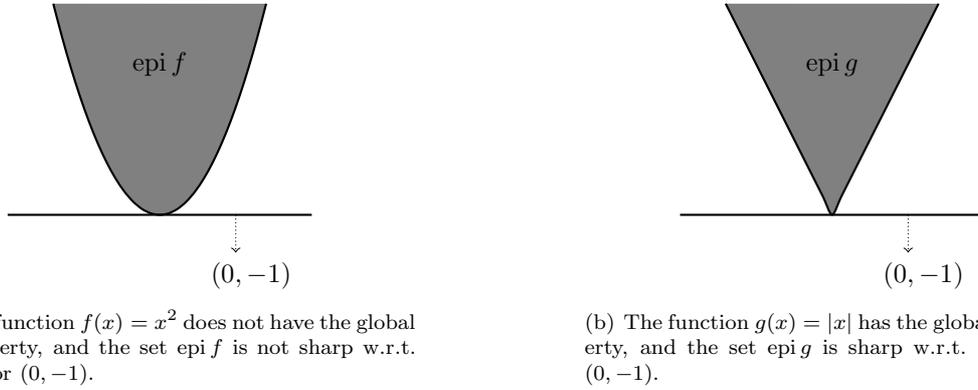

We illustrate the results in Proposition~\ref{EXP1} by two examples in Figure~\ref{F3}. Condition \eqref{eq:extended-KL} is closely related with the well-known Kurdyka-\L{}ojasiewicz inequality. To make this connection precise, we recall next the necessary definitions. 

{\color{black}
\begin{definition}[Kurdyka-\L{}ojasiewicz inequality {\cite[Section 2.3]{Bolte2017}}]\label{def:KL2}
    Let $f\colon\mathcal{H}\to\R_{\infty}$, and assume that $\mathbb{S}:=\argmin f\neq\emptyset$. Fix $\bar x\in \mathbb{S}$. The function $f$ satisfies the \emph{global Kurdyka-\L{}ojasiewicz (KL) property} at $\bar x$ if there exists a concave continuously differentiable function $\varphi\colon \R_+\to \R_+$ with $\varphi(0)=0$ and $\varphi'>0$ such that
    \begin{equation}\label{eq:KL}
    \varphi'(f(x)-f(\bar x))\,d(0,\partial f(x)) \geq 1\quad\forall x\notin  \mathbb{S}.
    \end{equation}
In this case, we say that $\varphi$ is a \emph{desingularizing function} for $f$ at $\bar{x}$.  
If $f$ satisfies the global KL property and admits the same desingularizing function $\varphi$ at every point $\bar x\in\mathbb{S}$, then we say that $f$ satisfies the \emph{global KL property} with \emph{global desingularizing function} $\varphi$.
\end{definition}
}
The next result establishes the connection between the global KL property and sharpness.

\begin{corollary}\label{cor:sharp and KL}
Let $f\colon\mathcal{H}\to\R_{\infty}$ be a proper lsc convex function. The following statements are equivalent.
\begin{itemize}
    \item[(i)] $\epi f$ is $\al$-sharp w.r.t.\ the vector $z^*=:(0_{\mathcal H},-1)$. 
    \item[(ii)] $\al<1$ and $f$ satisfies the global KL property with global desingularizing function $\varphi(t)=\frac{t\sqrt{1-\alpha^2}}{\al}$.
\end{itemize}
\end{corollary}
\begin{proof}
The claim of the corollary follows from Proposition~\ref{EXP1} because the global KL property (ii) is equivalent to condition (b) in Proposition \ref{EXP1}. 
\end{proof}
%aca
\begin{corollary}\label{cor: f polyhedral}
Let $f\colon\mathcal{H}\to\R_{\infty}$ be a polyhedral function (see Definition \ref{def:Polyhedral set}). Then there exists $\al<1$ such that $f$ satisfies the global KL property with global desingularizing function $\varphi(t)=\frac{t\sqrt{1-\alpha^2}}{\al}$. In particular, property \eqref{eq:extended-KL} holds for $\beta:=\frac{\al}{\sqrt{1-\alpha^2}}$.
\end{corollary}
\begin{proof}
By Definition \ref{def:Polyhedral set}, the epigraph of a polyhedral function is a polyhedral set. By Proposition \ref{active}, it is sharp with respect to any unit vector, in particular with respect to $z^*:=(0_{\mathcal{H}},-1)$. The two claims now follow from the part (i) implies (ii) in  Corollary \ref{cor:sharp and KL}, and from the part (i) implies (ii) in Proposition \ref{EXP1}.
\end{proof}

% \begin{figure*}[h]
% 	\centering
% 	\subfigure[The function $f(x) = x^2$ does not have KL property with exponent $0$, and the set $\epi f$ is not sharp w.r.t. vector $(0,-1)$.]{\makebox[6.5cm][c]{
% 	\begin{tikzpicture}
% 	\fill[scale=0.7, domain=-2:2, smooth, variable=\x, blue!30] plot ({\x}, {(\x)^2});
% 	\draw[thick,scale=0.7, domain=-2:2, smooth, variable=\x, blue] plot ({\x}, {(\x)^2});
% 	\draw[thick,blue] (2,0)--(-2,0);
% 	 \draw (0,2) node {$\epi f$};
% 	 \draw[->,red] (1,0)--(1,-0.5);
% 	 \draw (1.2,-0.8) node {$(0,-1)$};
% %  \addplot[domain=-3:3,thick,black,no marks,smooth,fill=blue!30] {x^2};
% \end{tikzpicture}
% }
% }
% \quad \quad \quad \quad\quad\quad
% 	\subfigure[The function $g(x) = \abs{x}$ has KL property with exponent $0$, and the set $\epi g$ is sharp w.r.t vector $(0,-1)$.]{\makebox[6.5cm][c]{
% 	\begin{tikzpicture}
% 	\fill[scale=0.7, domain=-2:2, smooth, variable=\x, blue!30] plot ({\x}, {2*abs(\x)});
% 	\draw[thick,scale=0.7, domain=-2:2, smooth, variable=\x, blue] plot ({\x}, {2*abs(\x)});
% 	\draw[thick,blue] (2,0)--(-2,0);
% 	 \draw (0,2) node {$\epi g$};
% 	 \draw[->,red] (1,0)--(1,-0.5);
% 	 \draw (1.2,-0.8) node {$(0,-1)$};
%     %\addplot[domain=-3:3,thick,black,no marks,smooth,fill=blue!30] {x^2};
%     \end{tikzpicture}
% 	}
% }
% 	\caption{Illustration of Proposition~\ref{EXP1}: $\epi f$ is sharp w.r.t. vector $(0,-1)$ if and only if $f$ has KL property with exponent $0$. The rationale for these figures is identical to the one given earlier for Figure \ref{F1}.}
% 		\label{F3}
% \end{figure*}

In the next result, we establish yet another connection between the sharpness property of a set and the KL property. The result shows that a set is sharp w.r.t.\ $x^*$ if and only if the function $\mathbbm{1}_A(\cdot)-\ang{x^*,\cdot}$ has the KL property with exponent of $0$. 

\begin{proposition}
\label{EXP2}
Fix $x^*\in \mathcal S$ and consider a closed and convex set $A$. The following statements are equivalent.
\begin{itemize}
\item[(i)] The modulus of sharpness of $A$ is w.r.t.\ $x^*$ is $\al>0$ or, equivalently, $\al:=\text{sr}[A,x^*] >0$.
\item[(ii)] The function $f(x):=\mathbbm{1}_A(x)- \ang{x^*,x}$ satisfies the global KL property with global desingularizing function $\varphi(t)=\frac{t}{\al}$.
\item[(iii)] The function $f$ as in (ii) satisfies \eqref{eq:extended-KL} with $\beta:={\al}$.
%\item[(c)] $f$ does not attain a minimum over $A$.
%\item[(d)] $x^*\notin N_A(x)$ for every $x\in A$. 
\end{itemize}
If, under any of the above conditions, we have that $\al<1$, then $\epi f$ is $(\al/\sqrt{1+\al^2})$-sharp w.r.t. $z^*:=(0_{\mathcal{H}},-1)$.
%The set if and only if the  
\end{proposition}

\begin{proof} The definition of $f$ yields
$$
\partial f(x) = N_A(x)-x^*,\quad \forall x\in \mathcal H.
$$
%{\color{red}[MKT: I think formula should be ``$\partial f(x) = N_A(x)-x^*$'']}
Thus, for any $x\in A$, $x^* \notin N_A(x)$ if and only if $0\notin \partial f(x)$. In other words, $x\notin \argmin f$ if and only if $x^* \notin N_A(x)$. The equality above implies that
%. Therefore, if $x$ is not a miminizer of $f$, we have $x^*\notin N_A(x)$, and hence
%{\color{red}[MKT: Why do we need to talk ``minmiziers of $f$'' in this proof? It seems sufficient to just use ``$0\in\partial f(x)$'' which follows even more directly.]}
\begin{equation}\label{equiv}
    d(0,\partial f(x))= d(0,N_A(x)-x^*) = d(x^*,N_A(x)),
\end{equation}
which, combined with the fact that $\varphi'(t)=1/\al$ gives, for all $x\notin \argmin f$, $\bar x\in \argmin f$,
\[
(1/\al)\,d(0,\partial f(x))= \varphi'(f(x)-f(\bar x))\,d(0,\partial f(x))= (1/\al) d(x^*,N_A(x)),
\]
The three equivalences follow directly from the expression above and the definitions. As for the last statement, assume that $\al<1$. Since all (i)-(iii) are equivalent, we can use Proposition \ref{EXP1}(ii) for $\beta:=\al$. Indeed, by (iii) we have that $\de(x)\ge \al$. By \eqref{eq:claim2} this implies that, whenever $\partial f(x)\not=\emptyset$ and $0\notin \partial f(x)$
\[
 d\left(z^*,N_E((x,f(x)))\right)   =  \frac{\delta(x)}{\sqrt{\delta(x)^2+1}}\ge \frac{\al}{\sqrt{1+\al^2}},
\]
where we used the fact that the function $g(s):=s/\sqrt{s^2+1}$ is increasing over $\R_+$. By \eqref{eq:case empty} we know that
\[
 d\left(z^*,N_E((x,f(x)))\right)   =  1, \hbox{ if }\partial f(x)=\emptyset.
\]
Altogether, and using the same set $K$ as in the proof of Proposition \ref{EXP1}, we can write
\[
\begin{array}{rcl}
  \inf_{(x,f(x))\in K}d\left(z^*,N_E((x,f(x))\right)  
 & =&  \min\left\{\!\! \inf_{\tiny\begin{array}{c}(x,f(x))\in K,\\ \partial f(x)=\emptyset\end{array}} d\left(z^*,N_E((x,f(x))\right), \inf_{\tiny\begin{array}{c}(x,f(x))\in K,\\ \partial f(x)\not=\emptyset\end{array}} d\left(z^*,N_E((x,f(x))\right) \right\}\\
     &&\\
     &\ge & \min\left\{1, \frac{\al}{\sqrt{1+\al^2}}\right\} = \frac{\al}{\sqrt{1+\al^2}},\\
\end{array}
\]
establishing the last claim.
\end{proof}
\subsection{Dual characterizations}\label{Dual}

As seen in Figure~\ref{F1}, the sharpness property of the set $A$ w.r.t.\ the vector $x^*$ is related to the condition that $x^*$ belongs to the interior of the set $\cup_{x\in F_A(x^*)}N_A(x)$. The next result explores this connection in general settings. 

We recall that the open ball of radius $r>0$ and center $x_0$ is written as $x_0+r\, (\mathcal B\!\setminus{\mathcal S})$, and that the corresponding closed ball is $x_0+r\, \mathcal B$. 
%5 June

% \textcolor{red}{EN 4 june: We already have $\mathcal B$ and $\mathcal S$, why not
% $v + \alpha {\mathcal B}$ and $v + \alpha ({\mathcal B} \setminus \mathcal S) $ ?
% Or introduce in Preliminaries the open unit ball
% $\mathcal B^0 = {\mathcal B} \setminus {\mathcal S}$ and $v +  \mathcal B^0$ ?
% (ii) will be easier to grasp.}

% \textcolor{blue}{RSB 4 June: OK, I will do so.}

% \textcolor{red}{ Thnx, Regina. Will quit for a moment not to mess with you editing}

\begin{proposition}\label{PP3}
%\label{PP3}\textcolor{magenta}{MKT 27/06: Do we need to assume $A$ is nonempty? For instance, (ii)  can never hold when $A$ is empty because we would have ``$x^*+\alpha(\mathcal{B}\setminus\mathcal{S})\subseteq\emptyset$''.}
%\textcolor{blue}{HB I add nonempty assumption}
Consider a nonempty closed convex set $A$ of a Hilbert space $\mathcal H$, a vector $x^*\in \mathcal S$, and a positive constant $\al \in (0,1]$. Consider the following statements:
\begin{enumerate}[(i)]
    \item The set $A$ is $\al$-sharp w.r.t.\ $x^*$.
%    \item $B_\al(x^*)\subset \cl\bigcup_{x\in F_A(x^*)}N_A(x)$.
    \item $x^*+\alpha (\mathcal B\!\setminus{\mathcal S})\subset \bigcup_{x\in F_A(x^*)}N_A(x)$.
\end{enumerate}
Then, (ii) $\Rightarrow$ (i). If the set $A$ is bounded, then (i) $\Rightarrow$ (ii). 
%If $H$ finite dimensional and $A$ is bounded, then all three statements are equivalent.
\end{proposition}

\begin{proof}

%\item The implication (iii) $\Rightarrow$ (ii) is trivial. 
We first prove (ii) $\Rightarrow$ (i). 
Suppose (ii) holds.
%Suppose $x^*+\alpha (\mathcal B\!\setminus{\mathcal S})\subset \bigcup_{x\in F_A(x^*)}N_A(x)$ and 
We now prove that, for any vector $y\in A$ satisfying
\begin{equation}
    \label{PP3.1}
    d(x^*,N_A(y))<\al,
\end{equation}
we must have $x^*\in N_A(y)$. This statement then implies that inequality \eqref{DC} must hold for any $x\in A$ with $x^*\notin N_A(x)$.

Suppose to the contrary that there is $y\in A$ such that \eqref{PP3.1} holds and $x^*\notin N_A(y)$. From the inequality \eqref{PP3.1}, we can choose $y^*\in N_A(y)$ such that $0<\norm{x^*-y^*}<\al$. Hence, $y^*\in x^*+\alpha (\mathcal B\!\setminus{\mathcal S})$.
% and assumption (ii) yields
% \[
% y^*\in \bigcup_{x\in F_A(x^*)}N_A(x).
% \]
%If $x^*\in N_A(y)$, then the set $A$ must be $\al$-sharp w.r.t. vector $x^*$. We suppose to the contrary that $x^*\notin N_A(y)$
%
Because $x^*+\alpha (\mathcal B\!\setminus{\mathcal S})$ is an open set, then there exists $r>0$ such that $y^*+r\mathcal B\subset x^*+\alpha (\mathcal B\!\setminus{\mathcal S})$. Let $t_0:= 1+\frac{r}{\norm{y^*-x^*}}$ and $z^*:= x^*+t_0(y^*-x^*)$ (see Figure~\ref{P3.9PF1}). Using these definitions, we can write
\[
\begin{array}{rcl}
 \norm{z^*-y^*}    &= &\norm{\left(x^*+\left(1+\frac{r}{\norm{y^*-x^*}}\right)(y^*-x^*)\right) -y^*} \\
 &&\\
     &= &\norm{\frac{r}{\norm{y^*-x^*}}(y^*-x^*)}=r.
\end{array}
\]
This gives 
\begin{equation}\label{eq:ii}
z^*=x^*+t_0(y^*-x^*) \in y^*+r\mathcal B\subset x^*+\alpha (\mathcal B\!\setminus{\mathcal S}),    
\end{equation}
where we used the definition of $r$ in the last inclusion. 
%Therefore, by the definition of $\tau$, we must have $\tau \ge 1+\frac{r}{\norm{y^*-x^*}} >1$ (see Figure~\ref{P3.9PF1}).
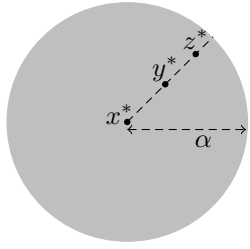
\begin{figure}[h]
	\centering
    \begin{tikzpicture}
     \filldraw[gray!50,opacity=0.2] (0,0) circle (45pt);
     \filldraw (0,0) circle (1pt);
     \filldraw (0.5,0.5) circle (1pt);
     \filldraw (0.9,0.9) circle (1pt);
     \draw[densely dashed,black] (0,0)--(1.13,1.13);
     \draw (-0.1,0.1) node {$x^*$};
     \draw (0.5,0.7) node {$y^*$};
     \draw (0.9,1.1) node {$z^*$};
     \draw[<->,densely dashed,black] (0,-0.1)--(1.56,-0.1);
     \draw (1,-0.25) node {$\al$};
    \end{tikzpicture}
	\caption{Here, $y^*$ belongs to the open segment $(x^*,z^*)$.}\label{P3.9PF1}
\end{figure}
Now we can use assumption (ii). Namely,
% \[
% y^*\in \bigcup_{x\in F_A(x^*)}N_A(x).
% \]
\[
z^*\in x^*+\alpha (\mathcal B\!\setminus{\mathcal S})\subset \bigcup_{x\in F_A(x^*)}N_A(x). 
\]
Hence, there exists $z\in F_A(x^*)$ such that $z^*\in N_A(z)$. Since $z\in F_A(x^*)$, Fact \ref{fact 1} yields $x^*\in N_A(z)$. Altogether, both of the vectors $x^*$ and $z^*$ belong to the normal cone $N_A(z)$. 
Using now the fact that $y\in A$ and $z^*\in N_A(z)$ we obtain 
\begin{equation}\label{eq:CC}
\ang{z^*,y- z}\le 0.
\end{equation}
% 
% $$
Recall that $x^*\notin N_A(y)$, and $x^*\in N_A(z)$. Using Fact~\ref{fact 1} the latter give $z\in \Argmax_A\ang{x^*,\cdot}$ and $y\notin\Argmax_A\ang{x^*,\cdot}$. Therefore, we can write
\begin{equation}\label{P3.7P1}
\ang{x^*,y}<\ang{x^*,z}=\sup_A\ang{x^*,\cdot}.
\end{equation}
From $y^*\in N_A(y)$ and $z\in A$, we also have
%Consider the following sets
%\begin{enumerate}
%    \item $F_0:=\{x\in A: x^*\in N_A(x)\}$ is the collection of all vectors in $A$ whose normal cones contain $x^*$;
%    \item $F_1:=\{x\in A: y^*\in N_A(x)\}$ is the collection of all vectors in $A$ whose normal cones contain $y^*$;
%    \item $F_2:=\{x\in A: z^*\in N_A(x)\}$ is the collection of all vectors in $A$ whose normal cones contain $z^*$.
%\end{enumerate}
%Then, there is $\bar z$ such that $\bar z\in F_0\cap F_1\cap F_2$. Because $y\notin F_0$, we have 
\begin{equation}
    \label{P3.P7P2}
    \ang{y^*,z} \le \ang{y^*,y}.
\end{equation}
We now use inequalities \eqref{P3.7P1} and \eqref{P3.P7P2} above, together with the facts that $z^*= x^*+t_0(y^*-x^*)$ and $t_0 >1$, to derive the following estimation
    \begin{align*}
        \ang{z^*,y}&=\ang{t_0y^*+(1-t_0)x^*,y} = t_0\ang{y^*,y}+(1-t_0)\ang{x^*,y}> t_0\ang{y^*,z}+(1-t_0)\ang{x^*,z}\\
        &= \ang{t_0y^*+(1-t_0)x^*,z}= \ang{z^*,z},
    \end{align*}
    which contradicts \eqref{eq:CC}. Therefore, we must have $x^*\in N_A(y)$. We have shown that, whenever \eqref{PP3.1} holds, we must have $x^*\in N_A(y)$. Equivalently, if $x^*\not\in N_A(y)$ then we must have $ d(x^*,N_A(y))\ge \al$. The latter statement implies that
 $A$ is $\al$-sharp w.r.t vector $x^*$. This completes the proof of (ii) $\Rightarrow$ (i).
 
 We now prove (i) $\Rightarrow$ (ii) when the set $A$ is bounded and $\al$-sharp w.r.t. $x^*$. Since the set $A$ is bounded, we can use the characterization of sharpness established in Proposition~\ref{PP01}. The latter implies that for all $y^*\in \mathcal H\setminus \{0\}$ such that $F_A(x^*)\cap F_A(y^*) =\emptyset$ we must have
\begin{equation}\label{eq:CC4}
    \norm{x^*-y^*} \ge \al.
\end{equation}
%From \cite[Corollary 21.25]{bauschke2011convex}, any vector $y^*\in H\setminus\{0\}$, the $\arg\max\ang{y^*,\cdot}$ exists, and hence $F_A(y^*)$ is nonempty.
We now show the inclusion in (ii). Take any $y^*\in x^*+\alpha (\mathcal B\!\setminus{\mathcal S})$, and assume that $y^*\notin \bigcup_{x\in F_A(x^*)}N_A(x)$. This implies that for every $x\in F_A(x^*)$, we must have $y^*\not\in N_A(x)$. By Fact~\ref{fact 1}, we deduce that for every $x\in F_A(x^*)$, we must have $x\not\in F_A(y^*)$. In other words, $F_A(y^*)\cap F_A(x^*)=\emptyset$. By (i) and Proposition~\ref{PP01} we deduce that \eqref{eq:CC4} holds. This contradicts the fact that $y^*\in x^*+\alpha (\mathcal B\!\setminus{\mathcal S})$. Therefore,
%we must have 
\(
y^*\in \cup_{x\in F_A(x^*)}N_A(x)
\),
which
%This
completes the proof of (i) $\Rightarrow$ (ii).
\end{proof}

\begin{remark}
In \cite{nurminski2015single}, a strict complementary condition is used to ensure that SPP solves the linear programming problem $\min_{x\in A}\ang{x^*,x}$ where $A\subset \mathcal H$ is a polyhedron, $x^*\in \mathcal H\setminus \{0\}$, and $\dim\mathcal{H}=n<+\infty$. This is established in \cite[Theorem 1]{nurminski2015single}, whose proof relies on the following two key conditions:
\begin{enumerate}[(i)]
    \item The minimization problem $\min_{x\in A}\ang{x^*,x}$ has an unique solution $\bx\in A$; and
    \item The normal cone $N_A(\bx)$ has nonempty interior and $-x^*\in \Int N_A(\bx)$ (see Figure~\ref{F1.1}(a)).
\end{enumerate}
Condition (ii) above implies that there must be at least $n$ linear independent constraints that are active at $\bar x$, where $n$ is the dimension of the primal problem. In fact, let $I(\bx)$ be the set of active constraints at $\bx$, then the normal cone of $A$ at $\bx$ is
$$
N_A(\bx)= \cone[\mathbb A^*\,e^*_i]_{i\in I(\bx)}.
$$
Since $\Int N_A(\bx) \neq \emptyset$, it follows that $\dim N_A(\bx) = n$. Therefore, $\abs{I(\bx)} \ge n$, and there are $n$ constraints $i_1,\ldots,i_n\in I(\bx)$ such that $\mathbb A^*\,e^*_{i_1},\ldots,\mathbb A^*\,e^*_{i_n}$ are linearly independent. 
Thus, any linear programming problem in $\R^n$ that satisfies the strict complementary condition must have at least $n$ linear independent constraints.
\textcolor{black}{By Definition~\ref{def:Polyhedral set}, a linear programming problem can only have finite number of constraints, hence the aforementioned requirement %that there are exactly $\dim \mathcal{H}$ linearly independent active constraints 
may only work for finite dimensions, and will never hold in infinite dimensional case.}%These arguments rely on finite dimensional arguments and cannot be extended to the infinite dimensional case. 
%\textcolor{magenta}{MKT 27/06: The last sentence is a little unclear to me. Eg Which arguement does ``These arguements'' refer to?}
%\textcolor{blue}{HB: I rewrite the paragraph}
%The next result explores the connections between the sharpness property and the nonemptiness of the interior of the normal cones.
\begin{figure*}[h]
	\centering
	\subfigure[Strict complementary condition in \cite{nurminski2015single} holds: the normal cone $N_A(\bx)$ has nonempty interior and $-x^*\in N_A(\bx)$.]{\makebox[6.5cm][c]{
% 	\begin{tikzpicture}[scale = 0.7]
% 	\fill[blue!30,opacity=1] (-2,3)--(0,0)--(2,3)--(-2,3);
% 	\draw[blue,thick] (-2,3)--(0,0)--(2,3);
% 	\fill[gray!30!red!20!blue!10!,opacity=1] (-2.25,-1.5)--(0,0)--(2.25,-1.5)--(-2.25,-1.5);
% 	\draw[->,thick,gray](0,0)--(-2.25,-1.5);
% 	\draw[->,gray,thick](0,0)--(2.25,-1.5);
% 	\draw[->,red!60](0,0)--(0,-1.5);
% 	\draw (0,2) node {$A$};
% 	\draw[red] (0.4,-1.2) node {$-x^*$};
% 	\draw (0.3,0) node {$\bx$};
% 	\draw (-1.2,-0.8) node {$N_A(\bx)$};
% 	\filldraw (0,0) circle (1pt);
%     %\addplot[domain=-3:3,thick,black,no marks,smooth,fill=blue!30] {x^2};
%     \end{tikzpicture}
	\begin{tikzpicture}[scale = 0.7]
	\fill[gray,opacity=.1] (-2,3)--(0,0)--(2,3)--(-2,3);
	\draw[black,opacity=0.5] (-2,3)--(0,0)--(2,3);
	\fill[gray,opacity=0.3] (-2.25,-1.5)--(0,0)--(2.25,-1.5)--(-2.25,-1.5);
	\draw[->,thick,gray](0,0)--(-2.25,-1.5);
	\draw[->,gray,thick](0,0)--(2.25,-1.5);
	\draw[->,densely dotted,black](0,0)--(0,-1.5);
	\draw (0,2) node {$A$};
	\draw (0.4,-1.2) node {$-x^*$};
	\draw (0.3,0) node {$\bx$};
	\draw (-1.2,-0.8) node {$N_A(\bx)$};
	\filldraw (0,0) circle (1pt);
    %\addplot[domain=-3:3,thick,black,no marks,smooth,fill=blue!30] {x^2};
    \end{tikzpicture}
}
}
\quad \quad \quad \quad\quad\quad
	\subfigure[Strict complementary condition in \cite{nurminski2015single} fails: the set $\bigcup_{x\in F_A(-x^*)}N_A(x)$ still has nonempty interior; and sharpness condition holds.]{\makebox[6.5cm][c]{
% 	\begin{tikzpicture}[scale = 0.7]
% 	\fill[blue!30,opacity=1] (-2.1,3)--(-1,0)--(1,0)--(2.1,3)--(-2.1,3);
% 	\draw[blue,thick] (-2.1,3)--(-1,0)--(1,0)--(2.1,3);
% 	\fill[gray!30!red!20!blue!10!,opacity=1] (-4.25,-1.5)--(-1,0)--(1,0)--(4.25,-1.5)--(-4.25,-1.5);
% 	\draw[->,thick,gray](-1,0)--(-4.25,-1.5);
% 	\draw[->,gray,thick](1,0)--(4.25,-1.5);
% 	\draw[->,red!60](0,0)--(0,-1.5);
% 	\draw (0,2) node {$A$};
% 	\draw[red] (0.4,-1.2) node {$-x^*$};
% 	\draw (0.3,-0.3) node {$\bx$};
% 	\draw (-1,-0.8) node {$N_A(\bx)$};
% 	\draw[blue,thick] (-1,0)--(1,0);
% 	\draw (0,0.3) node {$F_A(-x^*)$};
% 	\filldraw (0,0) circle (1pt);
%     %\addplot[domain=-3:3,thick,black,no marks,smooth,fill=blue!30] {x^2};
%     \end{tikzpicture}
	\begin{tikzpicture}[scale = 0.7]
	\fill[gray,opacity=0.1] (-2.1,3)--(-1,0)--(1,0)--(2.1,3)--(-2.1,3);
	\draw[black,opacity=0.5] (-2.1,3)--(-1,0)--(1,0)--(2.1,3);
	\fill[gray,opacity=0.3] (-4.25,-1.5)--(-1,0)--(1,0)--(4.25,-1.5)--(-4.25,-1.5);
	\draw[->,thick,gray](-1,0)--(-4.25,-1.5);
	\draw[->,gray,thick](1,0)--(4.25,-1.5);
	\draw[->,densely dotted,black](0,0)--(0,-1.5);
	\draw (0,2) node {$A$};
	\draw (0.4,-1.2) node {$-x^*$};
	\draw (0.3,-0.3) node {$\bx$};
	\draw (-1,-0.8) node {$N_A(\bx)$};
	\draw[black,thick] (-1,0)--(1,0);
	\draw (0,0.3) node {$F_A(-x^*)$};
	\filldraw (0,0) circle (1pt);
    %\addplot[domain=-3:3,thick,black,no marks,smooth,fill=blue!30] {x^2};
    \end{tikzpicture}
	}
}
	\caption{Comparison between sharpness condition and strict complementarity condition.}	\label{F1.1}
\end{figure*}
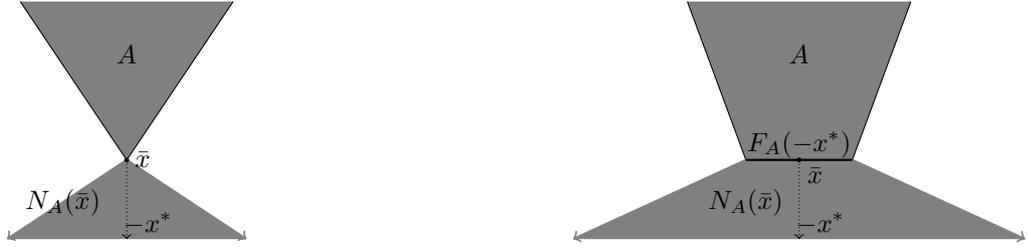

By contrast, for the sharpness condition to hold neither condition (i) nor condition (ii) above are required. Instead, Proposition \ref{PP3} shows that sharpness simply requires $-x^*$ to belong to the interior of the union of all normal vectors at every optimal solution (see Figure~\ref{F1.1}(b), where Condition (i) does not hold). Thus, the sharpness property requires a much weaker version of condition (ii).  Furthermore, Proposition~\ref{active} proves that a polyhedron in an arbitrary Hilbert space is sharp with respect to every unit vector; and we will show later in Section~\ref{Sec:LP} that the sharpness property is enough to ensure that the minimization problem $\min_{x\in A}\ang{x^*,x}$ can be solved by SPP. 
%In particular, when the feasible set is bounded, the union of all normal cones at every point from the face defined by $-x^*$ must have nonempty interior, and $-x^*$ belongs to the interior set (see Figure~\ref{F1.1}).
\end{remark}

We next establish the connection between the sharpness property and the subdifferential operator of the Fenchel conjugate of the indicator function at $x^*$. Recall that, for any proper, lower-semicontinuous convex function $f$, the \emph{Fenchel--Young characterization} of the subdifferential is given as follows:
\begin{equation}\label{eq:FY1}
x^*\in \partial f(x) \Longleftrightarrow \ang{x^*,x} = f(x)+f^*(x^*) \Longleftrightarrow x\in \partial f^*(x^*).
\end{equation}
In particular, 
\begin{equation}\label{eq:FY2}
x^*\in \partial f(x)\Longrightarrow x^*\in \dom f^*\text{~and~}x\in \dom f.
\end{equation}

\begin{remark}\label{rem:sigma} Given a convex and closed set $A$, the Fenchel conjugate of $f(\cdot):=\mathbbm{1}_A(\cdot)$, denoted by $\sigma_A(\cdot)$, is called the {\em support function of $A$}. Recall that $\partial f=N_A$ and $\dom f=A$. Applying \eqref{eq:FY1} to this $f$ and using the definitions, we obtain
\begin{equation}\label{def:supp}
  f^*(v)= \sigma_A(v)= \sup_{y\in A} \ang{v,y}.
\end{equation}
Using \eqref{eq:FY1} for this $f$,  we obtain 
\begin{equation}\label{eq:sa}
v\in N_A(x)\Longleftrightarrow x\in A\text{~and~}\ang{v,x}= \sigma_A(v)\Longleftrightarrow x\in \partial \sigma_A(v) \Longleftrightarrow  \partial \sigma_A(v)= \Argmax_A \ang{v,\cdot}=F_A(v),
\end{equation}
where we also used \eqref{def:supp} and Fact \ref{fact 1} in the rightmost equivalence. The equivalence above also shows that $D(\partial \sigma_A)={\text{R}(N_A)}$ (see Definition~\ref{def:pts}(a)(b)). These facts will be used in the next result.
\end{remark}

\begin{proposition}
\label{Prop:P1}
Let $A\subset \mathcal H$ be a closed and convex set. Fix $x^*\in \mathcal S\cap \text{R}(N_A)$ and $\al \in (0,1]$. 
Consider the following statements:
\begin{enumerate}[(i)]
    \item The set $A$ is $\al$-sharp w.r.t.\ $x^*$.
    \item $\emptyset\neq \partial \sigma_A(v)\subset \partial \sigma_A(x^*)$ for all $v\in x^*+\alpha (\mathcal B\!\setminus{\mathcal S})$.
\end{enumerate}Then (ii) implies (i). If $A$ is bounded, then (i) implies (ii). In the latter situation, $\partial \sigma_A(v)\subset F_A(x^*)$ for all $v\in x^*+\alpha (\mathcal B\!\setminus{\mathcal S})$.
\end{proposition}

\begin{proof}
% Since $A$ is convex and closed, we have $(\mathbbm{1}_A)^{**}(x)=\mathbbm{1}_A(x)$. The Fenchel-Young equality gives
% \[
% x^*\in \partial f(x) \Longleftrightarrow x\in \dom f, x^*\in \dom f^* \AND \ang{x^*,x} = f(x)+f^*(x^*).
% \]
For simplicity, denote $B(\alpha,x^*):=x^*+\alpha (\mathcal B\!\setminus{\mathcal S})$. Note also that the last statement will follow directly from (ii) and the equivalence \eqref{eq:sa} for $v:=x^*$. Indeed, the rightmost expression in \eqref{eq:sa} with $v:=x^*$ gives $F_A(x^*)=\partial \sigma_A(x^*)$.

Now, assume that (ii) holds. We will use Proposition~\ref{PP3} to show (i). More precisely, we will show that $B(\alpha,x^*)\subset \bigcup_{x\in F_A(x^*)}N_A(x)$ holds. Then (i) will follow from the part (ii) implies (i) in Proposition~\ref{PP3}. Indeed, take $v\in B(\alpha,x^*)$ we need to show that there exists $x\in F_A(x^*)$ such that $v\in N_A(x)$.  By (ii) we have $\partial \sigma_A(v)\neq \emptyset$ and $\partial \sigma_A(v)\subset \partial \sigma_A(x^*)$. Hence, we can write
$$
x\in \partial \sigma_A(v) \Longrightarrow  x\in \partial \sigma_A(x^*)=F_A(x^*).
$$
On the other hand, by \eqref{eq:sa} we have
$$
x\in \partial \sigma_A(v)\Longleftrightarrow x\in \Argmax_{y\in A}\ang{v,y}  \Longleftrightarrow v \in N_A(x).
$$
Combining the rightmost parts of the two expressions above gives that
% \[
% \hbox{If }v\in B(\alpha,x^*) \AND  x\in \partial \sigma_A(v), \quad \hbox{then}\quad x\in F_A(x^*) \AND v \in N_A(x).
% \]
\[
v\in B(\alpha,x^*)\text{~and~}x\in \partial \sigma_A(v)\implies x\in F_A(x^*)\text{~and~}v \in N_A(x).
\]
Therefore, for every $v\in B(\alpha,x^*)$, there exists $x\in F_A(x^*)$ s.t. $v \in N_A(x)$. In other words,
\[B(\alpha,x^*)\subset \bigcup_{x\in F_A(x^*)} N_A(x).
\]
Proposition~\ref{PP3} (part (ii) implies (i)), now implies that $A$ is $\al$-sharp w.r.t.\ $x^*$.

Next, assume now that $A$ is bounded. We will prove that (i) implies (ii).  Namely, we will use the boundedness and sharpness of $A$ and Proposition~\ref{PP3} to show that (ii) holds. 
Fix $v\in B(\alpha,x^*)$. The nonemptiness of $\partial \sigma_A(v)$ follows from the fact that the set $A$ is bounded and $F_A(v)=\partial \sigma_A(v)$ (from the rightmost expression in \eqref{eq:sa} for $v$).
%By assumption (i), Proposition~\ref{PP3}, we have
%\begin{equation}\label{eq:2ia}
% v\in B(\alpha,x^*)\subset \bigcup_{x\in F_A(x^*)} N_A(x).
%\end{equation}
%Hence, there exists $x\in F_A(x^*)$ such that $v\in N_A(x)$. By \eqref{eq:sa} this implies that $x\in \partial \sigma_A(v) = F_A(v)$.
%for every  $x\in \partial \sigma_A(x^*)$ we have $x\in \partial \sigma_A(v)$. Therefore, for every $v\in B(\alpha,x^*)$ we have $\partial \sigma_A(x^*)\subset \partial \sigma_A(v)$. We prove now the opposite inclusion for $v\in B(\alpha,x^*)$. 
 To show that $\partial \sigma_A(v) \subset \partial \sigma_A(x^*) $ we need to show that $F_A(v)\subset F_A(x^*)$. Suppose to the contrary that $F_A(v)\not\subset F_A(x^*)$. Hence, there is $y\in F_A(v)\setminus F_A(x^*)$. Using Fact \ref{fact 1}, the latter means that $v\in N_A(y)$ and $x^*\not\in N_A(y)$. Since (i) holds, we can use Proposition \ref{PP01} to write
        $$
  \norm{x^* - v}\ge       d(x^*,N_A(y)) \ge \al,
        $$
which contradicts the fact that $v\in B(\alpha,x^*)$. Therefore, we must have $F_A(v)\subset F_A(x^*)$, and hence (ii) holds. %Consequently, (ii) holds true.
% Let $F_A(x^*)$ be the face with respect to $x^*$. Then, by Fact~\ref{fact 1},
%  $\partial f^*(x^*)=F_A(x^*)$.
% %From Proposition~\ref{PP3}, we have $B_\al(x^*)\subset \bigcup_{x\in F_A(x^*)} N_A(x^*)$.
%       % Because $A$ is bounded, by Proposition~\ref{PP3}, $B_\al(x^*)\subset \cl\bigcup_{x\in F_A(x^*)}N_A(x) = \cl\bigcup_{x\in F_A(x^*)}\partial \mathbbm{1}_A(x)$. 
%       For any $v\in B_\al(x^*)$, and $u\in \partial f^*(v)$ we have $u\in \argmax_{A}\ang{v,\cdot}$, or $u\in F_A(v)$.
%Conversely, suppose statement (ii) holds. 
\end{proof}

\begin{remark}
\textcolor{black}{
We say that a proper lower semicontinuous convex function $f:\mathcal{H}\rightrightarrows \R_\infty$ is {\em quasi-polyhedral at $\bx\in \dom f$} (see \cite{Cha}) if there exists \textcolor{red}{$r>0$ such that $\partial f(x)\subset f(\bx)$} for all $x\in \bx + r\mathcal{B}$. Furthermore, from \cite[Proposition 3.4]{Cha}, if $f$ is continuous at $\bx$, then $f$ is quasi-polyhedral at $\bx$ if and only if $f$ is {\em conical} at $\bx$; meaning that there is $r>0$ and a sublinear function $p: \mathcal{H} \to \R_\infty$ such that $f(x) = f(\bx)+ p(x-\bx)$ for every $x\in  \bx + r\mathcal{B}$. Hence, Proposition \ref{Prop:P1} shows that $A$ is sharp w.r.t $x^*$ if and only if the function $\sigma_A$ is quasi-polyhedral at $x^*$; and if $\sigma_A$ is continuous at $x^*$, then $\sigma_A$ is also conical.}
%Note that the choice of $V$ clearly depends on the prescribed $W$, and in general $V$ will be smaller when $W$ is smaller. Maximally monotone maps are always USC in the interior of their domains (see \cite[Corollary 4.2.13]{BurachikIusemBook}). Being a subdifferential, the map $\partial \sigma_A(\cdot)$ is maximally monotone, and hence it is USC in the interior of its domain.  Since $D(\partial \sigma_A)={\rm R}(N_A)$, when $A$ is bounded and $\al$-sharp w.r.t vector $x^*$, \textcolor{blue}{then} $x^*\in {\rm int}({\text{R}(N_A)})={\rm int}(D(\partial \sigma_A))$. Hence, the sharpness property at $x^*$ and Proposition~\ref{PP3} imply that $\partial \sigma_A$ is USC at $x^*$. Moreover, Proposition \ref{Prop:P1} shows that $\partial \sigma_A$ is {\em uniformly} USC at $x^*$, in the sense that, no matter what $W$ we choose in \eqref{eq:T}, the required set $V=B(\alpha,x^*)$ can be always the same (it is independent of $W$).
\end{remark}
\subsection{Metric characterizations}
\label{metric}

The sharpness property has a strong connection with regularity-type properties of sets. In particular, we will show in this section that sharpness with respect to a vector $x^*$ is equivalent to subtransversality
%\textcolor{green}{EN: 18/6 Shall we put reference here ?}
%\textcolor{red}{Yes, Hoa can you please add a ref for the concept of subtransversality?}
 (also known as metric subregularity, see \cite[Definition~3.1]{li2018calculus} and \cite{KruThaoLuke}) between the set $A$ and its supporting hyperplane with respect to $x^*$.
Recall that, given two convex sets $A,B$ such that $A\cap B\not=\emptyset$, the pair $\{A,B\}$ is {\em subtransversal} if there is $\al \in (0,1)$ such that
$$
\al d(x,A\cap B)\le \max\{d(x,A),d(x,B)\},\quad \forall x\in \mathcal H.
$$

\begin{remark}\label{rem:subtrans}
Note that subtransversality for the sets $\{A,B\}$ is equivalent to the property
\begin{equation}
    \label{eq:ST}
    \al d(x,A\cap B)\le d(x,A),\quad \forall x\in B.
\end{equation}
Indeed, this follows from the fact that $d(x,A)\le \max\{d(x,A),d(x,B)\} $ when $x\not\in B$. The connection between subtransversality and sharpness arises when we specialize \eqref{eq:ST} for the pair $\{A,F\}$ with 
\begin{equation}\label{eq:F}
    F:=\{x\in \mathcal H\::\: \ang{x^*,x}=\sup_{y\in A} \ang{x^*,y}\}.
\end{equation}
In this situation, using the definition of $F_A(x^*)$, \eqref{eq:ST} becomes
\[
 \al d(x,F_A(x^*))\le d(x,A),\quad \forall x\in F.
\]
The next result shows that subtransversality of $\{A,F\}$ is equivalent to the sharpness of $A$.
\end{remark}

%The pair $\{A,B\}$ is called substransversal if there is $\al \in (0,1)$ such that $\{A,B\}$ is $\al$-subtransversal. \textcolor{magenta}{MKT 27/06: Some repetition here that needs to be fixed. Not sure if the second sentence should be different?} 
%We first provide a metric characterization for the sharpness property.
% \textcolor{green}{EN:18/6 $d(x, \emptyset) = \infty ?$ probably missed it somewhere.}
% \textcolor{red}{I assume the definition only holds when $A\cap B\not=\emptyset$. I will add it to the definition.}
\begin{theorem}
\label{P1}
Suppose $A$ is a closed convex set of a Hilbert space $\mathcal H$, $x^*\in \mathcal S$, {$F_A(x^*)\neq \emptyset$}, and $\al\in (0,1)$. Let $F$ be  as in \eqref{eq:F}. Assume that
\begin{equation}
    \label{P1.1}
    \al \,d(x,F_A(x^*))\le d(x,A),\quad\forall x\in F,
\end{equation}
and define $\gamma:= \al \sqrt{1-\tfrac{1}{4}\al^2}$. Then, the set $A$ is $\gamma$-sharp w.r.t. vector $x^*$. Conversely, define $\beta:=\textcolor{black}{2\al/(1-\al)}$. If $A$ is $\beta$-sharp w.r.t. $x^*$, then inequality \eqref{P1.1} holds.
\end{theorem}
 
%\textcolor{green}{EN 23/06 Here we define  $F_A(x^*)=A\cap F$ with $F:=\{x\in \mathcal H\::\: \ang{x^*,x}=\sup_{y\in A} \ang{x^*,y}\}$,
%and $A\cap F$ is used many times later on. Shell we replace them with $F_A(x^*)$ ? }
%\textcolor{blue}{HB: I replaced them}
\begin{proof}
First, assume that inequality \eqref{P1.1} holds. We will show that the set $A$ is $\gamma$-sharp w.r.t.\ $x^*$, where $\gamma:= \al \sqrt{1-\tfrac{1}{4}\al^2}$. If $A\setminus F_A(x^*)$ is empty, then $F_A(x^*)=A$. 
%\textcolor{green}{EN 01/07 That is $F_A(x^*) = A$, which is easier to digest: if $a-b = 0$ then $a=b$ :)}
In this case, by Remark \ref{rem:sharp 1}, we deduce that $A$ is trivially $\gamma$- sharp w.r.t. $x^*$ for any $\gamma>0$. So it is enough to assume that $A\setminus F_A(x^*)$ is not empty. By Fact~\ref{fact 1}, in this case there is $y\in A$ such that $x^*\notin N_A(y)$. Suppose that $N_A(y)\neq \{0\}$, since otherwise it always holds that $d(x^*, N_A(y)) = \norm{x^*} = 1 \ge \al$ for all $\al \in (0,1)$. So $N_A(y)\neq \{0\}$ and hence we can take $y^*\in \mathcal S\cap N_A(y)$, i.e.,  $\norm{y^*}=1$. We consider two cases. We claim that both cases lead to
\begin{equation}
    \label{eq:claim 1}
    d(x^*,\cone[y^*]) \ge \al\sqrt{1-\tfrac{\al^2}{4}}.
\end{equation}
{\bf Case 1.} Suppose that $\ang{x^*,y^*} \le 0$. Then, by Fact~\ref{L2}, $d(x^*,\cone[y^*])=1\ge \alpha$. Since $\al$ is always larger than the rightmost expression in \eqref{eq:claim 1}, in this case \eqref{eq:claim 1} holds. Now we show that \eqref{eq:claim 1} also holds for the next case.
{\bf Case 2.} Suppose that $\ang{x^*,y^*} > 0$. We first show that there is an $x\in F$ such that $x-y \in \cone[y^*]$ (see Figure~\ref{P3.12.F1}), and so $y$ is a projection of $x$ onto $A$. This is equivalent to the existence of $t\ge 0$ and $x\in F $ such that $x-y = ty^*$, or
    $$
    x=y+ty^*\in F .
    $$
    By the definition of $F $, the equality above is equivalent to stating that the equation $\ang{x^*,y+ty^*} = \sup_{A}\ang{x^*,\cdot}$ has a solution $t\ge 0$.
    Because $\ang{x^*,y^*}> 0$, and $y\in A$, the constant 
    $$t := \frac{\sup_{w\in A}\ang{x^*,w} - \ang{x^*,y}}{\ang{x^*,y^*}}$$ 
    is nonnegative and satisfies $\ang{x^*,y+ty^*} = \sup_{w\in A}\ang{x^*,w}$. 
Therefore, $x := y+ty^* \in F $, and $y=P_A(x)$. Hence, $x-y\in \cone[y^*]\subset N_A(y)$. Taking to account that $\norm{y^*} = 1$, we have $y^*= \frac{x-y}{\norm{x-y}}$.
    \begin{figure}[h]
\begin{center}
\begin{tikzpicture}[scale=0.6]
\draw[black, thick] (-3,0)--(4,0); 
\filldraw[gray,opacity=0.1] (-3,3)--(-3,0) -- (0,0) -- (3,3);
\draw[gray] (-3,3)--(-3,0) -- (0,0) -- (3,3)--(-3,3);
\draw (-1,1.5) node {$A$};
\filldraw[black] (2,2) circle (2pt);
\filldraw[black] (2.3,2) node {$y$};
\filldraw[black] (3.7,0) circle (2pt);
\filldraw[black] (3.7,-0.5) node {$x$};
\draw[dashed] (3.7,0)--(2,2);
\filldraw (0,0) circle (2pt);
\filldraw[black] (0,-0.5) node {$z$};
\filldraw (4.5,0) node {$F $};
\end{tikzpicture}
\end{center}
\caption{There is $x\in F $ such that $y\in P_A(x)$ and $y-x\in \cone[y^*]$.}\label{P3.12.F1}
\end{figure}
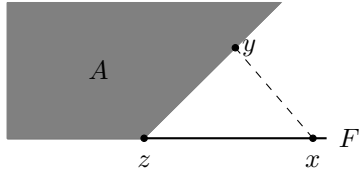

Let $z:= P_{F_A(x^*)}(x)$ be a projection of $x$ onto $F_A(x^*)$. Since $z\in A$ and $x-y\in N_A(y)$, the following inequality holds
    \begin{align}\label{P1P1}
        \ang{x-y,z-y} \le 0.
    \end{align}
    Additionally, inequality \eqref{P1.1} implies
    $$
    \al\norm{x-z}=\al d(x,F_A(x^*) ) = \al d(x,A\cap F )  \le d(x,A) = \norm{x-y}.
    $$
    Combine the expression above with \eqref{P1P1} to obtain 
    \begin{align*}
        \norm{x-z}^2 &=\norm{x-y}^2+\norm{y-z}^2+2\ang{x-y,y-z}\\
        &\ge \norm{x-y}^2+\norm{y-z}^2\\
        &\ge \al^2\norm{x-z}^2 +\norm{y-z}^2.
    \end{align*}
    Therefore, $(1-\al^2)\norm{x-z}^2 \ge \norm{y-z}^2.$ On the other hand,
    \begin{align*}
        \ang{x-y,x-z} &= \frac{1}{2}\left(\norm{y-x}^2 + \norm{x-z}^2-\norm{y-z}^2\right)\\
        &\ge \frac{1}{2}\left(\norm{y-x}^2 + \norm{x-z}^2-(1-\al^2)\norm{x-z}^2\right)\\
        &= \frac{1}{2}\left(\al^2\norm{x-z}^2 + \norm{x-y}^2\right).
    \end{align*}
    From $y^*= \frac{x-y}{\norm{x-y}}$, the inequality above, and Cauchy-Schwartz inequality, we obtain
    $$\ang{y^*,\frac{x-z}{\norm{x-z}}}= \ang{\frac{x-y}{\norm{x-y}},\frac{x-z}{\norm{x-z}}}\ge \frac{\al^2}{2} \frac{\norm{x-z}}{\norm{x-y}}+\frac{1}{2}\frac{\norm{x-y}}{\norm{x-z}}\ge \al,$$
    where the last inequality follows from the fact that the function $\eta(t):=\frac{\al^2\, t}{2} +\frac{1}{2t}$ over $\R_{++}$ attains a minimum at $t^*:=1/\al$ and its minimum value is $\eta(t^*)=\al$.
On the other hand, since $x,z\in F $, by the definition of $F $ we must have $\ang{x^*,\frac{x-z}{\norm{x-z}}} = 0$. Combining the latter equality with the last inequality yields
    $$
    \ang{y^*-x^*,\frac{x-z}{\norm{x-z}}}\ge \al,
    $$
    which, by the definition of dual norm, implies $\norm{y^*-x^*} \ge \al$. The latter inequality and Fact~\ref{L2} give
    \begin{align*}
         d(x^*,\cone[y^*])^2= 1 -\ang{y^*,x^*}^2 &= 1- \left(\frac{\norm{x^*}^2+\norm{y^*}^2-\norm{x^*-y^*}^2}{2} \right)^2\\
         &\ge 1- \left(\frac{2-\al^2}{2} \right)^2 = \frac{1}{4}\al^2(4-\al^2),
    \end{align*}
    where we also used the fact that $\|x^*\|=\|y^*\|=1$. Altogether, in both cases 1 and 2, we always have $d(x^*,\cone[y^*]) \ge \al\sqrt{1-\tfrac{\al^2}{4}}$, and because $y^*\in \mathcal S\cap N_A(y)$ is chosen arbitrarily, we have
    $$
    d(x^*, N_A(y)) = \inf_{y^*\in \mathcal S\cap N_A(y)} d(x^*,\cone[y^*]) \ge \al\sqrt{1-\tfrac{\al^2}{4}} =\gamma.
    $$
    Consequently, we have proved that the set $A$ is $\gamma$-sharp w.r.t. vector $x^*$, where $\gamma:=\al\sqrt{1-\tfrac{\al^2}{4}}>0$.
    To prove the converse implication, we assume by contradiction that the set $A$ is $\beta$-sharp w.r.t. $x^*$, where $\be := 2\al/(1-\al)$ and there is $\bar x\in F $ such that 
    \begin{equation*}
        \al d(\bar x,F_A(x^*)) > d(\bar x,A).
    \end{equation*}
    The strict inequality implies that there exists \textcolor{black}{$\hat \al \in (0,\al)$} such that
     \begin{equation}\label{P1.2}
       \textcolor{black}{\al d(\bar x,F_A(x^*))>\hat\al d(\bar x,F_A(x^*)) > d(\bar x,A).}
    \end{equation}
    Consider the closed half-space $F^\prime:= \{x\in \mathcal H:\; \ang{x^*,x} \ge \sup_{y\in A}\ang{x^*,y} \}$. Observe that, $F =\bd F^\prime$ and 
    $$
   F_A(x^*)= A\cap F  = A\cap F^\prime,
    $$
    and for any $x \in F $, we have $N_{F^\prime}(x) =\cone[-x^*]$. Set $\de:=d(\bar x,F_A(x^*))$. Since $\bar x\in F$, \eqref{P1.2} implies that $\delta>0$ and $\bar x\not\in  A$. Let $\bar y:= P_A(\bar x)$, so $d(\bar x,A)=\|\bar x- \bar y\|>0$. { With this notation, inequality \eqref{P1.2} becomes 
    \begin{equation}
    \label{eq:al}
    \al\de >\hat{\al}\de> \norm{\bar x- \bar y}.
    \end{equation}
We claim that \eqref{eq:al} implies that $\bar y \notin F_A(x^*)$.  Indeed, if $\bar y \in F_A(x^*)$ then the definition of $\de$ yields
  \[
  \de=d(\bar x,F_A(x^*))\le \norm{\bar x- \bar y}<\al\de<\delta,
  \]
  a contradiction. Therefore, our claim holds and $\bar y \notin F_A(x^*)$. Take now $\widetilde{y}:=P_{F_A(x^*)}(\bar y)$.   By triangle inequality,}
  \textcolor{black}{
    $$d(\bar y,F_A(x^*)) = \norm{\bar y-\widetilde{y}}\ge \norm{\widetilde{y}-\bar x} - \norm{\bar x-\bar y} \ge d(\bx,F_A(x^*))- \norm{\bar x - \bar y} >(1-\al)\de,$$
where we used the definition of $\widetilde{y}$ in the first equality, the fact that $\widetilde{y}\in F_A(x^*)$ in the second inequality, and \eqref{eq:al} in the last inequality.
Set $\hat \delta := (1-\al)\de$. Using the above expression, and the definitions of $F'$ and $\hat \delta$ we derive
    \begin{equation}
        \label{eq:yb}
     d(\bar y,F_A(x^*))   =d(\bar y,A\cap F')>\hat \delta.
    \end{equation}
With this notation, $\hat \al\de=\frac{\hat \al \hat\de}{(1-\al)}$ and the second inequality in \eqref{eq:al} rewrites as
    \begin{equation}
        \label{P1.2'}
        \textcolor{black}{\frac{\hat \al \hat\de}{1-\al}}> \norm{\bx-\by}.
    \end{equation}}
    Define the function $f:\mathcal H^2\to \R_{\infty}$ as
    \begin{equation}
        \label{def:f}
        f(x,y):=\norm{x-y}+\mathbbm{1}_{F^\prime}(x)+\mathbbm{1}_{A}(y),\quad (x,y)\in \mathcal H^2.
    \end{equation}
    We now use \eqref{P1.2'} and the Ekeland Variational Principle, Lemma \ref{EVP}, for $X:=\mathcal H^2$ equipped with the max norm  $\|(x,y)\|:=\max\{\|x\|,\|y\|\}$, and $\psi:=f$.  As mentioned in \eqref{def:dual norm}, the corresponding dual norm is $\|(x,y)\|_{*}=\|x\|+\|y\|$. Recall that $(\bar x, \bar y)\in F'\times A$, so $f(\bar x, \bar y)=\norm{\bar x -\bar y}$. This fact, combined with \eqref{P1.2'} and the definition of $f$ gives
    $$
     \inf_{\mathcal H\times \mathcal H} f(x,y)\ge 0 > f(\bar x,\bar y)-\textcolor{black}{\frac{\hat \al}{1-\al} \hat\de}.
    $$
   Hence we are in conditions of the Ekeland's variational principle with \textcolor{black}{$\varepsilon:=\frac{{\hat\al}\hat\de}{(1-\al)}$} and $\bar w:=(\bar x, \bar y)$. We apply the principle for the choice \textcolor{black}{$\lambda:=\hat\delta$}, for which there exists $(\hat x,\hat y)$ such that
    \begin{enumerate}[(i)]
        \item $\norm{\hat x -\bar x} <\textcolor{black}{\hat \de}$, $\norm{\hat y - \bar y} <\textcolor{black}{\hat \de}$;
        \item $f(\hat x,\hat y)=\norm{\hat x - \hat y} \le f(\bar x, \bar y)=\norm{\bar x - \bar y}$;
           \end{enumerate} 
    Moreover, our choices imply that \textcolor{black}{$\varepsilon/\lambda={\hat \alpha}/(1-\alpha)$}. Altogether, condition (iii) in Lemma \ref{EVP} with the max norm implies that for every $\mathcal H^2\ni (x,y)\not= (\hat x,\hat y)$ we have
        \begin{equation}\label{eq:h}
            f(\hat x,\hat y)< f(x,y) +\textcolor{black}{\frac{\hat\alpha}{(1-\alpha)}}\max\left\{\norm{x-\hat x},\norm{y-\hat  y}\right\}=:h(x,y). 
        \end{equation}
The above inequality implies that $(\hat x,\hat y)\in F'\times A$.  Consequently, the following statements hold. 
\begin{itemize}
    \item[(I)] By \eqref{eq:yb} we have that $\textcolor{black}{d(\by,A\cap F^\prime)>\hat \de}$. Using also (i)  and the fact that $d(\bx,A\cap F^\prime)=\textcolor{black}{\de > \hat \de}$,  we have $\hat x,\hat y\notin A\cap F^\prime$. Indeed, if we would have $\hat x\in A\cap F^\prime$ then this would imply
    \[
 \hat \de<   d(\bx,A\cap F^\prime)\le \norm{\bx-\hat x},
    \]
contradicting (i). A similar expression, mutatis mutandis, can be used to show that $\hat y\notin A\cap F^\prime$.
    \item[(II)] Because $\hat y\in A$, $\hat x\in F^\prime$, by (I) we must have $\hat x\neq \hat y$. The latter fact and (ii) yield $f(\hat x,\hat y) >0$.
    \item[(III)] By (II) and Remark \ref{rem:norm} we can write the subdifferential of $f$ as the sum of the subdifferentials at $(\hat x,\hat y)$
    $$
    \partial f(\hat x,\hat y) = \partial_1 f(\hat x,\hat y) \times \partial_2 f(\hat x,\hat y)=
    \left(\frac{\hat x-\hat y}{\norm{\hat x-\hat y}} +N_{F'}(\hat x),\frac{\hat y-\hat x}{\norm{\hat x-\hat y}} + N_A(\hat y)\right),
    $$
where $\partial_1$ and $\partial_2$ stand for the partial subdifferentials w.r.t. the first and second variables, respectively. Note that each of these partial subdifferentials exist because of the continuity of the first term in \eqref{def:f}.  
 \item[(IV)] By \eqref{eq:h}, it is clear that $(\hat x,\hat y)$ is the global minimizer of $h$.
 \end{itemize}
 By (IV), we have that $ 0\in \partial h(\hat x,\hat y)$. To compute the subdifferential of $h$, we note that the second term in \eqref{eq:h} is continuous everywhere. Hence, the subdifferential of $h$ can be expressed as the sum of the subdifferential of $f$ plus the subdifferential of the second term in \eqref{eq:h}. To write down the inclusion  $ 0\in \partial h(\hat x,\hat y)$, we use Remarks \ref{rem:norm} and \ref{rem:max} to write:
 % By the convex sum rule applied to the function $(f+\textcolor{red}{{\alpha}/(1-\alpha)}\max\{\norm{x-\hat x},\norm{y-\hat y}\})$  and the sum indicator functions $\mathbbm{1}_{F^\prime}(x)+\mathbbm{1}_{A}(y)$ at $(\hat x,\hat y)$ first; and applied to $f(x,y)$ and  $\textcolor{red}{{\alpha}/(1-\alpha)}\max\{\norm{x-\hat x},\norm{y-\hat y}\}$ at $\bar x$ second, and note that the subdifferential of  $\textcolor{red}{{\alpha}/(1-\alpha)}\max\{\norm{x-\hat x},\norm{y-\hat y}\}$ at $(\hat x,\hat y)$ is a subset of $\textcolor{red}{{\alpha}/(1-\alpha)}\B^*_{\mathcal H\times \mathcal H}$, and subdifferentials of $\mathbbm{1}_{F^\prime}(x)+\mathbbm{1}_{A}(y)$ at $(\hat x,\hat y)$ is $N_{F^\prime}(\hat x)\times N_A(\hat y)$, we have
    $$
    0\in \left(\frac{\hat x-\hat y}{\norm{\hat x-\hat y}} + N_{F^\prime}(\hat x),\frac{\hat y-\hat x}{\norm{\hat x-\hat y}}+N_A(\hat y)\right) +\textcolor{black}{\frac{\hat\alpha}{(1-\alpha)}} \B^*_{\mathcal H^2}, 
    $$
    where the ball in the rightmost term is the one induced by the sum norm as given in \eqref{eq:L1}. The inclusion above and the definition of the dual ball, give
    $$
    d\left(\frac{\hat y-\hat x}{\norm{\hat x-\hat y}},N_{F^\prime}(\hat x)\right)+d\left(\frac{\hat x-\hat y}{\norm{\hat x-\hat y}},N_{A}(\hat y)\right) = d\left(-\left(\frac{\hat x-\hat y}{\norm{\hat x-\hat y}},\frac{\hat y-\hat x}{\norm{\hat x-\hat y}} \right),N_{F^\prime}(\hat x)\times N_{A}(\hat y) \right) \le \textcolor{black}{{\hat\alpha}/(1-\alpha)},
    $$
    where we are using the definition of the sum ball (see \eqref{eq:L1}) in the first equality. Note also that $\hat x\in F^\prime = \{x\in \mathcal H:\; \ang{x^*,x} \ge \sup_{y\in A}\ang{x^*,y} \}$, so 
    $$ N_{F^\prime}(\hat x) = \begin{cases}\cone[-x^*] &\text{ if } \hat x\in F = \bd F^\prime,\\
    \{0\}& \text{ otherwise}.\end{cases}$$ 
    This gives, $N_{F^\prime}(\hat x)\subset \cone[-x^*]$.  Therefore,
    \begin{align}
 \notag   &d\left(\frac{\hat x-\hat y}{\norm{\hat x-\hat y}},N_{A}(\hat y)\right)+d\left(\frac{\hat x-\hat y}{\norm{\hat x-\hat y}},\cone[x^*]\right)\\
 \label{Eq:T2.21P2}   &\le d\left(\frac{\hat x-\hat y}{\norm{\hat x-\hat y}},N_{A}(\hat y)\right)+d\left(\frac{\hat y-\hat x}{\norm{\hat x-\hat y}},N_{F^\prime}(\hat x)\right) \le \textcolor{black}{{\hat\alpha}/(1-\alpha)}.
    \end{align}
To simplify notation, call $\omega_0:=\frac{\hat x-\hat y}{\norm{\hat x-\hat y}}$. By triangle inequality, we have
\begin{equation}
    \label{Eq:T2.21P1}
     d(x^*,N_A(\hat y)) \le \norm{x^*-\omega_0}+d\left(\omega_0,N_{A}(\hat y)\right) .
\end{equation}
We claim that
\begin{equation}\label{eq:claim}
    d\left(\omega_0,\cone[x^*]\right)\ge \frac{1}{2}\norm{x^*-\omega_0}.
\end{equation}
Indeed, from Fact~\ref{L2}, we have
\(
d\left(\omega_0,\cone[x^*]\right)^2 = 1-\max\left(0,\ang{x^*,\omega_0}\right)^2.
\)
We consider two cases:
\begin{enumerate}
    \item If $\ang{x^*,\omega_0}\le 0$, then $d\left(\omega_0,\cone[x^*]\right)=1$ and by triangle inequality
   \[
   \norm{\omega_0-x^*}\le \norm{\omega_0}+\norm{x^*}= 2= 2\, d\left(\omega_0,\cone[x^*]\right),
   \]
   as claimed in \eqref{eq:claim}.
    \item If $\ang{x^*,\omega_0}> 0$, then
    \begin{align*}
    d\left(\omega_0,\cone[x^*]\right)^2 &= 1- \ang{x^*,\omega_0}^2 \\
    &= \frac{1}{2}\left(1+ \ang{x^*,\omega_0}\right)\left(2-2 \ang{x^*,\omega_0}\right)\\
     &>  \frac{1}{2}\left(2-2 \ang{x^*,\omega_0}\right)\\
    &=  \frac{1}{2}\left(\norm{x^*}^2+\norm{\omega_0}^2-2 \ang{x^*,\omega_0}\right)\\
    &= \frac{1}{2} \norm{{\omega_0}-x^*}^2,
    \end{align*}
\end{enumerate}
where we used $(1+ \ang{x^*,\omega_0})> 1$ in the inequality, and the fact that $\norm{x^*}=\norm{\omega_0}=1$ in the third equality. The above expression yields
\[
  d\left(\omega_0,\cone[x^*]\right)> \frac{1}{\sqrt{2}} \norm{{\omega_0}-x^*}>\frac{1}{2}\norm{{\omega_0}-x^*}.
\]
In both cases, we proved that \eqref{eq:claim} holds. From \eqref{Eq:T2.21P1}, we have
\begin{equation}\label{eq:T2.21_3}
  \begin{array}{rcl}
d(x^*,N_A(\hat y))      & \le &  \norm{x^*-{\omega_0}}+d\left(\omega_0,N_{A}(\hat y)\right)\\
&&\\
     &\le  & 2 d\left(\omega_0,\cone[x^*]\right) + d\left(\omega_0,N_{A}(\hat y)\right)\\
     &&\\
     &=& \left[ d\left(\omega_0,\cone[x^*]\right) + d\left(\omega_0,N_{A}(\hat y)\right)\right] +d\left(\omega_0,\cone[x^*]\right)\\
     &&\\
     &\le & \textcolor{black}{\hat\al/(1-\al) +\hat\al/(1-\al)= 2\hat\al/(1-\al)}< 2\al/(1-\al),
 \end{array}
\end{equation}
where we used \eqref{Eq:T2.21P1} in the first inequality, and \eqref{eq:claim} in the second one. For the last inequality we used \eqref{Eq:T2.21P2} for the expression between square brackets, and for the remaining term we used the fact that disregarding the first term in \eqref{Eq:T2.21P2} implies that $d\left(\omega_0,\cone[x^*]\right)\le \hat\al/(1-\al)$.
From (I), we have $\hat y \notin A\cap F' = A\cap F=F_A(x^*)$, and so by Fact \ref{fact 1} it holds that $x^*\notin N_A(\hat y)$. By the  $\beta$-sharpness of $A$ w.r.t. $x^*$ and the assumption on $\beta$ we deduce that
\[
d(x^*,N_A(\hat y)) \textcolor{black}{\ge 2\al/(1-\al)=\beta},
\]
contradicting \eqref{eq:T2.21_3}. This implies that inequality \eqref{P1.1} holds and the proof is complete.
% and \eqref{DC} implies that $x^*\in N_A(\hat y)$. On the other hand, from  which yields a contradiction.
    %Because $A\cap H\neq \emptyset$, then $H\neq \emptyset$, so there is 
   % Take  {\color{red}[MKT: Why is the existence of $y$ ensured?]}. , since otherwise \eqref{DC} holds trivially {\color{red}[MKT: What is the qualification on $\alpha$? Because the ``trival'' reason is less clear if $\alpha$ is fixed. Indeed, if $N_A(y)=\{0\}$, then \eqref{DC} gives
   %  $$ d(x^*,N_A(x)) = \|x^*\| >0, $$
    %but the validity of the inequality ``$\|x^*\|\geq \alpha$'' depends on the value of $\alpha$
    %]}. 
    %  Take $y^*\in S\cap N_A(y)$,  so $\norm{y^*}=1$. {\color{red}[MKT: I'm confused by the notation here. In the equation before \eqref{DC}, $\widehat{N}_A$ is the set of unit normals. So why do we need ``$\setminus\{0\}$'' and ``$\|y^*\|=1$''?]}.
\end{proof}

The following result from \cite{Yang04} is a characterization of subtransversality. We will use this result to formally express the connection between subtransversality and sharpness. 

\begin{lemma}[Subtransversality {\cite[Theorem~3.1]{Yang04}}]
\label{L10}
Suppose X is a normed linear space, $A,B\subset X$ are closed convex sets and $A\cap B\neq \emptyset$. The pair $\{A,B\}$ is subtransversal if and only if there exists a number $\al\in (0,1)$ such that
\begin{equation}\label{L1.1}
\al d(x,A\cap B)\le d(x,A),\quad \forall x\in B\setminus A.
\end{equation}
\end{lemma}
Theorem~\ref{P1} and Lemma~\ref{L10} yield the following result.
\begin{corollary}
\label{C2}
Suppose $A$ is a closed convex set of a Hilbert space $\mathcal H$, and $x^*\in \mathcal S$. Then, $A$ is sharp w.r.t.\ $x^*$ if and only if the pair $\{A,F\}$ is subtransversal, where $F:=\{x\in \mathcal H \::\, \ang{x^*,x}=\sup_{A} \ang{x^*,\cdot}\}$. %Conversely, if $A$ is $\al$-sharp w.r.t. vector $x^*$ where $\al \in (0,1/3)$, then the pair $\{A,F\}$ is subtransversal.
\end{corollary}
\begin{proof}
Suppose $\{A,F\}$ is subtransversal. Then, from Lemma~\ref{L10}, we have that \eqref{L1.1} holds for some $\al \in (0,1)$. By Theorem~\ref{P1}, the set $A$ is $\al^\prime$-sharp w.r.t. $x^*$, where $\al^\prime := \al\sqrt{1-\tfrac{1}{4}\al^2}\in (0,1)$. Conversely, assume that  $A$ is $\al$-sharp w.r.t. $x^*$ for some $\al \in (0,1)$. Then, by Theorem~\ref{P1}, we have
    $$
    \al^\prime d(x,A\cap F) \le d(x,A),\quad\forall x\in F,
    $$
    with $\al^\prime := \frac{\al}{2+\al}\in (0,1)$. Note that 
    $$2\al^\prime/(1-\al^\prime)= 2\frac{\tfrac{\al}{2+\al}}{1-\tfrac{\al}{2+\al}}=2\frac{\al}{2+\al - \al} = \al.$$ Because $\al^\prime \in (0,1)$, Lemma~\ref{L10} yields that the pair $\{A,F\}$ is subtransversal.
\end{proof}
\section{Optimization problems under sharpness condition}\label{S2}
%\texcolor{red}{MKT: Perhaps we should have a section focusing on studying the regularity properties themselves, and then another section on their application to the project result?
%
%HB: Yes. It sounds great.}

We consider now constrained convex problems of the following type.
\begin{equation}\label{eq:min f_A}\tag{CP}
 \min_{x\in A} f(x),
 \end{equation}
where $f:\mathcal{H}\to \R_{\infty}$ is a proper, lsc convex function, and $A$ is a closed convex set. We provide in this section sufficient conditions under which Problem~\eqref{eq:min f_A} can be solved by SPP. As expected,  sharpness plays a crucial role in our analysis. Namely, under the sharpness assumption, if (i) $b\notin A$, and (ii) the difference $(\inf_{x\in A}f) - f(b)$ is sufficiently large, then $P_A(b)$  solves Problem~\eqref{eq:min f_A}. In such situation, instead of solving Problem~\eqref{eq:min f_A}, we can solve the (hopefully) simpler problem of finding $P_A(b)$. Before establishing the main results of this section, we first find an upper bound on the distance between a point and a set using normal cones.
%following result is the cornerstone of this section, providing an upper bound on the distance between a point and a set.
%we first answer the question: given a set $A$, and a subset $\bar A\subset A$, and a point $b\notin A$, which conditions ensure the projections of $b$ onto $A$ will be in the set $\bar A$, in other words, $P_{\bar A}(b)\subset P_A(b)$.
%We first provide a simple characterization of 
%\todo[inline]{Say something here...}
\subsection{Upper bound on the distance}
Given a set $A$ and a point $b\notin A$, how can we estimate the distance $d(b,A)$? We address this question next, our analysis holds in a general Banach space (not necessarily Hilbert). A main tool in our proof is again the Ekeland Variational Principle. In the result below, we denote by $B[\rho,b]$ the closed ball of radius $\rho$ and centre $b$, and by $B(\rho,b)$ the corresponding open ball .

\begin{theorem}
\label{sl}	
Consider a Banach space $X$, a nonempty closed convex set $A$, points $a\in A$, $b\notin A$ with $\rho:=\norm{a-b}$, and $\eps >0$. Then, 
%\todo{Is it the same $\delta$ in both parts? IF so, should we write the quantifiers differently?}
$d(b,A)\le \norm{a-b}-\eps$ if
there is $\delta>0$ such that
	\begin{equation}\label{sl1}
	\inf\left\{d\left(x^*,N_A(x)\right):\; \norm{x^*} = 1,\;\ang{x^*,b-x} =\norm{b-x},\; x\in B[\rho,b]\cap B(\delta,a) \cap A \right\}\ge\frac{\eps}{\de}.
	\end{equation}
\end{theorem}

\begin{proof}
    Assume \eqref{sl1} holds for some $\eps,\de >0$. For contradictory purposes, assume also that $\eps$ is such that $d(b,A)> \norm{a-b}-\eps$. Consider the function $f:X\to \R_{\infty}$ defined by $f(y) := \norm{y-b}+ \mathbbm{1}_A(y)=\varphi_b(y)+ \mathbbm{1}_A(y)$. Here $\varphi_b(y):=\|y-b\|$ as in Remark \ref{rem:norm}.  Then, the assumption on $\eps$ implies that
    $f(a)=\norm{a-b}< \inf_A f+\eps$. Take $0<\eps^\prime<\eps$, such that
    $$
    \inf_A f +\eps> \inf_A f +\eps^\prime > f(a).
    $$
    By Ekeland's Variational Principle (Lemma~\ref{EVP}) {applied to $\bar w:=a$, $\psi:=f$, $\eps:=\eps^\prime$ and $\la:=\delta$}, there exists a vector $\hat x\in A\cap B(a,\de)$ such that
	\begin{align}
	\label{T1.P1}
		f(\hat x) &\le f(a),\\
	\label{T1.P2}
		f(\hat x)    &\le f(y) +\tfrac{\eps^\prime}{\de}\norm{y-\hat x},\quad \forall y \in A.
	\end{align}
Due to \eqref{T1.P1} and the fact that $a\in A$, we have that $\hat x\in A$ and  $\norm{\hat x-b}\le \norm{a-b} =\rho$, or $\hat x\in A\cap B[\rho,b]$. 
	Define $h(y):= \tfrac{\eps^\prime}{\de}\norm{y-\hat x}$.
	By \eqref{T1.P2}, it follows that $\hat x$ is a global minimizer of the sum function $f+h$ {and hence the definition of subdifferential yields the inclusion}
	$$
	0 \in \partial (f+h)(\hat x)=\partial f(\hat x)+\partial h(\hat x)
	$$
	%Take $\epsilon>0$ such that $$\epsilon<\min\left\{\al - \al', \rho-\norm{x_0-b}, \de - \norm{x_0-a}\right\}.$$
	Note that the subdifferential sum formula can be used to differentiate $(f+h)$ because $h$ is continuous everywhere. By Remark \ref{rem:norm}, $\partial h(\hat x) = \tfrac{\eps^\prime}{\de} \B^\ast$, and we obtain
	\begin{align*}
	0 \in \partial \varphi_b(\hat x)+ \partial \mathbbm{1}_A(\hat x)+ \partial h(\hat x)
	 &= \partial \varphi_b(\hat x)+ N_A(\hat x) + \tfrac{\eps^\prime}{\de} \B^\ast.
	\end{align*}
	Therefore, 
	\begin{equation}
	    \label{eq:inclusion}
	    \left[\partial \varphi_b(\hat x)+ N_A(\hat x)\right]\cap (\tfrac{\eps^\prime}{\de} \B^*) \neq \emptyset.
	\end{equation}
	Furthermore, since $\hat x\in A$ we have that
	$\varphi_b(\hat x)= \norm{\hat x-b} \ge d(b,A)>0$.  By \cite[Corollary 2.4.16]{Za} and the chain rule, we have
	$$
	\partial \varphi_b(\hat x) = \left\{u:\; \ang{u,\hat x-b} = \norm{\hat x-b},\; \norm{u} = 1 \right\}.
	%\hbox{\textcolor{red}{Why?}} 
%	\hbox{\textcolor{blue}{Hoa. Because of \eqref{E1}, and the chain rule}}
	$$
	Altogether, by \eqref{eq:inclusion} there is $x^*\in(-\partial \varphi_b(\hat x))$ such that 
	$$d\left(x^*,N_A(\hat x)\right)\le \frac{\eps^\prime}{\de}< \frac{\eps}{\de},$$
	{where the first inequality holds because, by \eqref{eq:inclusion}, there exists $w\in N_A(\hat x)$ such that $(w-x^*)\in \tfrac{\eps^\prime}{\de} \B^*$,  and the last inequality holds because $\eps^\prime < \eps$.}
% 	\textcolor{purple}{RSB 12Sept: Why (55) implies the previous claim? It seems to me that (55) implies that $(u+x^*)\in \tfrac{\eps^\prime}{\de} \B^*$}
% 	\textcolor{blue}{HB 13Sept: I fixed the mistake}
Noting that \eqref{sl1} holds, and using the facts that $x^*\in(-\partial \varphi_b(\hat x))$, i.e., $\norm{x^*} = 1, \ang{x^*,b-\hat x} =\norm{b-\hat x}$,  and $\hat x \in B[\rho,b]\cap B(\delta,a) \cap A$ we can write
	$$
	\frac{\eps}{\de}\le \inf\left\{d\left(x^*,N_A(x)\right):\; \norm{x^*} = 1,\;\ang{x^*,b-x} =\norm{b-x},\; x\in \ B[\rho,b]\cap B(\delta,a) \cap A \right\}< \frac{\eps}{\de},
	$$
	which is a contradiction. This implies that we must have $d(b,A)\le \norm{a-b}-\eps$.
\end{proof}
%{\color{magenta}RSB 16/02: Hoa please write the comment below on the Hilbert space case in full as a corollary of Theo 4.1, because it is used in this form in Theo 5.1. It is not clear what the full statement is.}
%Recall that in Hilbert spaces a point $a \in A$ is a projection of $b\notin A$ onto the set $A$ if and only if $b-a \in N_A(a)$, which is equivalent to $d\left(\tfrac{b-a}{\norm{b - a}}, N_A(a)\right) = 0$.

%Regi 6 October

When the space $X$ is Hilbert, Theorem~\ref{sl} can be simplified as follows.
\begin{corollary}\label{cor:sl3}
Suppose $A$ is a closed convex set of a Hilbert space, $a\in A$,  $b\notin A$ with $\rho:=\norm{a-b}$, and $\eps >0$. If there is $\de >0$ such that the following inequality holds
\begin{equation}\label{sl2}
	\inf\left\{d\left(\frac{b-x}{\norm{b-x}},N_A(x)\right):\; x\in B[\rho,b]\cap B(\delta,a) \cap A \right\}\ge\frac{\eps}{\de},
	\end{equation}
then $d(b,A)\le \norm{a-b}-\eps$.
\end{corollary}
\begin{proof}
In a Hilbert space, the element $x^*$ we find in the proof of Theorem \ref{sl} can be taken as $x^*:=(\hat x-b)/\|\hat x-b\|$, so the expression in \eqref{sl1} becomes \eqref{sl2}.
\end{proof}

\begin{figure}
\begin{center}
% \begin{tikzpicture}[scale=0.5]
% %\draw (0,3) node {$B=\{x:\mathbb{A}x\leq b\}$};
% \filldraw[blue,opacity=0.1] (0,0) circle (9em);
% %\draw[blue,opacity=0.2] (0,0) circle (9em);
% %\draw[->,blue](2,2) -- (3,3);
% %\filldraw[gray] (0,0) circle (2pt);
% \filldraw (4,4) circle (2pt);
% \filldraw (4.5,4.5) node {$b$};
% \filldraw (3.2,0) circle (2pt);
% \filldraw (3.6,-0.4) node {$a$};
% %\draw[dashed] (0,0)-- (4,4);
% \draw (3.2,0) -- (4,4);
% \draw[->,blue](3.2,0) -- (4.5,0);
% \filldraw[gray,opacity=0.09] (3.2,0) circle (5em);
% \filldraw (2.7,1.6) circle (2pt);
% \draw[gray] (2.7,1.6) -- (4,4);
% \draw[->,blue] (2.7,1.6) -- (4,2.5);
% \filldraw (2.4,0.5) node {$\delta$};
% \filldraw[red] ([shift=(34:1cm)]2.7,1.6) arc (34:60:1cm);
% \filldraw (4,3) node {$\al$};
% \draw[gray,<->] (2.7,1.6) -- (3.2,0);
% \end{tikzpicture}
\begin{tikzpicture}[scale=0.5]
%\draw (0,3) node {$B=\{x:\mathbb{A}x\leq b\}$};
\filldraw[gray,opacity=0.1] (0,0) circle (9em);
%\draw[blue,opacity=0.2] (0,0) circle (9em);
%\draw[->,blue](2,2) -- (3,3);
%\filldraw[gray] (0,0) circle (2pt);
\filldraw (4,4) circle (2pt);
\filldraw (4.5,4.5) node {$b$};
\filldraw (3.2,0) circle (2pt);
\filldraw (3.6,-0.4) node {$a$};
\filldraw (2.3,1.8) node {$x$};
\filldraw (6,2.6) node {$N_A(x)$};
\filldraw (0,0) node {$A$};
%\draw[dashed] (0,0)-- (4,4);
%\draw (3.2,0) -- (4,4);
%\draw[->](3.2,0) -- (4.5,0);
\filldraw[gray,opacity=0.3] (3.2,0) circle (5em);
\filldraw (2.7,1.6) circle (2pt);
\draw[->] (2.7,1.6) -- (4,4);
\draw[->] (2.7,1.6) -- (5.3,3.2);
\filldraw (2.4,0.5) node {$\delta$};
\filldraw ([shift=(34:1cm)]2.7,1.6) arc (34:60:1cm);
\filldraw (3.6,2.6) node {$\al$};
\draw[gray,<->] (2.7,1.6) -- (3.2,0);
\end{tikzpicture}
\end{center}
\caption{The distance from $b$ to the set $A$ is bounded {above} by $\norm{b-a} - \de\al$ (so $\eps=\de\al$). Here, $\al$ is the smallest cosine of all the angles between  vectors $x-b$ and normal cone $N_A(x)$ for $x\in A$, and $\delta$ is the size of the neighbourhood around $a$, i.e., we are taking $x \in A\cap B_\de(a)$. Note here that we only consider $x\in A\cap B_\de(a)$ such that $\norm{x-b}\le \norm{b-a}$.\label{fig1b}}
\end{figure}
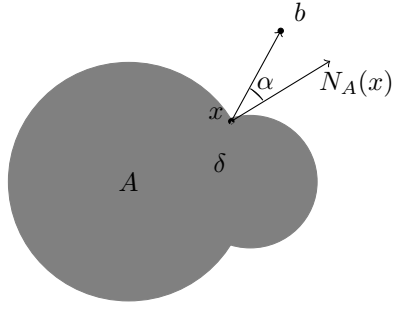

\begin{remark}
The aim of Corollary \ref{cor:sl3} is to establish a sufficient condition for $a$ being ``far enough" from being a projection of $b$ onto $A$. Note that $a=P_A(b)$ if and only if $b-a \in N_A(a)$. Equivalently, $d\left(\tfrac{b-a}{\norm{b - a}}, N_A(a)\right) = 0$. Hence, to ensure we are far from the latter situation, we require \eqref{sl2} to hold, not merely at $a$, but at every point $x\in B[\rho,b]\cap B(\delta,a) \cap A$. We quantify this property by showing that, if \eqref{sl2} holds, then we must have $\norm{a-b} > d(b,A)+\eps$. The latter, in turn, means that the difference $\norm{a-b}-d(b,A)$ is bounded from below by a constant $\eps$. Note also that the opposite inequality to \eqref{sl2} is the optimality condition for $\eps$-projections, where the later define points in $A$ that are within distance $d(b,A)+\eps$ from $b$. Geometrically, \eqref{sl2} ensures that the cosine of the angle between $b-x$ and a vector in $N_A(x)$ is always bigger than a positive constant $\frac{\eps}{\de}>0$, where $x\in B[\rho,b]\cap B(\delta,a) \cap A$. 
Figure~\ref{fig1b} illustrates an example on estimating the distance from a point to a set.
% \textcolor{purple}{RSB 15/08: Hoa, I edited this, but I am not sure I am quite clear, can you please check? Also, figure \ref{fig1b} is not quite clear to me, Hoa, can you make it more clear?}
% \textcolor{blue}{HB: I fixed the figure. I hope it is clearer now?}
\end{remark}

%\todo[inline]{MKT: The figure is a bit unclear to me. For instance, I don't see the distance from $b$ to $A$ illustrated. I'm not not sure what the other lines should represent.}
%\begin{example}
%\label{Ex1}
%Consider $A$ is the unit circle in $(\R^2,\norm{\cdot}_2)$ and $b = (x,y) \notin %A$ (with $x,y>0$). Then, for any $z\in A$, we have $N_A(z) = \R_+(z)$ and 
%$$d\left(\frac{b-z}{\norm{b-z}},N_A(z)\right)=\tfrac{1}{\norm{b-z}}.$$ Consider %$a = (1,0)$, we have $\norm{a-b} = \sqrt{(x-1)^2+y^2}$.
%\end{example}

 %Recall that for any $\eps >0$, we denote the set of $\eps$-projections of $b$ onto $A$ is $P_A(b,\eps):=\{x\in A:\; \norm{x-b}< d(b,A)+\eps \}$. Note that $P_A(b) = \bigcap_{\eps >0} P_A(b,\eps)$.
%It is well-known that
%$$
%a\in P_A(b) \Longleftrightarrow b-a \in N_A(a).
%$$

Since the projection of $b$ onto $A$ is an element $a\in A$ such that $d(b,A)=\|a-b\|$, an  $\eps$-projection of $b$ onto $A$ can be understood as an element $a'\in A$ such that $d(b,A)>\|a'-b\|-\eps$. This motivates the following corollary of Theorem \ref{sl}.

\begin{corollary}
\label{C2b}
Consider a Banach space $X$, a nonempty closed convex set $A$, $b\notin A$, $a\in A$, $\rho = \norm{a-b}$ and $\eps \ge 0$. For any given $\de>0$, define the set
\[
C(\delta):=\{ 
(x,x^*)\in X\times X^*\::\: x\in B[{\rho},b]\cap B(\delta,a)\cap A \hbox{ and }\norm{x^*} = 1,\ang{x^*,b-x} =\norm{b-x} \}.
\]
Assume that $a$ is an $\eps$-projection of $b$, in the sense that $\norm{a-b}< d(b,A)+\eps$.
Then, for every $\de >0$, there exists $(x_0,x_0^*)\in C(\delta)$ such that
\[
d(x_0^*,N_A(x_0)) <\frac{\eps}{\de}.
\]
%Then, we must have 
% if for any $\de >0$, there exist $x_0\in B[{\rho},b]\cap B(\delta,a)\cap A$ and $x^*\in X^*$ s.t.
% \begin{align}\label{T1.1}
% \norm{x^*} = 1,\ang{x^*,b-x_0} =\norm{b-x_0},\AND d(x^*,N_A(x_0)) <\frac{\eps}{\de}.
% \end{align}
\end{corollary}
\begin{proof}
Assume that the conclusion of the corollary is not true. Namely, assume that there is $\de>0$ such that for all $(x,x^*)\in C(\delta)$ we have
    $$
 d(x^*,N_A(x)) \ge \frac{\eps}{\de}.
    $$
Hence,
\[
\inf_{(x,x^*)\in C(\de)} d(x^*,N_A(x)) \ge \frac{\eps}{\de},
\]
which, by definition of $C(\de)$, the above inequality is exactly \eqref{sl1}. By Theorem~\ref{sl}, we must have $\norm{a-b}\ge d(b,A)+\eps$, contradicting the fact that $a$ is an $\eps$-projection of $b$.
\end{proof}

%{\color{magenta}RSB 12/01: I think we need to explain how this corollary follows from Theo 4.1.}\\[2mm]

\subsection{Solving problems with SPP: the case with a linear objective}\label{Sec:LP}
%\todo[inline]{MKT: Throughout this section, we consider a set $A$ which is $2\alpha$-sharp for some $\alpha\in(0,1)$. This is equivalent to saying $A$ is $\beta-$sharp with $\beta\in(0,2)$. However, this doesn't make complete sense because we know from the definition of sharpness that we must have $\beta\leq 1$.}
%\textcolor{purple}{Reply to MKT: Thanks, $\al$ must be in $(0, 1/2]$ for this condition to hold, will change it.}
We start this section by considering the following minimization problem with linear objective:
\begin{align}
    \tag{P}\label{LP1} 
    \min \quad & \ang{x^*,x}\\
    \notag
    \text{s.t.}\quad & x\in A,
\end{align}
here $x^*\in \mathcal S$ and $A$ is a closed convex set. Denote by $\mathbb{S}$ the set of solutions of problem \eqref{LP1}. We assume that $\mathbb{S}\not=\emptyset$. The optimality conditions for Problem \eqref{LP1} imply that 
$$
\mathbb{S}= \left\{x: 0\in x^*+N_A(x) \right\}.
$$
The following theorem shows that Problem~\eqref{LP1} can be solved by projecting an infeasible point onto the feasible region $A$ if the set $A$ is sharp. We will use the following fact, which is a consequence of Fact~\ref{L2}. For any $p,q$ nonzero vectors, we have 
\begin{equation}\label{eq:cones}
d(p,\cone[q])=d(q,\cone[p]).
\end{equation}

%{\color{magenta}RSB 8/2/22: In Theo below is $\al\in (0,1)$?}
\begin{theorem}
\label{T3}
Let $A$ be a closed convex set of a Hilbert space $\mathcal H$, and $x^*\in \mathcal S$. Suppose that $A$ is $\al$-sharp w.r.t.\ $-x^*$ for some $\al\in (0,1]$. Suppose also that $v\in \mathcal H$ satisfies following conditions:
\begin{enumerate}
    \item[1.] $\ang{x^*,v} < \inf_{x\in A}\ang{x^*,x}$;
    \item[2.] $(1- ({\al}/{2})^2) d(v,A) < \inf_{x\in A}\ang{x^*,x-v}$.
%  \item[2.] $(1- \al^2) d(v,A) < \inf_{x\in A}\ang{x^*,x-v}$,
  \end{enumerate}
Then the projection of $v$ onto $A$ is a solution of Problem~\eqref{LP1}.
\end{theorem}
\begin{proof}
Suppose to the contrary that there exists $v\in \mathcal H$ satisfying Conditions~1 and 2 such that the projection of $v$ onto $A$ is not a solution of \eqref{LP1}. Since $A$ is $\al$-sharp w.r.t.\ $-x^*$, we have
\begin{equation}
    \label{T3.1}
    \inf_{\substack{x\in A\\ -x^*\notin N_A(x)}} d(-x^*,N_A(x)) \ge \al.
\end{equation}
Take $y= P_A(v)$, then $v-y\in N_A(y)$. Because $y$ is not a solution of the convex problem \eqref{LP}, we must have  $0\notin x^*+N_A(y)$.
%, and because $v-y\in N_A(y)$, we have $-x^*\in \cone[v-y]$. 
From \eqref{T3.1} and the fact that $-x^*\notin N_A(y)$, we have 
$$
d(-x^*, N_A(y)) \ge \al >0.
$$
%{\color{magenta}RSB 12/01, Question: the above has minus sign in front of $N_A$, how does this inequality follows from \eqref{T3.1}, where there is no minus sign?}
Combining this inequality with the inclusion $v-y\in N_A(y)$ yields
\begin{equation}\label{eq:A1}
    d\left(\frac{y-v}{\norm{y-v}}, \cone[x^*] \right)=d(x^*,\cone[y-v])= d(-x^*,\cone[v-y])\ge d(-x^*, N_A(y))  \ge \al,
\end{equation}
where we also used \eqref{eq:cones} in the first equality. Consider the closed half space $F:=\left\{x\in \mathcal{H}:\; \ang{x^*,x-v} \le 0 \right\}$.
From Condition~1, $\ang{x^*,v} < \inf_{x\in A}\ang{x^*,x}$, hence the sets $A$ and $F$ are disjoint, i.e., $A\cap F =\emptyset$.
%Recall that $y-b\notin \R_+(x^*) = N_B(b)$, so $b$ is not the projection of $y$ onto $B$. 
Now, we are going to apply Corollary \ref{cor:sl3}. Namely, we will show that \eqref{sl2} holds with lower bound $\alpha/2$. Indeed, we apply this corollary to the set $F$, $v\in F$, and $y\notin F$ to estimate the distance $d(y,F)$ relative to $\norm{v-y}$. Set $\de :=(\al/2)\norm{v-y}$, and $\rho := \norm{v-y}$, we claim that
\begin{equation}
    \label{eq:claim1}
    \inf\left\{d\left(\frac{y-z}{\norm{y-z}},N_F(z)\right):\;z\in F\cap B(\de,v)\cap B[{\rho},y]\right\} \ge \al/2.
\end{equation}
Indeed, take $z\in F\cap B(\de,v)\cap B[{\rho},y]$. We consider two cases.

Case 1. If $z\in \Int F$, then $N_F(z) = \{0\}$, hence $d\left(\frac{y-z}{\norm{y-z}},N_F(z)\right) = 1$.

Case 2. If $z\notin \Int F$, then the definition of $F$ implies that $N_F(z) = \cone[x^*]$.
Since $z\in B[\rho,y]$, we use \eqref{eq:A1} to write
\begin{align*}
d\left(\frac{y-z}{\norm{y-z}},N_F(z)\right)& = d\left(\frac{y-z}{\norm{y-z}},\cone[x^*] \right)=\inf_{t\ge 0} \left\|\frac{y-z}{\norm{y-z}} - t x^*\right\| \\
&\\
&= \frac{\norm{y-v}}{\norm{y-z}} d\left(\frac{y-z}{\norm{y-v}},\cone[x^*]  \right)\\
&\\
&\ge d\left(\frac{y-z}{\norm{y-v}},\cone[x^*]\right)= d\left(\frac{y-v}{\norm{y-v}}+\frac{v-z}{\norm{y-v}},\cone[x^*] \right)\\
&\\
&= \inf_{t\ge 0} \norm{\frac{y-v}{\norm{y-v}}+\frac{v-z}{\norm{y-v}}- t x^*} \ge \inf_{t\ge 0} \norm{\frac{y-v}{\norm{y-v}}- t x^*} -\frac{\norm{v-z}}{\norm{y-v}} \\
&\\
 &= d\left(\frac{y-v}{\norm{y-v}},\cone[x^*] \right) -\frac{\norm{v-z}}{\norm{y-v}} \ge \al - \al/2 = \al/2,
\end{align*}
where in the fourth equality we used the change of variable $t\to \Tilde{t}:=t \|y-z\|/\|y-v\|$, and we used the fact that $\|y-z\|\le \rho=\|y-v\|$ in the first inequality. As for the last inequality, we use \eqref{eq:A1} for obtaining the lower bound of the first term. For the second term, recall that $z\in B(\de,v)$, so $\norm{v-z}< \delta=(\al/2)\norm{v-y}$. These establish the last inequality. Hence, our claim \eqref{eq:claim1} holds. By Theorem~\ref{sl}, we have 
\begin{equation}
    \label{eq:d1}
    d(y,F) \le \norm{y-v} - (\al/2)\de = \norm{y-v} - (\al/2)^2\norm{y-v}=(1-(\al/2)^2) \norm{y-v}.
\end{equation}
To arrive at a contradiction, we will use Condition 2. {Suppose $w\in F$ is the projection of $y$ onto $F$, then $y-w\in N_F(w)$. By definition of $F$ we know that $ N_F(w)= \cone[x^*]$. Hence, $x^*=(y-w)/\|y-w\|$ and therefore, $d(y,F)=\norm{y-w}=\ang{x^*,y-w}$. Taking into account that $x^*\in N_F(v)\cap N_F(w)$, we obtain $\ang{x^*,w-v}\le 0$, and $\ang{x^*,v-w}\le 0$, and hence $\ang{x^*,w-v}=0$. From the previous equality $d(y,F)=\ang{x^*,y-w}$, we have $\ang{x^*,y-v}=\ang{x^*,y-w}+\ang{x^*,w-v}= \ang{x^*,y-w} =d(y,F)$}. So, the following estimation holds
$$ \inf_{x\in A}\ang{x^*,x-v}\le \ang{x^*,y-v}= d(y,F) \le (1-(\al/2)^2) \norm{y-v}= (1-(\al/2)^2) d(v,A),$$ 
where we also used \eqref{eq:d1} in the second inequality, and the definition of $y$ in the last equality. The above expression contradicts condition 2. Therefore, we must have that $P_A(v)$ solves problem (P).
\end{proof}

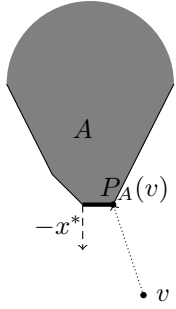
\begin{figure}
\begin{center}
% \begin{tikzpicture}[scale=0.5]
% %\draw (0,3) node {$B=\{x:\mathbb{A}x\leq b\}$};
% \fill[blue,opacity=0.1] (0,0)--(2,2)--(3,4)--(-2,4)--(-1,1)--(0,0);
% \fill[blue,opacity=0.1] (-2,4)--(3,4) arc(0:180:2.5) --cycle;
% \draw[blue,opacity=0.2] (0,0)--(2,2)--(3,4);
% \draw[blue,opacity=0.2] (-2,4)--(-1,1)--(0,0);
% \draw[blue] (3,-2) -- (-3,2);
% \draw[->,blue] (0,0) -- (1,1.5); 
% \filldraw (0,2.5) node {$A$};
% \filldraw (0.5,-2) circle (2pt);
% \filldraw (0,0) circle (2pt);
% \draw[dashed] (0.5,-2)--(0,0);
% %\draw[blue] (3,-3.7) -- (-3,0.3);
% \filldraw (0,-2.5) node {$v$};
% \filldraw (3,-1.5) node {$x^*$};
% \end{tikzpicture}
% \end{center}
\begin{tikzpicture}[scale= 0.4]
		\fill[gray,opacity=0.1] (0,0)--(1,0)--(2,2)--(3,4)--(-2.5,4)--(-1,1)--(0,0);
        \fill[gray,opacity=0.1] (-3,4)--(3,4) arc(0:180:2.75) --cycle;
        \draw[opacity=0.5] (0,0)--(1,0)--(2,2)--(3,4);
        \draw[black,ultra thick,opacity=0.5] (0,0)--(1,0);
        \draw[opacity=0.5] (-2.5,4)--(-1,1)--(0,0);
       % \draw[black] (-3,0) -- (3,0);
        \draw[->,black,densely dashed] (0,0) -- (0,-1.5); 
        \filldraw (0,2.5) node {$A$};
        \filldraw (-0.8,-0.7) node {$-x^*$};
        \node[circle,draw=black, fill = black, inner sep=0pt,minimum size=2pt, label=right:{$v$}] (b) at (2,-3) {};
        \draw[densely dotted,black,->](2,-3) -- (1,0);
        \node[circle,draw=black, fill = black, inner sep=0pt,minimum size=2pt] (b) at (1,0) {};
        \filldraw (1.7,0.5) node {$P_A(v)$};
        %\draw[densely dotted,black,->](2,0.5) -- (1.4,0.8);
		\end{tikzpicture}
\end{center}
\caption{Illustration of Theorem~\ref{T3}: the set $A$ is sharp w.r.t vector $-x^*$. Vector $v$ satisfies conditions 1\&2, hence the projection $P_A(v)$ of $v$ onto the set $A$ is the solution of the minimization problem $\inf_{A}\ang{x^*,\cdot}$.\label{fig2}}
% \textcolor{purple}{RSB 12 Sept: Can we please explain this figure better? Where is $-x^*$? Let's also indicate $P_A(v)$. What does the black line stand for? }
\end{figure}

As a consequence of Theorem \ref{T3}, if $A$ is $\alpha$-sharp w.r.t. a vector $-x^*$, Problem \eqref{LP1} can be solved by projecting onto $A$ an infeasible point $v$ s.t. conditions 1 and 2 hold (see Figure~\ref{fig2}).  Hence, it is important to be able to construct such vectors. It is clear that condition 1 in Theorem \ref{T3} follows from condition 2. The next lemma shows that, once we have a vector verifying Condition 1, we can always construct a translation of the vector that verifies Condition 2.

\begin{lemma}
\label{rem:trans}
With the notation of Theorem \ref{T3}, assume that $v\in \mathcal H$ verifies Condition 1 and fix $\alpha\in (0,1]$. Assume that Condition 2 with parameter $\al$ does not hold for $v$. Define
\begin{equation}\label{eq:hypo}
\theta(v):= \inf_{x\in A}\ang{x^*,x-v},\,\hbox{ and }  \mu_0:= \frac{(1-(\alpha/2)^2)d(v,A)-\theta(v)}{(\al/2)^2}.
\end{equation}
Then $\theta(v)>0$ and $\mu_0\ge 0$. Moreover, if $\mu>\mu_0$, then $u:=v - \mu x^*$ verifies conditions 1 and 2 from Theorem \ref{T3}.
\end{lemma}
\begin{proof}
The fact that $\theta(v)>0$ is equivalent to the validity of Condition 1 for $v$, so it holds by assumption. The fact that $\mu_0\ge 0$ is equivalent to the assumption that $v$ fails to verify Condition 2 for the given $\alpha$. Altogether, we have that
\begin{equation}\label{eq:C2 false}
0<  \theta(v)\le  (1-(\al/2)^2) d(v,A).
\end{equation}
We proceed to prove that Conditions 1 and 2 hold for $u$ if $\mu>\mu_0$. Use the definition of $\theta(v)$ to write, for all $x\in A$,
$$
\ang{x^*,x-u} = \ang{x^*,x-(v-\mu x^*)} = \ang{x^*,x-v} +\mu\ge \theta(v)+ \mu > \mu>0,
$$
where we used the definition of $u$ in the first equality, and the fact that $x^*\in {\cal S}$ in the second one. Therefore, $\inf_{x\in A}\ang{x^*,x-u}\ge \mu>0$ and Condition 1 holds for $u$. The above expression also yields
\begin{equation}\label{eq:mua}
\ang{x^*,x-u}\ge \theta(v)+ \mu>0,
\end{equation}
for all $x\in A$.  Let us check now that Condition 2 holds for $u$. Using again the fact that $x^*\in {\cal S}$ gives
\begin{equation}\label{eq:mu1}
d(u,A) = d(v-\mu x^*,A) =\inf_{x\in A } \| v-\mu x^* -x\| \le \inf_{x\in A } 
\| v -x\| +\mu \|x^*\|= d(v,A)+\mu.
\end{equation}
Using the definition of $\mu_0$, we rewrite the inequality $\mu>\mu_0$ as
\begin{equation}\label{eq:mu0}
    d(v,A)<\dfrac{(\al/2)^2 \mu +\theta(v)}{(1-(\al/2)^2)}.
\end{equation}
Using \eqref{eq:mu0} in \eqref{eq:mu1} yields
$$
d(u,A) \le d(v,A)+\mu< \frac{((\al/2)^2  \mu +\theta(v))}{1-(\al/2)^2}+\mu=\frac{\mu +\theta(v)}{1-(\al/2)^2}\le \frac{\ang{x^*,x-u}}{1-(\al/2)^2},
$$
where we also used \eqref{eq:mua} in the last inequality. Since the above inequality holds for every $x\in A$, we deduce that Condition 2 holds for $u$. 
  
\end{proof}

The argument in Lemma \ref{rem:trans} is the main idea behind the next proposition. Namely,
%{\color{red}[MKT: Don't read past here.]}
if a vector $v\in \mathcal H$ satisfies Condition 1 and not Condition 2, then we translate $v$ by a large enough multiple of  $-x^*$, so that Condition 2 holds for the translated vector.

\begin{proposition}\label{1s}
Suppose $A$ is a closed convex set of a Hilbert space $\mathcal H$, and $x^*\in \mathcal S$ such that the set $A$ is $\al$-sharp w.r.t. $-x^*$, for some $\al\in (0,1]$. 
Consider $v \in \mathcal H$ such that $\ang{x^*,v} < \min_{x\in A}\ang{x^*,x}$. Then, the projection of $u:=v - \mu x^*$ onto  $A$, where $\mu \ge \frac{4-\al^2}{\al^2}d(v,A)$, is a solution of \eqref{LP1}.
\end{proposition}
\begin{proof}
The proof follows by noting that $\frac{4-\al^2}{\al^2}d(v,A)=\frac{1-(\al/2)^2}{(\al/2)^2}d(v,A)> \mu_0$, with $\mu_0$ as in Lemma \ref{rem:trans}. Using the lemma, we see that $u$ verifies conditions 1 and 2 in Theorem \ref{T3}. Therefore, the claim follows directly from the theorem.
%{\color{magenta} RSB 11/2 Isn't the previous argument enough to complete the proof, maybe we don't need what follows.}
\end{proof}

%\textcolor{blue}{HB: Proposition 5.3 is fixed. Don't read past this}
%Corollary~\ref{1s} establishes that from an infeasible point $x_0\in X$ such that $\ang{c,x_0} < \min_{x:\,\mathbb{A}x\leq b}\ang{c,x}$, we can shift $x_0$ by the vector $- \mu c$ in such a way that its projection onto $B :=\{x\in \R^n: \mathbb{A}x \le b\}$ is the solution of LP. To determine $\mu$, we need to estimate the upper bound for $d(x_0,B)$. This is possible if we know a feasible point $\bx \in B$ in advance, then we can take $\mu > \frac{(1-\al^2)}{\al^2}\norm{\bx - x_0}$.

We illustrate Proposition~\ref{1s} with the following two examples.

\begin{example}\label{exa 3.8}
Consider Problem (P) with $\mathcal{H}=\R^2$, \textcolor{black}{$x^*=(0,1)$} and 
\(
A:=\{x\in \R^2\::\: \mathbb{A}x\le b\},
\)
where $\mathbb{A}=\begin{bmatrix}1& -1 \\ -1 & -1\end{bmatrix}$ and $b=(0,0)$. Our first task is to determine the modulus of sharpness of $A$. It is easy to check that, for every $x\in {\rm bd}A$, 
\[
N_A(x)=\left\{
\begin{array}{cc}
  {\rm cone}[(-1,-1)]=:K_1  & \hbox{ if } x_1<0,  \\
  &\\
     {\rm cone}[(1,-1)]=:K_2   & \hbox{ if } x_1>0, \\
     &\\
     {\rm cone}[(-1,-1),(1,-1)]  & \hbox{ if } x_1=0, 
\end{array}
\right.
\]
which can be graphically verified from Figure \ref{fig1}. Note that we have $-x^*\not\in N_A(x)$  if and only if $x_1\not=0$. Using \eqref{T3.1} in Theorem \ref{T3}, we have
\[
\begin{array}{rcl}
   \inf_{\substack{x\in A\\ -x^*\notin N_A(x)}} d(-x^*,N_A(x))    & = &  \min\{ d(-x^*, K_1), d(-x^*,  K_2)  \}=\sqrt{2}/2= 2(\sqrt{2}/4).
\end{array}
\]
Hence, \eqref{T3.1}  holds with $\alpha:=\sqrt{2}/2$. Consequently, $x_0$ will solve (P) if it verifies conditions 1 and 2 in Theorem \ref{T3}. If only condition 1 holds, then we can use Proposition~\ref{1s} and find $\mu_0$
s.t. $P_A(x_0-\mu x^*)$ solves the LP for $\mu>\mu_0$. For instance, take $x_0:=(-1,-1/2)$. It is easy to check that $x_0$ verifies condition~1 in Theorem~\ref{T3}, but not condition~2, and that $d(x_0, A)=3\sqrt{2}/4$. With the notation of Proposition~\ref{1s} and $\al=\sqrt{2}/4$, we need to take $\mu$ such that 
\[
\mu> d(x_0,A)(4-\al^2)/\al^2= 21\sqrt{2}/4.
\]
Take $\mu=10>21\sqrt{2}/4$. So $u:=x_0-\mu x^* = (-1,(-1/2)-10 )=(-1,-21/2)$ with $P_A(u)=(0,0)$, the solution of (P). An illustration of this example is shown in Figure~\ref{fig1}.
\begin{figure}
\begin{center}
\begin{tikzpicture}[scale=0.5]
\draw (0,3) node {$A=\{x:\mathbb{A}x\leq b\}$};
\filldraw[gray,opacity=0.1] (-5,5) -- (0,0) -- (5,5);
\draw[<->,thick] (-5.5,5.5) -- (0,0) -- (5.5,5.5);

%\filldraw[blue,opacity=0.1] (-5,-5) -- (0,0) -- (5,-5);
\draw[->,thick] (3.5,-1) -- (3.5,-3);
\draw (3.5,-1.5) node[right] {$-x^*$};
\draw[fill] (-3,-1) circle (0.25ex) node[above left] {$x_0$};
\draw[<->,thick] (-5.5,-1) -- (5.5,-1);
\filldraw[gray,opacity=0.1] (-5,-1) -- (5,-1) -- (5,-5) -- (-5,-5);
\draw (0,-3) node {$F=\{x:\langle x^*,x\rangle\leq M\}$};
\end{tikzpicture}
\qquad\quad\quad\quad\quad\quad
\begin{tikzpicture}[scale=0.5]
\draw (0,3) node {$A=\{x:\mathbb{A}x\leq b\}$};
\filldraw[gray,opacity=0.1] (-5,5) -- (0,0) -- (5,5);
\draw[<->,thick] (-5.5,5.5) -- (0,0) -- (5.5,5.5);

\filldraw[gray,opacity=0.1] (-5,-5) -- (0,0) -- (5,-5);
\draw[gray,->,thick] (0,0) -- (0,-2);
\draw (0,-1.5) node[right] {$-x^*$};
\draw[gray] (2.5,-4) node {$N_A(\bar x)$};

\draw[fill] (-3,-1) circle (0.25ex) node[right] {$x_0$};
\draw[->,dotted] (-3,-1) -- (-1,1);
\draw (-1,1) node[above right] {$P_A(x_0)$};

\draw[fill] (-3,-4) circle (0.25ex) node[right] {$x_0-\mu x^*$};
\draw[->,dotted] (-3,-4) -- (0,0);
\draw[fill] (0,0) circle (0.25ex) node[right] {$\bar x=P_A(x_0-\mu x^*)$};
\end{tikzpicture}
\end{center}
\caption{Illustration for Example \ref{exa 3.8}: (left) the disjoint sets $A$ and $F$ and (right) the projection of the point $x_0-\mu x^*$ onto $A$ is the unique solution, $\bar x$, of the LP.\label{fig1}} 
%\textcolor{purple}{RSB 12 Sept: the location of this $x_0-\mu x^*$ does not match the one in Ex 4.8}
%
%\textcolor{blue}{HB I think the figure should be consistent now?}
\end{figure}
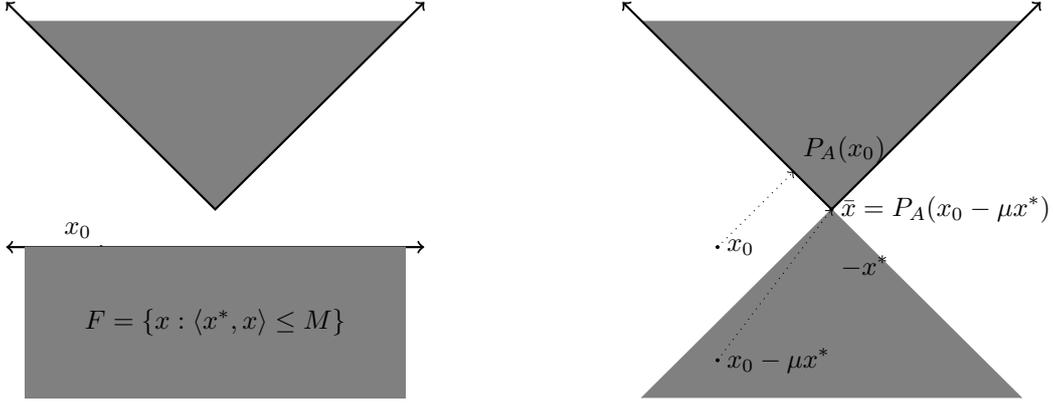
\end{example}

\begin{example}
Consider Problem (P) with $\mathcal{H}=\R^3$, $x^*=(0,1/\sqrt{2},1/\sqrt{2})$ and 
\(
A:=\{x\in \R^2\::\: \mathbb{A}x\le b\},
\)
where $\mathbb{A}=-I$ and $b=(0,0,0)$. Namely, $A=\R^3_{+}$. The solution of (P) is the set $\R_{+}\times \{0\}\times \{0\}$ and zero is its optimal value.  Again, we first determine the modulus of sharpness of $A$. It can be checked that the only cases in which $-x^*\in N_A(x)$ for $x\in {\rm bd}A$ is when $x_1\ge 0$ and $x_2=x_3=0$. Hence, we need to compute
$d(-x^*,N_A(x))$ for $x$ in the following set:
\[
T:=\{x\in {\rm bd}A\::\: x_1\ge 0 \hbox{ and } x_2,x_3 \hbox{ are not simultaneously zero}\}
\]
It can be checked that
% It is easy to check that, for every $x\in {\rm bd}A$, 
% \[
% N_A(x)=\left\{
% \begin{array}{cc}
%   {\rm cone}[(-1,0,0)]=:K_1  & \hbox{ if } x_1=0,  \\
%   &\\
%      {\rm cone}[(0,-1,0)]=:K_2   & \hbox{ if } x_2=0, \\
%      &\\
%      {\rm cone}[(0,0,-1)]=:K_3  & \hbox{ if } x_3=0, 
% \end{array}
% \right.
% \]
% We have that in all the three cases it holds that $-x^*\not\in N_A(x)$. Since $d(-x^*, K_1)=1$ and $d(-x^*, K_2)= d(-x^*, K_3)=1/\sqrt{2}$ we have
\[
\begin{array}{rcl}
   \inf_{{x\in T}} d(-x^*,N_A(x))    & = &  1/\sqrt{2},
   \end{array}
\]
so \eqref{T3.1}  holds with $\alpha:=1/\sqrt{2}$. Take $x_0:=(1,-1,0)$ then 
$d(x_0,A)=1$, and $x_0$ verifies condition 1 in Theorem \ref{T3} but not condition 2. With the notation of Proposition~\ref{1s}, we need $\mu$ such that 
\[
\mu> \dfrac{(4-\al^2)}{\al^2} d(x_0,A)= 7.
\]
Take $\mu=7\sqrt{2}$ so $u:=x_0-\mu x^* = (1,-8,-7)$ with $P_A(u)=(1,0,0)$ a solution of (P).
\end{example}

\subsection{SPP for a general case}

We extend next Theorem~\ref{T3} to a problem where the objective function is an arbitrary convex and lsc function $f:\mathcal H\to \R_{\infty}$. Namely, we consider the convex problem
\begin{align}
    \label{P2}\tag{CP} \min_{x\in A}f(x),
\end{align}
where $A$ is a closed and convex set. For this problem, we will assume that $f$ and $A$ are such that $\partial (f+\mathbbm{1}_A)=\partial f + N_A$ over $\dom f \cap A$. The latter is true, for instance, when some standard constraint qualification holds (see \cite[Corollary~16.38]{bauschke2011convex}), e.g., when $\Int{A}\cap \dom f\not=\emptyset$ or $\Int{(\dom f)}\cap A\not=\emptyset$. The function $f_A:=(f+\mathbbm{1}_A)$ will have a crucial role in the next result. Note that \begin{equation}
    \label{fA}
\epi f_A=\{ (x,t)\in A\times \R\::\: f(x)\le t\,\}=(A\times \R) \cap \epi f   
\end{equation}
By imposing a sharpness condition on the set $\epi f_A$ w.r.t. the vector $(0_{\mathcal{H}},-1)$, we can recover a solution of problem \eqref{P2} by using Theorem~\ref{T3} in the extended space $\mathcal H\times \R$. 
\medskip

%\textcolor{purple}{RSB: 23/08  Theo below UNDER CONSTRUCTION, please disregard until ready.}
\medskip

\begin{theorem}
\label{T6}
Suppose that the convex problem \eqref{P2} has solutions with optimal value $M:=\inf_{x\in A}f(x)$ and nonempty set of optimal solutions denoted by $\mathbb{S}$. Assume that the following conditions hold. 
\begin{itemize}
\item[(i)] The set $\epi f_A=:\Tilde{A}$ given in \eqref{fA} is $\al$-sharp w.r.t $z^*:=(0_{\mathcal{H}},-1)$ for $\al \in (0,1)$. 
\item[(ii)] Let $(v,t)\in \mathcal H\times \R$ be such that
\begin{enumerate}
    \item[(a)] $t < M$; and
    \item[(b)] $(1- (\al/2)^2) d((v,t),\Tilde{A}) < (M - t)$.
\end{enumerate}
\end{itemize}
In this situation, consider $P_{\Tilde{A}}(v,t)=(w,f(w))$. Then $w\in \mathbb{S}$ and hence $f(w)=M$.
\end{theorem}
\begin{proof}
With the notation of (i), problem (CP) is equivalent to the following problem
\begin{align}
    \label{Ep}\tag{EP} \min_{(x,s)\in \Tilde{A} }\,s,
\end{align}
which has the same optimal value as (CP), and linear objective $\psi: \mathcal{H}\times \R \to \R$ defined as $\psi(x,s):=s=\ang{(0_{\mathcal{H}},1),(x,s)}$. Note that (EP) is an optimization problem with a linear objective and by (i), its constraint set is $2\al$ sharp w.r.t $z^*:=(0_{\mathcal{H}},-1)$. Take  $\Tilde{v}:=(v,t)\in \mathcal{H}\times \R$ verifying assumptions (a)(b). We claim that this implies that conditions 1 and condition 2 in Theorem \ref{T3} hold, with $x^*:=(0_{\mathcal{H}},1)$, and $A:=\Tilde{A}$. Indeed, condition (a) re-writes as
\[
t=\ang{(0,1),(v,t)}< M=\inf_{A} f(x)= \inf_{(x,s)\in \Tilde{A}} \, \ang{(0,1),(x,s)},
\]
which is Condition 1 in Theorem \ref{T3} for $\Tilde{v}:=(v,t)$ and $x^*:=(0_{\mathcal{H}},1)$. Condition 2 in Theorem \ref{T3} follows directly from (b) and the definitions. Therefore, we are in conditions of Theorem \ref{T3}, and $P_{\Tilde{A}}(v,t)$ solves (EP). By \cite[Proposition 29.35]{bauschke2011convex}, $P_{\Tilde{A}}(v,t)=(w,f(w))$ where $w\in A$ is the unique solution of the inclusion
\[
\dfrac{v-w}{(f(w)-t)}\in \partial f(w).
\]
Note that $t<M\le f(w)$ so $(f(w)-t)>0$. Since $P_{\Tilde{A}}(v,t)=(w,f(w))$ solves (EP) this means that $f(w)=M$ and since $w\in A$ we must have $w\in \mathbb{S}$.
\end{proof}

\begin{remark}\label{rem:KL}
We give in Proposition \ref{EXP1} a necessary and sufficient condition for assumption (i) in Theorem \ref{T6} to hold for $f_A$.
\end{remark}

\begin{remark}\label{rem: f polyhedral 2}
Let $A$ be a closed and convex set and assume that $f_A:=f+\mathbbm{1}_A$ is polyhedral. By Corollary \ref{cor: f polyhedral}, we know that there exists $\al<1$ such that property \eqref{eq:extended-KL} holds for $\beta:=\al/\sqrt{1-\al^2}$. The latter fact, combined with Proposition \ref{EXP1}, imply that 
$\epi f_A=:\Tilde{A}$ is $\al$-sharp w.r.t. $(0_{\mathcal{H}},-1)$. Hence, we are in the situation of Theorem \ref{T6}.% with $\al/2$ instead of $\al$. 
Therefore, if $(v,t)$ is such that $t<M$ and
\[
\left(1-\frac{\al^2}{4}\right) d((v,t), \Tilde{A})<M-t,
\]
then $P_{\Tilde{A}}(v,t)=(\bar x, f(\bar x))$ is such that $\bar x$ solves (CP).
\end{remark}

The following result considers problem (CP) on its own, and establishes a sharpness condition under which a projection onto $A$ solves (CP).

\begin{theorem}
\label{T4}
Suppose that the convex problem \eqref{P2} has solutions and denote by $\mathbb{S}$ the set of optimal solutions, and $\al \in (0,1)$. Assume that the following conditions hold. 
\begin{itemize}
%    \item[(i)] The set $\Tilde{A}$ does not contain any unconstrained global minimizer of $f$ (if any). Namely, $\Tilde{A}\cap (\partial f)^{-1}(0)=\emptyset$. 
%\item[(i)] The set $\bigcup_{x\in \Tilde{A}}\partial f(x)\cap S$ is nonempty.   
\item[(i)] 
%For any $x^*\in D:=\bigcup_{\substack{x\in \Tilde{A}}}\partial f(x)$, the set $A$ is $(2\al)$-sharp w.r.t. $-x^*/\norm{x^*}$ for a fixed constant $\al\in (0,1)$. 
%The following inequality holds
\begin{equation}
\label{T5.51}
\inf_{\substack{x\in A\setminus \mathbb{S}\\ x^*\in \partial f(y),y\in \mathbb{S}}} d\left(-x^*/\norm{x^*},N_A(x)\right) \ge \al >0,
\end{equation}
with the convention $\inf\emptyset = +\infty$ and $\frac{0}{\norm{0}} = 0$.
    \item[(ii)] There is $v\in \mathcal H$ such that
\begin{enumerate}
    \item[(a)] $f(v) < \inf_{x\in A}f(x)$; and
    \item[(b)] $(1- (\al/2)^2)d(v,A) < \inf_{x\in A}f(x) - f(v)$.
\end{enumerate}
\end{itemize}
In this situation, the projection of $v$ onto $A$ solves problem \eqref{P2}.
\end{theorem}

\begin{proof} Our assumption on problem \eqref{P2} implies that, at a solution $\bx \in \mathbb{S}$ there exists $x^*\in \partial f(\bx)$ such that $-x^* \in N_A(\bx)$. By (i) we know that $x^*\not=0$, otherwise the left hand side of \eqref{T5.51} equals $0$. 
    %By (ii), $A$ is  $(2\al)$-sharp w.r.t. vector $-x^*/\norm{x^*}$ for some $\al\in (0,1)$. 
    %Then,
    %$$
   % \inf_{x\in A\setminus  F_A(x^*)}d\left(\frac{x^*}{\norm{x^*}},-N_A(x)\right) \ge 2\al.
   % $$
   % We first show that $\Tilde{A}\subset F_A(x^*)$. Suppose there exists $x\in \Tilde{A}\setminus F_A(x^*)$.
  %  The latter can be re-written as
  %  \begin{equation}\label{T4.P1}
  %   \inf_{x\in A, -x^*\notin N_A(x)}d\left(\frac{x^*}{\norm{x^*}},-N_A(x)\right) \ge 2\al.
  %  \end{equation}
    % {\color{magenta}RSB 15/02: Why is the above inf over $A\setminus \tilde{A}$? We just know that $-x^*\in N_A(\bx)$.}
    % $\bigcap_{\substack{x\in \Tilde{A}}}\partial f(x)\setminus \{0\}\neq \emptyset$.
    % Take $x^*\in \bigcap_{\substack{x\in \Tilde{A}}}\partial f(x)\setminus \{0\}$, and $\bx \in \Tilde{A}$ such that $x^*\in \partial f(\bx)$. By the $\al$-sharpness condition of $A$ w.r.t vector $-x^*/\norm{x^*}$, we have
    Let $M :=\inf_A f$, we consider the sublevel set of $f$ at value $M$
    $$
    F:=\left\{x\in \mathcal H\,:\; f(x)\le M \right\}.
    $$
    Then it is clear that $A\cap F = \mathbb{S}$, the set of optimal solutions of problem \eqref{P2}.  
    % Since $\bx \in \Tilde{A}$ and $x^*\in \partial f(\bx)$, we must have
    % \[
    % 0\ge f(x)-f(\bx) \ge \ang{x^*,x-\bx},
    % \]
   % so $x^*\in (-N_A(\bx))$. 
From $x^*\in (-N_A(\bx))\cap \partial f(\bx)$, we have
 \begin{equation}\label{eq:n2}
    \ang{x^*,x-\bx} \ge 0,\;\forall x\in A\AND 0\ge f(x')-M=f(x')-f(\bx) \ge \ang{x^*,x'-\bx},\; \forall x'\in F,
     \end{equation}
    where the first inequality in the rightmost expression follows because $x'\in F$ and $\bx \in \mathbb{S}$. In other words, $x^*$ separates the (closed and convex) sets $A$ and $F$. Therefore, $A\cap F = \mathbb{S} \subset \left\{x:\; \ang{x^*,x} = \ang{x^*,\bx} \right\}$. In particular, this implies that
    \begin{equation}\label{eq:n1}
     \ang{x^*,z-\bx}=0,
    \end{equation}
    for all $z\in \mathbb{S}$.  
    Take now any $y\in A$,  $z\in \mathbb{S}$. Use the left hand side of \eqref{eq:n2} and \eqref{eq:n1} to deduce that
    \[
    \ang{-x^*,y-z}=  \ang{-x^*,y-\bx}+\ang{-x^*,\bx-z} \le 0,
    \]  
and so
$$
-x^*\in N_A(z),\quad\forall z\in \mathbb{S}.
$$
The inclusion above implies that
$
\ang{x^*,z}\le \ang{x^*,y}$,
for any $y\in A$,  $z\in \mathbb{S}$.  In other words, $\mathbb{S}\subset \argmin_A \ang{x^*,\cdot}$.  Namely,
 \[
 \begin{array}{rcl}
      \mathbb{S}&\subset&\{z\in A\::\: \ang{x^*,z}\le \ang{x^*,y} \forall\, y\in A\} 
      = \{z\in A\::\: \ang{-x^*,z}\ge \ang{-x^*,y} \forall\, y\in A\}
      = F_A(-x^*),
 \end{array}
 \]
 and therefore $A\setminus \mathbb{S}\supset A\setminus F_A(-x^*)$. Combining the latter inclusion with inequality \eqref{T5.51} yields
     \begin{equation}\label{T4.P1}
     \inf_{x\in A\setminus F_A(-x^*)}d\left(-x^*/\norm{x^*},N_A(x)\right) \ge \al.
    \end{equation}
The expression above means that $A$ is $\al$-sharp w.r.t. vector $-x^*/\norm{x^*}$.
 %Denote the latter set by $S$.%together with \eqref{T4.P1} imply 
% $$
% \inf_{x\in A, -x^*\notin N_A(x)} d\left(\frac{x^*}{\norm{x^*}},-N_A(x)\right) \ge 2\al>0.
% $
 Take $v\notin A$ such that $f(v) < M$ and $d(v,A) < \frac{M - f(v)}{1-(\al/2)^2}$. Namely, $v$ verifies condition (ii).
 %By construction, $v$ verifies condition 2 in Theorem Theorem~\ref{T3}. We claim that $b$ also verifies condition 1 in Theorem Theorem~\ref{T3}. Indeed, o
By construction, $v\in F$, and using this fact in the rightmost side of \eqref{eq:n2} yields
\[
\ang{x^*,v}\le \ang{x^*,\bx}+(f(v)-M) < \ang{x^*,\bx}\le \inf_{x\in A}\ang{x^*,x}.
\]
Hence $\ang{x^*,v}< \inf_{x\in A}\ang{x^*,x}$ and so $v$ verifies condition 1 in Theorem~\ref{T3}. We also have that $x^*\in \partial f(\bx)$, so
\begin{equation}
    \label{eq:subdiff ineq}
    \ang{x^*,\bx} -\ang{x^*,v} \ge f(\bx)- f(v) = M - f(v)>0,
\end{equation}
where we also used the fact that $f(\bar x)=M$. By definition of $v$, we have
$$
d(v,A) < \frac{M - f(v)}{1-(\al/2)^2} \le \frac{\ang{x^*,\bx} -\ang{x^*,v}}{1-(\al/2)^2}= \frac{\inf_{x\in A}\ang{x^*,x} -\ang{x^*,v}}{1-(\al/2)^2},
$$
where we used \eqref{eq:subdiff ineq} in the second inequality and the fact that $\bx\in \mathbb{S}\subset \argmin_A \ang{x^*,\cdot}$ in the equality. The expression above implies that $v$ verifies condition 2 in Theorem~\ref{T3}. Since $v$ satisfies both conditions in the theorem, we deduce that $P_A(v)$ is a solution of $\inf_{x\in A}\ang{x^*,x}$. Equivalently, $P_A(v) \in F_A(-x^*)$. To complete the proof, 
we will show that $\mathbb{S}=F_A(-x^*)$. This will establish that $P_A(v)\in \mathbb{S}$, as wanted. We already know that $\mathbb{S}\subset F_A(-x^*)$, so it is enough to show that $\mathbb{S}\supset F_A(-x^*)$. Indeed, assume that there is $z\in A\setminus \mathbb{S}$ such that $z\in F_A(-x^*)$. The fact that $z\in F_A(-x^*)$  means that $\ang{x^*,z}=\min_{x\in A}\ang{x^*,x}$, or equivalently, $-x^*\in N_A(z)$. The latter inclusion gives
\[
0=d(-x^*/\norm{x^*}, N_A(z))\ge 
\inf_{\substack{x\in A\setminus \mathbb{S}\\ x^*\in \partial f(y),y\in \mathbb{S}}} d\left(-x^*/\norm{x^*},N_A(x)\right)\ge \al>0,
\]
where we are using \eqref{T5.51} and the fact that $z\in A\setminus \mathbb{S}$ and $x^*\in \partial f(\bar x)$ with $\bar x\in \mathbb{S}$ in the first inequality. The above expression entails a contradiction and hence we must have $\mathbb{S}=F_A(-x^*)$. Therefore, our claim is true and  $P_A(v) \in F_A(-x^*)=\mathbb{S}$ and thus $P_A(v)$ solves problem \eqref{P2}.
\end{proof}

\section{Conclusions and open questions}
\label{sec:conclusion}
In this work, we introduce the notation of a sharp set and use it to analyze the single projection procedure for solving convex optimimization problem. To conclude, we outline the directions for future work.
\begin{enumerate}
  %  \item Finite termination of projection-type methods (e.g., alternating projections method, Douglas Rashford) often relies on ``one-step convergence'' of the final iteration. Namely, the projection of the last iteration must be in the solution set. 
%    \textcolor{purple}{RSB 28 Sept: The next sentence is not clear.} Thus, if one of the sets is sharp, and its arrangement with respect to other set is such that condition 2 in Theorem~\ref{T3} to hold. We conjecture that the method will stop if, e.g., $(1-\al^2)d(x_n,A) < \inf_{x\in A}\ang{x^*,x-x_n}$, where $A$ is one of the sets, $x^*$ is in the normal cone of some points in the solution set, and $x_n$ is the last iteration before convergence. Theorem~\ref{T3} may also be useful in providing an estimation on the number of iterations required for convergence, extending some convergence results in \cite{BehBellSan}.
    \item 
    The paper \cite{BehBellSan} proved the finite convergence of projection-type methods (e.g., alternating projections method, Douglas Rashford) between a closed half space and a polyhedral set for the cases when the two sets do not intersect. However, the authors in \cite{BehBellSan} do not provide an estimation of how many steps are required for the convergence. Can Theorem~\ref{T3} be used to estimate the number of steps required for convergence of the projection-type algorithms analyzed in \cite{BehBellSan}?

    \item The same results in Theorem~\ref{T3} still hold for small perturbations of the linear function, namely $\hat x^{\star} \in \mathcal{H}$ with  $\|\hat x^{\star} - x^*\|$ small enough so that $\hat x^\star \in \Int \bigcup_{x\in F_A(x^*)}N_A(x)$. This allows inexact projections, reducing the computational effort. This observation is not trivial for nonlinear functions. However, it is worth to notice that SPP may ensure finite termination for nonlinear problems with inexact gradient oracles, provided that we consider an oracle with small error and implement projections with sufficiently high accuracy as a finite operation. The precise conditions for ensuring these properties are the topic of future research.
    \item Bundle methods \cite{Oliveira_Sagastizabal_Lemarechal_2013} are designed to minimize nonsmooth convex functions. These methods approximate the original function by a suitable piece-wise function, and the iterates are minimizers of these approximations. Since these approximations have a polyhedral epigraph, they will always be sharp sets (i.e., condition (i) in Theorem~\ref{T6} will always hold for some $\alpha>0$). So, if the modulus of sharpness of the current approximation is known, conditions (a)(b) in  Theorem~\ref{T6} could potentially provide a simple way of computing the iterates. 
   
\item From Proposition~\ref{PP01}, for any polyhedral set, there is $\al>0$ such that it is $\al$-sharp w.r.t. every unit vector. It is interesting to ask whether the converse statement is true, namely: If a set is $\al$-sharp w.r.t. every unit vector for some $\al>0$, then is it true that the set must be polyhedral?
\end{enumerate}
\section*{Acknowledgements}
The research of EAN was supported by the Ministry of Science and Higher Education of the Russian Federation, project 075-02-2022-880 from 31.01.2022. HTB is supported by the Australian Research Council through the Centre for Transforming Maintenance through Data Science (grant number IC180100030). MKT would like to thank HTB for her hospitality during his stay at Curtin University supported by a Small Grant from Faculty of Science and Engineering, Curtin University.

%\bibliographystyle{plain}
%bibliography{Hoa}
\end{document}